\documentclass[11pt]{article}

\usepackage{amsmath, amsfonts, amssymb,latexsym,wasysym,graphicx}
\input epsf.sty

\newtheorem{Theorem}{Theorem}[section]
\newtheorem{Lemma}{Lemma}[section]
\newtheorem{Proposition}[Lemma]{Proposition}
\newtheorem{Corollary}[Lemma]{Corollary}
\newtheorem{Definition}[Lemma]{Definition}
\newtheorem{Remark}[Lemma]{Remark}

\newcommand{\BEQ}{\begin{equation}}     
\newcommand{\BEA}{\begin{eqnarray}}
\newcommand{\BD}{\begin{displaymath}}
\newcommand{\EEQ}{\end{equation}}       
\newcommand{\EEA}{\end{eqnarray}}
\newcommand{\ED}{\end{displaymath}}

\newcommand{\del}{\delta}

\newcommand{\eps}{\varepsilon}          

\newcommand{\Arg}{{\mathrm{Arg}}}

\newcommand{\R}{\mathbb{R}}
\newcommand{\C}{\mathbb{C}}
\newcommand{\Z}{\mathbb{Z}}
\newcommand{\N}{\mathbb{N}}

\def\proba{{\mathbb{P}}}
\def\esper{{\mathbb{E}}}

\newcommand{\eop}{\hfill $\Box$}        

\newcommand{\II}{{\rm i}}               
\renewcommand{\Re}{{\rm Re\ }}          
\renewcommand{\Im}{{\rm Im\ }}          
\newcommand{\half}{{1\over 2}}          
\newcommand{\wit}[1]{\widetilde{#1}}    
\renewcommand{\vec}[1]{\boldsymbol{#1}} 

\catcode`\@=11
\def\numberbysection{\@addtoreset{equation}{section}
        \def\theequation{\thesection.\arabic{equation}}}
\numberbysection

\begin{document}

\vspace*{1.5cm}
\begin{center}
{\Large \bf A central limit theorem for the rescaled L\'evy area
of two-dimensional fractional Brownian motion with Hurst index
$H<1/4$}
\end{center}

\vspace{2mm}
\begin{center}
{\bf  J\'er\'emie Unterberger}
\end{center}

\vspace{2mm}
\begin{quote}

\renewcommand{\baselinestretch}{1.0}
\footnotesize
{Let $B=(B^{(1)},B^{(2)})$ be a two-dimensional fractional Brownian motion with Hurst index $\alpha\in (0,1/4)$.
Using an analytic approximation $B(\eta)$ of $B$ introduced in \cite{Unt08}, we  prove
 that the rescaled L\'evy area process $(s,t)\to
\eta^{\half(1-4\alpha)}\int_s^t dB_{t_1}^{(1)}(\eta) \int_s^{t_1}
dB_{t_2}^{(2)}(\eta)$ converges in law to $W_t-W_s$ where $W$ is a  Brownian motion independent from $B$.
 The method relies on 
a very general scheme of analysis of singularities of analytic functions,  applied to  the moments of 
 finite-dimensional distributions of the L\'evy area.
 }
\end{quote}

\vspace{4mm}


\newpage


\section{Introduction}


The (two-sided) fractional Brownian motion $t\to B_t$, $t\in\R$  (fBm for short) with Hurst exponent $\alpha$, $\alpha\in(0,1)$, defined as the centered Gaussian process with covariance
\BEQ \esper[B_s B_t]=\half (|s|^{2\alpha}+|t|^{2\alpha}-|t-s|^{2\alpha}), \EEQ
is a natural generalization in the class of Gaussian processes of
the usual Brownian motion, in the sense that it exhibits two fundamental properties shared with Brownian motion, namely,
it has stationary increments, viz. $\esper[(B_t-B_s)(B_u-B_v)]=\esper[(B_{t+a}-B_{s+a})(B_{u+a}-B_{v+a})]$ for
every $a,s,t,u,v\in\R$, and it is self-similar, viz. 
\BEQ 
\forall \lambda>0, \quad (B_{\lambda t}, t\in\R) \overset{(law)}{=} (\lambda^{\alpha} B_t,
t\in\R).
\EEQ 
One may also define  a $d$-dimensional 
vector Gaussian process (called: {\it $d$-dimensional fractional Brownian motion}) by setting $B_t=(B_t(1),\ldots,B_t(d))$ where $(B_t^{(i)},t\in\R)_{i=1,\ldots,d}$ are $d$ independent (scalar) fractional Brownian motions.

Its theoretical interest lies in particular in the fact that it is (up to normalization) the only Gaussian process satisfying these two properties.

A standard application of Kolmogorov's theorem shows that fBm has a version with $(\alpha-\eps)$-H\"older paths 
for every $\eps>0$. In particular, all its paths possess finite $q$-variation for every $q>\frac{1}{\alpha}$, in the sense that
\BEQ \sup_{n\ge 1} \sup_{s=t_0<\ldots<t_n=t} \left( \sum_{l=0}^n |B_{t_l}-B_{t_{l-1}}|^q \right)<\infty 
\quad {\mathrm{a.s.}}  \EEQ
where the sum ranges over all partitions $(s=t_0<t_1<\ldots<t_n=t)$ of any order $n$ of the interval $[s,t]$. 

There has been a widespread interest during the past ten years in constructing a stochastic integration theory
with respect to fBm and solving stochastic differential equations driven by fBm, see for instance
\cite{LLQ02,GraNouRusVal04,CheNua05,RusVal93,RusVal00}. The multi-dimensional case is
very different from the one-dimensional case. When one tries to integrate for instance a stochastic differential
equation driven by a two-dimensional fBm $B=(B^{(1)},B^{(2)})$ by using any kind of Picard iteration scheme, one
encounters very soon the problem of defining the L\'evy area of $B$ which is the antisymmetric part
of ${\cal A}_{s,t}:=\int_s^t dB_{t_1}^{(1)} \int_s^{t_1} dB_{t_2}^{(2)}$. This is the simplest occurrence
of iterated integrals $\vec{B}^k_{s,t}(i_1,\ldots,i_k):=\int_s^t dB_{t_1}^{(i_1)}\ldots \int_s^{t_{k-1}} dB_{t_k}^{(i_k)}$,
$i_1,\ldots,i_k\le d$
for $d$-dimensional fBm $B=(B^{(1)},\ldots,B^{(d)})$ 
which lie at the heart of the rough path theory due to T. Lyons, see \cite{Lyo98,LyoQia02}. 
Let us describe this briefly. Assume
$\Gamma_t=(\Gamma_t^{(1)},\ldots,\Gamma_t^{(d)})$ is some non-smooth $d$-dimensional path with bounded $q$-variation
for some $q>1$ (take for instance an $\alpha$-H\"older path with $\alpha=1/q<1$). Integrals such as
$\int f_1(\Gamma_t)d\Gamma^{(1)}_t+\ldots+f_d(\Gamma_t)d\Gamma^{(d)}_t$ do not make sense a priori because
$\Gamma$ is not differentiable (Young's integral \cite{Lej03} works for $\alpha>\half$ but not beyond). In order to define
the integration of a differential form along $\Gamma$, it is enough to define a truncated {\it multiplicative functional} (consisting of increments)
$(\vec{\Gamma}^{1},\ldots,\vec{\Gamma}^{\lceil q\rceil})$, $\lceil  q\rceil$=entire part of $q$,
  where $\vec{\Gamma}^{1}_{s,t}=\Gamma_t-\Gamma_s$ and
each
$\vec{\Gamma}^k=(\vec{\Gamma}^k(i_1,\ldots,i_k))_{1\le i_1,\ldots,i_k\le d}$, $k\ge 2$ -- a matrix of continuous paths --
is a {\it substitute} for the iterated integrals $\int_s^t d\Gamma_{t_1}^{(i_1)}\int_s^{t_1} d\Gamma_{t_2}^{(i_2)}
\ldots \int_{s}^{t_{k-1}} d\Gamma_{t_k}^{(i_k)}$ with the following two properties:

\begin{itemize}
\item[(i)] each component of $\vec{\Gamma}^k$ has bounded $\frac{q}{k}$- variation semi-norm
$$||\vec{\Gamma}^k_{s,t}||^{q/k}_{q/k}:=
\sup_{n\ge 1} \sup_{s=t_0<t_1<\ldots<t_n=t} \sum_{l=1}^{n}
 |\vec{\Gamma}^k_{t_{l-1},t_{l}}|^{q/k};$$ 
\item[(ii)] ({\it multiplicativity}) letting $\vec{\Gamma}^k_{s,t}:=\vec{\Gamma}_t^k-\vec{\Gamma}_s^k$, one requires
\BEA
 \vec{\Gamma}^k_{s,t}(i_1,\ldots,i_k) &=& \vec{\Gamma}^k_{s,u}(i_1,\ldots,i_k)+\vec{\Gamma}^k_{u,t}(i_1,\ldots,i_k) \nonumber\\
&+& \sum_{k_1+k_2=k} \vec{\Gamma}_{s,u}^{k_1}(i_1,\ldots,i_{k_1}) \vec{\Gamma}_{u,t}^{k_2}(i_{k_1+1},\ldots,i_k).
\nonumber\\ \EEA
\end{itemize}

The multiplicativity property implies in particular the following identity for the L\'evy area:
\BEQ {\cal A}_{s,t}={\cal A}_{s,u}+{\cal A}_{u,t}+(B^{(1)}_u-B^{(1)}_s)(B^{(2)}_t-B^{(2)}_u). \label{eq:0:mult} \EEQ

Letting $\vec{\Gamma}=1\oplus\vec{\Gamma}^1 \oplus \vec{\Gamma}^2\oplus\ldots\oplus \vec{\Gamma}^{\lceil q \rceil}$
live in the truncated tensor algebra $\R\oplus \R^d\oplus (\R^d)^{\otimes 2}\oplus \ldots \oplus (\R^d)^{\otimes \lceil q \rceil}$,
the latter property reads simply $\vec{\Gamma}_{s,t}=\vec{\Gamma}_{s,u}\otimes \vec{\Gamma}_{u,t}$. Then there
is a standard procedure which allows to define out of these data iterated integrals of any order and
to solve differential equations driven by $\Gamma$.

The multiplicativity property is satisfied by smooth paths, as can be checked by direct computation. So the most
natural way to construct such a multiplicative functional is to start from some smooth approximation
$\Gamma(\eta)$, $\eta\overset{>}{\to} 0$ of $\Gamma$ such that each iterated
 integral $\vec{\Gamma}^k_{s,t}(\eta)(i_1,\ldots,i_k)$, $k\le \lceil q \rceil$ converges in the $\frac{q}{k}$ -variation 
semi-norm.

This general scheme has been applied to fBm in a paper by L. Coutin and Z. Qian \cite{CQ02} and later in a paper
by the author \cite{Unt08}. L. Coutin and Z. Qian used the standard $n$-dyadic piecewise linear approximation
$B^{CQ}(2^{-n})$ of $B$. Our approximation consists in seeing $B$ as the real part of the boundary value of 
an analytic process  $\Gamma$ 
living on the upper half-plane $\Pi^+=\{z\in\C\ |\  \Im z>0\}$. The time-derivative of this
centered Gaussian process has the following hermitian positive-definite covariance kernel:
\BEQ \esper[\Gamma'(z) \overline{\Gamma'(w)}]=: K^{',-}(z,\bar{w})=\frac{\alpha(1-2\alpha)}{2\cos\pi\alpha}
(-\II(z-\bar{w}))^{2\alpha-2},\quad z,w\in\Pi^+,\EEQ
where $z^{2\alpha-2}:=e^{(2\alpha-2)\ln z}$ (with the usual determination of the logarithm) is defined and analytic
on the cut plane $\C\setminus\R_-$. Also, by construction, 
\BEQ \esper[\Gamma'(z)\Gamma'(w)]\equiv 0 \label{eq:0:0} \EEQ  identically.
It is essential to understand that $K'$ is a singular multivalued function on $\C\times\C$; for $z,w\in\Pi^+$,
$\Re (-\II(z-\bar{w}))>0$ so the kernel $K'$ is well-defined. 
Then 
\BEQ B_t(\eta):=\Gamma_{t+\II\frac{\eta}{2}}+\overline{\Gamma_{t+\II\frac{\eta}{2}}}=2\Re \Gamma_{t+\II\frac{\eta}{2}}  \label{eq:0:BXi} \EEQ
 is a {\it good} approximation of fBm, namely, $B(\eta)$ converges a.s. in
the $q$-variation distance for every  $q>\frac{1}{\alpha}$ to a process $B_t(0)=2\Re \Gamma_t$ with the same law as fBm.

Both approximation schemes lead to the same semi-quantitative result, namely:

-- when $\alpha>1/4$, the L\'evy area and volume (in other words, the  multiplicative functional truncated
to order $3$) converge a.s. in the correct variation norm. The heart of the proof lies in the study of the
piecewise linear approximation  ${\cal A}^{CQ}_{s,t}(2^{-n})$, resp. analytic
approximation ${\cal A}_{s,t}(\eta)$ of the L\'evy area; one may prove in
particular that $\esper[\left({\cal A}^{CQ}_{s,t}(2^{-n})\right)^2]$ and $\esper[\left({\cal A}_{s,t}(\eta)\right)^2]$
converge to the same limit when $2^{-n}$ and $\eta$ go to 0;

-- when $\alpha<1/4$, $\esper[\left({\cal A}^{CQ}_{s,t}(2^{-n})\right)^2]$ and $\esper[\left({\cal A}_{s,t}(\eta)\right)^2]$
diverge resp. like $n^{(1-4\alpha)}$ and $\eta^{-(1-4\alpha)}$. Hence the above method fails.

The latter result is of course unsatisfactory, and constitutes by no means a proof that no coherent
stochastic integration theory
with respect to fBm may exist when $\alpha<1/4$.

\vskip 1 cm

We are interested in this paper in the singular case $\alpha<1/4$. We give no construction of a rough path
but provide results which, hopefully, may be a first step in that direction.

First of all (see Section 1), we give a precise analysis of the singular terms appearing in the moments of
 the L\'evy area  constructed out of the analytic approximation of  a two-dimensional fractional Brownian motion
 $B(\eta)=(B^{(1)}(\eta),B^{(2)}(\eta))$
when $\eta\to 0$.  The main results are Theorem \ref{th:2:divergence-moments} and Corollary \ref{cor:2:2Nth-moment}
 which
give in particular an equivalent of $\esper[( {\cal A}_{s,t}(\eta))^{2N}]$ when $\eta\to 0$. They are in fact
much more precise in that they provide a general method to find 
 the exponents $(4\alpha-1)N=\beta_0<\beta_1<\beta_2<\ldots$ of the asymptotic
 expansion of the moments, namely,
\BEQ \esper[( {\cal A}_{s,t}(\eta))^{2N}]=\sum_{j=0}^J c_j |t-s|^{4\alpha N-\beta_j} \eta^{\beta_j}+
o(\eta^{\beta_J}) \EEQ
where the $c_j$ are coefficients which depend only on $\alpha$, and $c_0$ may be evaluated explicitly. 
 But an easy generalization
 yields the same kind of results for $\esper[{\cal A}_{s_1,t_1}(\eta)\ldots {\cal A}_{s_{2N},t_{2N}}(\eta)]$ where
 $s_1<t_1,\ldots,s_{2N}<t_{2N}$  are arbitrary arguments (see Lemma \ref{lemma:3:convergence-moments2} in
 Section 2). The method we use is sufficiently general to be applied with some extra efforts to any kind of
iterated integral of any order.

Section 1 may be seen as a long exercise in complex analysis. Let us give a simple example coming from
\cite{Unt08}. 
By definition (recall $\esper[ \Gamma'(z)\Gamma'(w)] \equiv 0$
identically, see (\ref{eq:0:0}))
\begin{eqnarray} \esper [({\cal A}_{s,t}^{\eta})^2] &=& 2\esper \left( \int_s^t d\Gamma_{x_1+\II\frac{\eta}{2}}^{(1)} \int_s^{x_1}
d\Gamma_{x_2+\II\frac{\eta}{2}}^{(2)} \right)\left( \int_s^t d\bar{\Gamma}_{y_1+\II\frac{\eta}{2}}^{(1)}\int_s^{y_1}
d\bar{\Gamma}_{y_2+\II\frac{\eta}{2}}^{(2)} \right) \nonumber\\
&+& 2\Re \esper  \left( \int_s^t d\Gamma_{x_1+\II\frac{\eta}{2}}^{(1)} \int_s^{x_1}
d\bar{\Gamma}_{x_2+\II\frac{\eta}{2}}^{(2)} \right)\left( \int_s^t d\bar{\Gamma}_{y_1+\II\frac{\eta}{2}}^{(1)}\int_s^{y_1}
d\Gamma_{y_2+\II\frac{\eta}{2}}^{(2)} \right) \nonumber\\
&=:& {\cal V}_1(\eta)+{\cal V}_2(\eta).
\end{eqnarray}

The first term in the right-hand side writes (using the stationarity of the increments)

\begin{eqnarray}
&& {\cal V}_1(\eta)=C \int_0^{t-s} dx_1 \int_0^{x_1} dx_2 \int_0^{t-s} dy_1 \int_0^{y_1} dy_2
(-\II(x_1-y_1)+\eta)^{2\alpha-2} (-\II(x_2-y_2)+\eta)^{2\alpha-2} \nonumber\\
&& = C' \int_0^{t-s} dx_1 \int_0^{t-s} dy_1 (-\II(x_1-y_1)+\eta)^{2\alpha-2} \nonumber\\
&& \quad  
 \left[(-\II(x_1-y_1)+\eta)^{2\alpha}-
(-\II x_1+\eta)^{2\alpha}-(\II y_1+\eta)^{2\alpha} \right] \nonumber\\
\end{eqnarray}

while the second term writes
\begin{eqnarray}
&&{\cal V}_2(\eta)=\nonumber\\
&& C' \int_0^{t-s} dx_1 \int_0^{t-s} dy_1 (-\II(x_1-y_1)+\eta)^{2\alpha-2} 
 \left[(\II(x_1-y_1)+\eta)^{2\alpha}-
(\II x_1+\eta)^{2\alpha}-(-\II y_1+\eta)^{2\alpha} \right] \nonumber\\
\end{eqnarray}

Both integrals look the same {\em except} that ${\cal V}_2$ (contrary to ${\cal V}_1$) involves
both $-\II x_1$ and $\II x_1$, and similarly for $y_1$. This seemingly insignificant difference
is essential, since ${\cal V}_1$ can be shown to have a bounded limit when $\eta\to 0$ by using a 
contour deformation in $\Pi^+\times\Pi^-$ (where $\Pi^-$ denotes the {\em lower} half-plane)
 which avoids the real axis where singularities live,
while this is impossible for ${\cal V}_2$ (namely, $(-\II(x_1-y_1)+\eta)^{2\alpha-2}$ is well-defined
if $(x_1,y_1)$ are in the closure of $\Pi^+\times\Pi^-$, while $(\II(x_1-y_1)+2\eta)^{2\alpha}$
for instance is well-defined on the closure of $\Pi^-\times\Pi^+$).
  In fact, explicit computations using Gauss' hypergeometric function
 prove that ${\cal V}_2$ diverges in the limit
$\eta\to 0$ when $\alpha<1/4$. More general results are given in Lemmas \ref{lemma:2:lemma3} and \ref{lemma:2:lemma3bis}. Let us state a simple consequence of them. Let
\BEQ I_-(\beta_1,\beta_2;0,t)(a,b):=\int_0^t (-\II(u-a)+\eta)^{\beta_1} (-\II(u-b)+\eta)^{\beta_2}\ du 
\label{eq:0:I-} \EEQ
and
\BEQ  I_+(\beta_1,\beta_2;0,t)(a,b):=\int_0^t (+\II(u-a)+\eta)^{\beta_1} (-\II(u-b)+\eta)^{\beta_2}\ du 
\label{eq:0:I+} \EEQ
for $\beta_1,\beta_2\in\R$ such that $\beta_2>-1,\beta_1+\beta_2+1<0$ and $a,b\in(0,t)$. (Notice
one retrieves terms contained in ${\cal V}_1(\eta)$ or ${\cal V}_2(\eta)$ when one sets $a=b$). Take $\eta\to 0$ and
$a-b\to 0$.
  Then the integral $I_-$ converges (which follows again from a deformation of contour), while $I_+$ diverges
like $C(\II(b-a)+2\eta)^{\beta_1+\beta_2+1}$.

When evaluating the $2N$-th moment of ${\cal A}_{s,t}(\eta)$ for instance, iterated integrals of the same type
as $I_+$ produce multivalued power functions $(\pm\II(x-y)+k\eta)^{\beta}$ with $k\in\N$ and various exponents $\beta=
4\alpha-1,6\alpha,8\alpha-1,10\alpha,\ldots$. The
singularities of the moments come exactly from those non-analytic terms with {\em negative} exponent when evaluated on the 
diagonal $x=y$. Aside from these terms, iterated
integrals also produce {\em analytic terms} of the form

$$F(z):=\int_a^b (-\II(z-u)+\eta)^{\gamma} (u-a)^{\beta} f(u)\ du$$
($\gamma=2\alpha$ or $2\alpha-2$, $\eta>0$, $\beta>-1$) where $f$ is analytic on some appropriate complex
 domain containing $[a,b]$. The
question is to give a maximal domain where $F$ is analytic and to understand its singularities around $a$. Although
one is mainly interested in the behaviour of $F$ on the real axis, computations show clearly that
complex analytic methods (including contour shifts e.g.) are the most appropriate in this setting, and circumvent
the heavy arguments one would inevitably get using a pathwise linear approximation.
 A detailed analysis shows that all these terms are regular in the limit $\eta\to 0$. A shorcut
(which does not give the values of the exponents though) is provided in Lemma \ref{lemma:2:refinement}.

The reader who is interested more in applications than in the details should essentially have in mind 
the two following results:

{\bf Lemma \ref{lemma:1:exp}}

{\it
The generating function $\phi_{s,t}(\eta;\lambda):=
\esper[ e^{\II\lambda {\cal A}_{s,t}(\eta)}]$ of the L\'evy area
 ${\cal A}_{s,t}(\eta)$  is the exponential of the generating function of
connected diagrams, i.e.
\BEQ \phi_{s,t}(\eta;\lambda)=\exp \phi_{s,t}^{(c)}(\eta;\lambda).\EEQ
}

(see subsection 1.1 for the definition of connected moments which are Gaussian cumulants of a certain type);

\medskip

{\bf Theorem \ref{th:2:divergence-moments}}

{\it
The $2N$-th connected moment of  ${\cal A}_{s,t}(\eta)$ is given
by the sum of two terms: the first one is regular in the limit $\eta\to 0$ and equal to $C_{reg,N} t^{4N\alpha}+O(\eta^{2\alpha})$ for
some constant $C_{reg,N}$; the second one is equal to $C_{irr,N} t \eta^{4N\alpha-1}$ with
\BEQ C_{irr,N}=\left( \frac{\pi/2}{\cos \pi\alpha \Gamma(-2\alpha)} \right)^{2(N-1)} \sin \pi\alpha \frac{\Gamma(2\alpha+1)}{\Gamma(2-2\alpha)} \Gamma(1-4\alpha N) (2N)^{4\alpha N-1}.\EEQ
}

\bigskip

\bigskip

In Section 2, we apply this analysis of singularities to convergence results concerning the L\'evy area. Since
the second moment of the  L\'evy area ${\cal A}_{s,t}(\eta)$
diverges like $\eta^{-(1-4\alpha)}$, it is natural to introduce the
{\em rescaled L\'evy area} 
 \BEQ \tilde{\cal A}_{s,t}(\eta)=\eta^{\half(1-4\alpha)}{\cal A}_{s,t}(\eta).\EEQ
The main result of this paper is a kind of central limit theorem   which we state here:

\newpage

{\bf Theorem A.}

{\it The three-dimensional  process $(B^{(1)}(\eta),B^{(2)}(\eta),\tilde{\cal A}(\eta))$  converges in law to
 $(B^{(1)},B^{(2)},\sqrt{C_{irr,1}} \del W)$
where $\del W_{s,t}:=W_t-W_s$ are the increments of  a standard one-dimensional Brownian motion {\em
independent} from $B^{(1)}$ and $B^{(2)}$.
}

The value of the constant $C_{irr,1}=\lim_{\eta\to 0} \esper[ (\tilde{\cal A}_{0,1}(\eta))^2]$ is given
in Theorem \ref{th:2:divergence-moments}.  We feel the
exact value is not important though (different schemes of approximation lead to different values, whereas the value
of $\lim_{\eta\to 0}\esper[ (\tilde{\cal A}_{0,1}(\eta))^2]$ when $\alpha>1/4$ seems to be more universal,
see above).

Other central limit theorems have been obtained for sums of integer powers of increments 
of one-dimensional fBm (see \cite{BreMaj83,NuaOrt08,NouPec08} for instance). The most closely related result is maybe 
that of I. Nourdin \cite{Nou08} which shows
an It\^o-type formula for a two-dimensional fBm with Hurst index $H=1/4$ with a 'bracket-term' involving
an independent Brownian motion, but the results are of a very different nature than ours (also, they hold precisely
for $H=1/4$, whereas our results concern the case $H<1/4$). Note also the paper by Y. Hu and D. Nualart \cite{HuNu05} which
gives a Brownian scaling limit for the self-intersection local time for $\alpha$ large enough.

One would then like to say that the L\'evy area ${\cal A}_{s,t}(\eta)$ writes (up to some finite coefficient)
as a counterterm $\eta^{-\half(1-4\alpha)}(W_t-W_s)$ where $W$ is a standard Brownian motion, plus some
term which is finite in the limit $\eta\to 0$. Unfortunately this statement is not true with the Brownian motion $W$
constructed in Theorem A because it is independent from the L\'evy area ${\cal A}_{s,t}^{\eta}$. But one may
conjecture that some related counterterm
yields a corrected L\'evy area which is finite.

Let us also mention the last result of this paper, which yields a uniform exponential bound
 for the rescaled L\'evy area  when $\alpha\in(\frac{1}{8},\frac{1}{4})$ (see
Corollary \ref{cor:2:uem}). We give no application, but it allows for instance to use the Berry-Ess\'een Lemma to get
 precise estimates of
the rate of convergence of $\tilde{\cal A}_{s,t}$ to $\sqrt{C_{irr,1}}\del W_{s,t}$. Also, Stein's method combined
with the Malliavin calculus \cite{NouPec08},
\cite{NouPec08bis}  may
easily be applied to our setting to give convergence rates since Section 1 gives  estimates of all
 cumulants of the L\'evy area. 
We hope to come back to this in the future.

\bigskip

We shall be using a number of times the following integral representation
of Gauss' hypergeometric function $_2 F_1$ (see \cite{Abra84}, 15.3.1):

\BEQ
_2 F_1(a,b,c;z)=\frac{\Gamma(c)}{\Gamma(b)\Gamma(c-b)} \int_0^1 t^{b-1}
(1-t)^{c-b-1} (1-tz)^{-a}\ dt, \label{hyper1}
\EEQ
valid if $\Re c>\Re b>0$ and $z\in\C\setminus[1;+\infty)$. Recall that $_2 F_1(a,b,c;z)$ $(c\not=0,-1,\ldots)$
 is defined around $z=0$ by an infinite series
with radius of convergence 1 and has an analytic extension to the cut plane $\C\setminus[1;+\infty[$.
 The {\it connection formulas} give $_2 F_1(a,b,c;z)$ in terms of a linear combination of hypergeometric
 functions in the transformed argument $\phi(z)$, where $\phi$ is any projective transformation of the Riemann sphere
 preserving the set of singularities of the hypergeometric differential equation, namely $\{0,1,\infty\}$. They relate
 the behaviour of the hypergeometric functions around $0$ with their behaviour around $1$ and $\infty$.    We  reproduce here
 for the convenience
of the reader three connection formulas, relating the behaviour around $\infty$ with the behaviour around $0$ with $0$ for
the first two ones, and the behaviour around $1$ with the behaviour around $0$ for the third one
 (see \cite{Abra84}, section 15.3):
\BEA
_2 F_1(a,b,c;z)&=& \frac{\Gamma(c)\Gamma(b-a)}{\Gamma(b)\Gamma(c-a)} (-z)^{-a}
\ _2 F_1(a,1-c+a ;1-b+a;\frac{1}{z}) \nonumber \\ &+& 
 \frac{\Gamma(c)\Gamma(a-b)}{\Gamma(a)\Gamma(c-b)} (-z)^{-b}
\ _2 F_1(b,1-c+b ;1-a+b;\frac{1}{z}),\quad z\not\in\R_+ \nonumber\\  \label{eq:0:1/z}
\EEA

\BEA
_2 F_1(a,b,c;z)&=&  (1-z)^{-a} \frac{\Gamma(c)\Gamma(b-a)}{\Gamma(b)\Gamma(c-a)} 
\ _2 F_1(a,c-b ;a-b+1;\frac{1}{1-z}) \nonumber \\ &+&  (1-z)^{-b}
 \frac{\Gamma(c)\Gamma(a-b)}{\Gamma(a)\Gamma(c-b)} 
\ _2 F_1(b,c-a ;b-a+1;\frac{1}{1-z}),\quad z\not\in[1;+\infty) \nonumber \\ \label{eq:0:1/1-z}
\EEA

\BEA
&&_2 F_1(a,b,c;z)= \frac{\Gamma(c)\Gamma(c-a-b)}{\Gamma(c-a)\Gamma(c-b)} 
\ _2 F_1(a,b ;a+b-c+1;1-z) +  (1-z)^{c-a-b} \nonumber\\
&&  \frac{\Gamma(c)\Gamma(a+b-c)}{\Gamma(a)\Gamma(b)} 
\ _2 F_1(c-a,c-b ;c-a-b+1;1-z),\quad z\not\in[1;+\infty) \nonumber \\ \label{eq:0:1-z}
\EEA

Let us also recall that $_2 F_1(a,b,c;z)$ is symmetric in the arguments $a,b$, constant (equal to 1) if $a=0$ or $b=0$, and that
\BEQ _2 F_1(a,b;a;z)=\ _2 F_1(b,a;a;z)=(1-z)^{-b}\EEQ
which is a consequence of \cite{Abra84}, 15.3.4.


\section{Moments of the L\'evy area}


Section 1 (by far the longest one of the article) is organized as follows.

Subsection 1.1 contains the main definitions, followed by a diagrammatic expansion of the moments of
the L\'evy area (see in particular Definition \ref{def:1:generating-functions} for the definition of 
the generating function of the connected diagrams).

Subsection 1.2 gives the general scheme for subsequent computations and contains explicit closed formulas
(see Lemmas \ref{lemma:2:lemma3} and \ref{lemma:2:lemma3bis}) for the functions $I_{\pm}$ defined in the
Introduction, see equations (\ref{eq:0:I-}), (\ref{eq:0:I+}). 

Subsection 1.3 (dedicated to the computations of the singularity exponents) is the heart of the section, and
(unfortunately) the most technical one. The reader who is interested more in applications than in the details of
the proofs may skip it, since the terms evaluated in this paragraph (called {\it admissible functions})
turn out in the end to be regular in the limit $\eta\to 0$.

Finally, subsection 1.4 gives the asymptotic behaviour of the connected moments  of the L\'evy area when
$\eta\to 0$, from which one deduces easily the asymptotic behaviour of $\esper[ ({\cal A}_{s,t}(\eta))^{2N}]$.


\subsection{Definitions and combinatorial arguments}


Let us start by introducing three kernels which will be the fundamental objects of study in this article.

\begin{Definition}

{\it
Let, for $\eta>0$,
\begin{enumerate}
\item 
\BEQ K^{',\pm}(\eta;x,y)=\frac{\alpha(1-2\alpha)}{2\cos\pi\alpha} (\pm\II(x-y)+\eta)^{2\alpha-2}; \EEQ
\item
 \BEA
&&  K^{\pm}(\eta;x,y)=\int_0^x du \int_0^y dv\  K^{',\pm}(\eta;u,v) \nonumber\\
&& =\frac{1}{4\cos\pi\alpha} 
 \left((\pm\II x+\eta)^{2\alpha}+(\mp\II y+\eta)^{2\alpha}-(\pm\II(x-y)+\eta)^{2\alpha}\right); \nonumber\\  \EEA

\item 
\BEQ
K^{*,\pm}(\eta;x,y)=-\frac{1}{4\cos\pi\alpha}   (\pm\II(x-y)+\eta)^{2\alpha}.
\EEQ
\end{enumerate}

Similarly, let $K'(\eta;x,y):=2\Re K^{',\pm}(\eta;x,y)$, $K(\eta;x,y):=2\Re K^{\pm}(\eta;x,y)$ and
$K^*(\eta;x,y):=2\Re K^{*,\pm}(\eta;x,y)$ be the real parts (up to a coefficient $2$) of the previous kernels.

}
\label{def:1:K'K}
\end{Definition}

As showed in \cite{Unt08}, the kernel $K'(\eta)$ is  positive and represents (for every fixed $\eta>0$)
 the covariance of 
of a real-analytic centered Gaussian process with real time-parameter $t$. The easiest way to see it is to make use of the
following explicit series expansion: letting (for $k\ge 0$)
\BEQ f_k(z)=2^{\alpha-1} \sqrt{\frac{\alpha(1-2\alpha)}{2\cos\pi\alpha}}
 \sqrt\frac{\Gamma(2-2\alpha+k)/\Gamma(2-2\alpha)}{k!}
\left( \frac{z+\II}{2\II} \right)^{2\alpha-2} \left( \frac{z-\II}{z+\II}\right)^k, \quad z\in\Pi^+ \EEQ
one has
\BEQ \sum_{k\ge 0} f_k(x_1+\II \frac{\eta_1}{2}) \overline{f_k(y+\II\frac{\eta_2}{2})} = K^{',-}(\half(\eta_1+\eta_2);x,y).\EEQ 
Define more generally  a Gaussian process with time parameter $z\in \Pi^+$ as follows:
\BEQ \Gamma'(z)=\sum_{k\ge 0} f_k(z)\xi_k \EEQ
 where $(\xi_k)_{k\ge 0}$ are independent standard complex
Gaussian variables, i.e. $\esper[\xi_j \xi_k]=0$, $\esper[\xi_j \bar{\xi}_k]=\del_{j,k}$.
The Cayley transform $\Pi^+ \to {\cal D}, z\to  \frac{z-\II}{z+\II}$ ($\cal D$=unit disk of
the complex plane) makes the series defining $\Gamma'$ into a random entire series which may
be shown to be analytic on the unit disk by standard arguments. Hence the process $\Gamma'$ is analytic on $\Pi^+$.
Note that (restricting to the horizontal line $\R+\II\frac{\eta}{2}$)
 $\Re \esper[\Gamma'(x+\II\eta/2)\overline{\Gamma'(y+\II\eta/2)}]=K^{',-}(\eta;x,y)$.
 
One may now integrate the process $\Gamma'$ over any path $\gamma:(0,1)\to\Pi^+$ with endpoints $\gamma(0)=0$
and $\gamma(1)=z\in\Pi^+\cup\R$ (the result does not depend on the particular path but only on the endpoint $z$).
The result is a process $\Gamma$ which is still analytic on $\Pi^+$. As mentioned in the Introduction,
one may retrieve the fractional Brownian motion by considering the real part of the boundary value of $\Gamma$
on $\R$. Another way to look at it is to define $\Gamma_t(\eta):=\Gamma(t+\II\eta)$
 as a regular process living on $\R$, and to remark that the real part of $\Gamma(\eta)$ converges
when $\eta\to 0$ to fBm.  In the following Proposition, we give precise statements which
summarize what has been said up to now:

\begin{Proposition}[see \cite{Unt08}]

{\it

\begin{enumerate}
\item
Let $\gamma:(0,1)\to\Pi^+$ be a continuous path with endpoints $\gamma(0)=0$ and $\gamma(1)=z$,
 and set $\Gamma_z=\int_{\gamma}\Gamma'_u \, du$. Then $\Gamma$ is an analytic process on $\Pi^+$.
 Furthermore, as $z$ runs along any path in $\Pi^+$ going to $t\in\R$, the random variables $\Gamma_z$
 converge almost surely to a random variable called again $\Gamma_t$.

\item

The family $\{\Gamma_t;\, t\in\R\}$ defines a  centered Gaussian complex-valued process whose
paths are almost surely $\kappa$-H\"older for any $\kappa<\alpha$. Its real
part $B_t:=2\Re\Gamma_t$  has the same law as fBm.

\item
The family of centered Gaussian real-valued processes $B(\eta)_t:=\Re \Gamma_{t+\II\eta}$ converges a.s. to $B_t$ in the
$q$-variation norm for very $q>\frac{1}{\alpha}$. Its covariance kernel is $K(\eta)$.

\end{enumerate}

}
\label{prop:1}
\end{Proposition}

\bigskip

Let us introduce the L\'evy area for a two-dimensional fBm $B=(B^{(1)},B^{(2)})$.

\begin{Definition}

{\it
Let
\BEQ {\cal A}_{s,t}(\eta):=\int_s^t dB_x^{(1)}(\eta) \int_s^t dB_y^{(2)}(\eta). \EEQ
}
\label{def:1:Area}
\end{Definition}

In order to evaluate the moments of the L\'evy area, we first need some combinatorial arguments. Recall
 to begin with the classical

\begin{Proposition}

{\it
Let $(X_1,\ldots,X_{2N})$ be a Gaussian vector with zero means. Then
\BEQ \esper[X_1\ldots X_{2N}]=\sum_{(i_1,i_2),\ldots,(i_{2N-1},i_{2N})}
 \prod_{j=1}^N \esper[X_{i_{2j}} X_{i_{2j+1}}] \EEQ
where the sum ranges over the $(2N-1)!!=1.3.5\cdots (2N-1)$ couplings of the indices
$1,\ldots,2N$.
}
\label{prop:1:Wick}
\end{Proposition}

\begin{Lemma}

{\it
\BEA \esper[{\cal A}_{s,t}(\eta)^{2N}]&=&  \int_0^{t-s} dx_1  \ldots
 \int_0^{t-s} dx_{2N}  \sum_{(i_1,i_2),\ldots,(i_{2N-1},i_{2N})}
\sum_{(j_1,j_2),\ldots,(j_{2N-1},j_{2N})}  \nonumber\\
&& \prod_{k=1}^N K'(\eta;x_{i_{2k-1}},x_{i_{2k}}) \ .\
\prod_{k=1}^N K(\eta;x_{j_{2k-1}},x_{j_{2k}}). \EEA
}
\label{lemma:1:Area2N}
\end{Lemma}

{\bf Proof.}

By stationarity of the increments, one may assume that $s=0$. By Definition \ref{def:1:Area},
\BEA \esper[{\cal A}_{0,t}(\eta)^{2N}] &=& \left( \int_0^t dx_1 \int_0^{x_1} dy_1 \right) \ldots
\left( \int_0^t dx_{2N} \int_0^{x_{2N}} dy_{2N} \right) \nonumber\\
&&  \esper\left[B_{x_1}^{'(1)}(\eta) B_{y_1}^{'(2)}(\eta)
\cdots B_{x_{2N}}^{'(1)}(\eta) B_{y_{2N}}^{'(2)}(\eta) \right] \nonumber\\
\EEA

Now apply Proposition \ref{prop:1:Wick} and Definition \ref{def:1:K'K}. \hfill \eop

\bigskip

Every product of $K,K'$ in the above sum may be represented by a diagram.
Draw  a simple line $x$ \textemdash\  $y$ for the infinitesimal kernel $K'(\eta;x,y)$,
a dashed line $x\  - - - -\  y$ for  the integrated kernel
$K(\eta;x,y)$, and a double line $x\  = \   y$ for the kernel 
 $K^*(\eta;x,y)$.

 Each point
$x_1,\ldots,x_{2N}$ is connected to two points, with a simple line for one of them and
a dashed line for the other  (the two points  coincide  in the case of the {\it trivial diagram}
with only two points). Hence each diagram falls into
a number of connected components, each one consisting of a simple bipartite closed polygonal line
whose 2n edges, $n\ge 1$ are alternatively simple and dashed lines (see figure below).

\begin{figure}[h]
  \centering
   \includegraphics[scale=0.7]{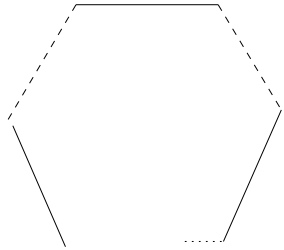}
   \caption{Connected diagram.}
  \label{Fig1}
\end{figure}

There exist $(2N-1)!$ bipartite closed polygonal lines with fixed vertices $(x_1,\ldots,x_{2N})$, i.e. $(2N-1)!$
{\em connected diagrams}
(namely, take as first point $x_1$, then choose any  vertex among the $(2N-1)$ remaining and
connect it to $x_1$ by a simple line, and so on).

\medskip

\begin{Definition}[generating functions]
{\it

Let 
\BEQ \phi_{s,t}(\eta;\lambda):=\sum_{N\ge 0} (-1)^N \esper[{\cal A}_{s,t}(\eta)^{2N}] \frac{\lambda^{2N}}{(2N)!}
=\esper[e^{\II \lambda {\cal A}_{s,t}(\eta)}] \EEQ
be the generating function of  ${\cal A}_{s,t}(\eta)$, and
 $\phi^{(c)}(\eta;\lambda)$ be the generating function of all {\em connected diagrams} (i.e. bipartite
closed polygonal lines). In other words,
\BEQ \phi_{s,t}^{(c)}(\eta;\lambda)=\sum_{N\ge 1} (-1)^N \frac{ \phi^{(c)}_{2N}(\eta;s,t)}{2N} \lambda^{2N},\EEQ
where
\BEA
 \phi^{(c)}_{2N}(\eta;s,t)&& = \int_0^{t-s} dx_1 \ldots \int_0^{t-s}  dx_{2N} \nonumber\\
&&  
\left[ K(\eta;x_1,x_2)K'(\eta;x_2,x_3)\cdots
K(\eta;x_{2N-1},x_{2N})\right] K'(\eta;x_{2N},x_1). \nonumber\\ \label{eq:1:phic} \EEA
}
\label{def:1:generating-functions}
\end{Definition}

Note one has $2N$ in the denominator instead of $(2N)!$ because of the $(2N-1)!$ equivalent connected diagrams.

\medskip

\begin{Lemma}
{\it
The generating function of the L\'evy area is the exponential of the generating functional of
connected diagrams, i.e.
\BEQ \phi_{s,t}(\eta;\lambda)=\exp \phi_{s,t}^{(c)}(\eta;\lambda).\EEQ
}
\label{lemma:1:exp}
\end{Lemma}

{\bf Proof.}

A general bipartite diagram with $2N$ vertices $x_1,\ldots,x_{2N}$ 
 may be decomposed into its connected components, which define a partition
of the set $(x_1,\ldots,x_{2N})$ into $N_1$ subsets of two vertices, $N_2$ subsets of
four vertices, $N_3$ subsets of six vertices and so on. The number of such partitions is
 $\frac{(2N)!}{ (2!)^{N_1} (4!)^{N_2} (6!)^{N_3} \cdots} \frac{1}{N_1! N_2! N_3!\cdots}.$
Now
\BEA \phi_{s,t}(\eta;\lambda)&=&
 \sum_{N_1,N_2,N_3\ldots\ge 0} \frac{(\II\lambda)^{2(N_1+2N_2+3N_3+\ldots)}}{(2!)^{N_1} (4!)^{N_2}(6!)^{N_3}\cdots}
\frac{1}{N_1 !N_2 ! N_3 !\cdots} ( \phi_2^{(c)})^{N_1} (3!\ \phi_4^{(c)})^{N_2} (5!\ \phi_6^{(c)})^{N_3}\cdots \nonumber\\
&=& \exp \left(- \frac{\lambda^2}{2} \phi_2^{(c)} +\frac{\lambda^4}{4} \phi_4^{(c)}-\frac{\lambda^6}{6} \phi_6^{(c)} +\ldots \right) \nonumber\\
&=& \exp \phi^{(c)}(\eta;\lambda).
\EEA
 \hfill \eop

These Feynmann diagram techniques are standard in quantum field theory, see \cite{LeBellac} section 5.3.2 for instance.

\bigskip

Turning now to the connected diagram of order $2N$ -- which is the main object of this section --, it may be split
into the sum of a number of terms by decomposing $K(\eta;x,y)$ into $K(\eta;x,y)=K^*(\eta;x,y)+\frac{1}{2\cos\pi\alpha}
\Re (-\II x+\eta)^{2\alpha} +\frac{1}{2\cos\pi\alpha} \Re (-\II y+\eta)^{2\alpha}$. Replace the $N$ dashed lines
of the closed bipartite polynomial line by a double line whenever $K^*$ replaces $K$, and draw a bullet ($\newmoon$)
at each point $x$
where the function $\Re (-\II x+\eta)^{2\alpha}$ has been inserted instead. Then one has:

-- {\it one closed connected  bipartite diagram} with alternating simple and double lines;

-- and {\it a number of open diagrams} with $n$ components, $n\le N$, each component consisting of alternating simple
and double lines of one of the following three types:
\BEQ (\emptyset \emptyset) \quad \quad   -\ =\ -\ =\ \ldots\ =\ - \EEQ
\BEQ (\emptyset \newmoon) \quad \quad  -\ =\ -\ =\ \ldots \ =\ -\ \newmoon \EEQ
\BEQ (\newmoon \newmoon) \quad \quad \newmoon\  -\ =\ -\ =\ \ldots \ =\ -\ \newmoon \EEQ


\subsection{Preliminary computations and general scheme}


As mentioned in the Introduction and in the discussion preceding Proposition \ref{prop:1}, 
the kernel $K^{',\pm}(\eta)(x,y)$ $(x,y\in\R)$ is the trace on the horizontal
line $\R+\II\frac{\eta}{2}$  of a positive-definite kernel $K^{',\pm}$ defined on $\Pi^{\mp}\times\Pi^{\pm}$:

\begin{Definition}

{\it
For $z\in\Pi^{\mp}$ and $\bar{w}\in\Pi^{\pm}$, we let
\BEQ K^{',\pm}(z,\bar{w})=\frac{\alpha(1-2\alpha)}{2\cos\pi\alpha} (\pm\II(z-\bar{w}))^{2\alpha-2}, \EEQ
\BEQ K^{\pm}(z,\bar{w})=\frac{1}{4\cos\pi\alpha} \left( (\pm \II z)^{2\alpha}+(\mp \II \bar{w})^{2\alpha}-  (\pm\II(z-\bar{w}))^{2\alpha} 
\right) \EEQ
and
\BEQ K^{*,\pm}(z,\bar{w})=-\frac{1}{4\cos\pi\alpha} (\pm\II(z-\bar{w}))^{2\alpha} \EEQ
so that 
\BEQ K^{',\pm}(\eta)(x,y)=K^{',\pm}(x\mp\II\frac{\eta}{2},y\pm\II\frac{\eta}{2}) \EEQ
and similarly for $K^{\pm}(\eta)$ and $K^{*,\pm}(\eta)$.
}

\label{def:2:K'K+-}
\end{Definition}

\begin{Definition}

{\it
For $f\in L^1([a,b],\C)$,  $a,b\in\R$,  define
\BEA
 (K^{',\pm}_{[a,b]}f)(z) &:=& \int_a^b f(u) K^{',\pm}(z,u)\ du \nonumber\\
&=& \frac{\alpha(1-2\alpha)}{2\cos\pi\alpha} \int_a^b f(u) (\pm \II(z-u))^{2\alpha-2}\ du, \quad z\in \Pi^{\mp} \nonumber\\
\EEA
and
\BEA
 (K^{*,\pm}_{[a,b]}f)(z) &:=& \int_a^b f(u) K^{*,\pm}(z,u)\ du \nonumber\\
&=& - \frac{1}{4\cos\pi\alpha}  \int_a^b f(u) (\pm \II(z-u))^{2\alpha}\ du, \quad z\in \Pi^{\mp}. \nonumber\\ \EEA
Both $K^{',\pm}_{[a,b]}f$ and $K^{*,\pm}_{[a,b]}f$ are analytic functions on $\Pi^{\mp}$.
 Similarly, the operators $K^{',\pm}_{[a,b]}(\eta)$, resp.
$K^{*,\pm}_{[a,b]}(\eta)$ are obtained by integrating against the $\eta$ approximations of the kernels
$K^{',\pm}$, resp. $K^{*,\pm}$, so that they may be extended analytically to a complex neighbourhood of the real axis,
\BEQ  (K^{',\pm}_{[a,b]}(\eta)f)(z)=\frac{\alpha(1-2\alpha)}{2\cos\pi\alpha} \int_a^b f(u) (\pm \II(z-u)+\eta)^{2\alpha-2}\ du, \quad z\in \bar{\Pi}^{\mp}
\EEQ  with $\bar{\Pi}^+:=\{\Im z\ge 0\}$, $\bar{\Pi}^-:=\{\Im z\le 0\}$,
and
\BEQ
 (K^{*,\pm}_{[a,b]}(\eta)f)(z)= -\frac{1}{4\cos\pi\alpha}  \int_a^b f(u) (\pm \II(z-u)+\eta)^{2\alpha}\ du, \quad z\in \bar{\Pi}^{\mp}.
\EEQ

Finally, set $K'_{[a,b]}=K^{',+}_{[a,b]}+K^{',-}_{[a,b]}=2\Re K^{',+}_{[a,b]}$ and similarly for the five
other kernels.
}
\label{def:2:K'K+-ab}
\end{Definition}

It is clear from the definition of $K^{',\pm}$ and $K^{*,\pm}$ that $K^{',\pm}_{[a,b]}f$ and $K^{*,\pm}_{[a,b]}f$ are
 well-defined and
analytic on the domain $\C\setminus\{s\pm\II y\ |\ a \le s\le b, y\ge 0\}$, but one needs larger domains of convergence (including
if possible the closed interval $[a,b]$ or at least the open interval $(a,b)$) since eventually one is interested in {\em real}
 variables.
Appropriate maximal domains of analytic extension for $K^{',\pm}_{[a,b]} f$
and $K^{*,\pm}_{[a,b]} f$ are given in the Appendix
under some conditions on the function $f$. The reader interested in the details of the proofs should first look
 at the results of the Appendix,
since we shall constantly be refering to them.
Let us give the general idea for the convenience of the reader. Though the results in the Appendix apply to somewhat
 more general functions
somtimes, we shall only need to consider functions $f$ which are holomorphic on a complex neighbourhood of $(a,b)$, 
with a possible power
behaviour at one of the ends ($a$, say) of the interval, namely $f(z)=(z-a)^{\beta}F(z)$ 
($\beta=0$ or $\beta\in\R\setminus\Z$) for some function $F$
 which is holomorphic in a neighbourhood
of $a$.  The results in the Appendix show that $K^{',\pm}_{[a,b]}f$, $K^{*,\pm}_{[a,b]}f$ are then holomorphic on a complex neighbourhood
of $(a,b)$ and multivalued at $a,b$ with prescribed exponents (see the first lines of the Appendix and next
 subsection 
for more precise statements, in particular
when the interval $[a,b]$ is not fixed). The computation of these exponents is fundamental for our study.

\bigskip

Here is the general scheme of this section. Formula  (\ref{eq:1:phic}) may be rewritten in the following way. Let
$b=x_{2N}$ and $u=x_1$,
\BEA
 \phi_{2N}^{(c)}(\eta;s,t) &&= -\frac{1}{2\cos\pi\alpha} \Re \int_0^{t-s} du K'(\eta;b,u) \ .\
 (K_{[0,t-s]}(\eta) K'_{[0,t-s]}(\eta))^{N-1} \nonumber\\
&& \left(u\mapsto
(\pm\II(u-b)+\eta)^{2\alpha}-(\pm\II u+\eta)^{2\alpha}-(\mp\II b+\eta)^{2\alpha}  \right)(u). \nonumber\\
\EEA
Replacing each occurrence of $K(\eta;x,y)$ by $K^*(\eta;x,y)+\frac{1}{2\cos\pi\alpha} \Re(-\II x+\eta)^{2\alpha}
+\frac{1}{2\cos\pi\alpha} \Re (-\II y+\eta)^{2\alpha}$, see end of subsection 1.1,
 leads (up to some coefficient) either to the single closed diagram
\BEQ \int_0^{t-s} du\  K'(\eta;b,u) \ .\ (K^*_{[0,t-s]}(\eta) K'_{[0,t-s]}(\eta))^{N-1}  \left(u\mapsto
(\pm\II (u-b)+\eta)^{2\alpha} \right)(u) \label{eq:2:K*K'} \EEQ
or to products of terms of the type
\BEQ
 \int_0^{t-s} du  (\II \sigma u+\eta)^{\gamma} (K'_{[0,t-s]}(\eta) K^*_{[0,t-s]}(\eta))^n  \left(u\mapsto
(\pm\II(u-b)+\eta)^{2\alpha-2} \right)(u) (\II \sigma' u+\eta)^{\gamma'} \nonumber\\
 \EEQ
where $\gamma,\gamma'=0$ or $2\alpha$,  and $\sigma,\sigma'\in\{\pm 1\}$. 

The main task will be to estimate the iterated integrals 
$$(K^*_{[0,t-s]}(\eta) K'_{[0,t-s]}(\eta))^n   \left(u\mapsto
(\pm\II(u-b)+\eta)^{2\alpha} \right)(u)$$ and 
 $$(K'_{[0,t-s]}(\eta) K^*_{[0,t-s]}(\eta))^n   \left(u\mapsto
(\pm\II(u-b)+\eta)^{2\alpha-2} \right)(u).$$

We shall content ourselves in this introductory subsection with evaluating 
 the expression $K^{',\sigma_1}_{[0,t]}(\eta) (u\mapsto (\II\sigma_2(u-b)+\eta)^{2\alpha})(a)$
($\sigma_1,\sigma_2\in\{\pm 1\}$), which is (up to a coefficient) equal to the integral $I_{\pm}$ defined in the
Introduction. We shall actually need more general formulas involving arbitrary powers. It turns out that the case $\sigma_1=-\sigma_2$ -- leading to $I_-$ -- 
(see Lemma \ref{lemma:2:lemma3}) is very different from the case $\sigma_1=\sigma_2$ -- leading to
$I_+$ -- (see Lemma \ref{lemma:2:lemma3bis}), the latter case involving
a multivalued term of the form $(\pm\II(a-b)+2\eta)^{\gamma}$ for some $\gamma\in \R$. 

Note in the following formulas that we have in mind $\beta_1=2\alpha$ or $2\alpha-2$, while $\beta_2$ may be much more general (see
comments just before subsection 1.3).

\begin{Lemma}

{\it
Let, for $\beta_1,\beta_2\in\R$ with $\beta_2>-1$, 
\BEQ I_-(\beta_1,\beta_2;0,t)(a,b):=\int_0^t (-\II(u-a))^{\beta_1} (-\II(u-b))^{\beta_2} \ du,\EEQ
defined a priori, {\em for every fixed complex number $b$ with $\Im b\le 0$},  as an analytic function of $a$  on $\Pi^-$. We restrict
to $0<\Re a<t$ and $0<\Re b<t$, $\Im b\le 0$. Let $\Omega^-_t:= \{a\in \C\ |\ 0<\Re a<t, \Im a<\Im b\}$.  Then the following results hold:

\begin{itemize}
\item[(i)] On the domain $\Omega^-_t$, one has
\BEQ I_-(\beta_1,\beta_2;0,t)(a,b)=\frac{\II}{\beta_1+\beta_2+1} \left[ \Phi(\beta_1,\beta_2;t)(a,b)-\Phi(\beta_1,\beta_2;0)(a,b)
\right] \EEQ
with, for $s\in[0,t]$,
\BEQ \Phi(\beta_1,\beta_2;s)(a,b)=(-\II(s-b))^{\beta_1+\beta_2+1} \ _2 F_1(-\beta_1,-\beta_1-\beta_2-1;-\beta_1-\beta_2;
\frac{a-b}{s-b}).\EEQ

The function $a\mapsto \Phi(\beta_1,\beta_2;0)(a,b)$, resp. $a\mapsto \Phi(\beta_1,\beta_2;t)(a,b)$ given
by the above expression has an analytic extension to the domain $\{0<c<|a/b|<C\}\cap \{0<\Re a<t\}$,
resp.  $\{ 0<c<|\frac{t-b}{t-a}|<C\}\cap \{0<\Re a<t\}$ with arbitrary constants 
$c<1$, $C>1$. Both functions extend analytically to 
the whole domain $\{0<\Re a<t\}$, with different expressions given below.

\item[(ii)]

Suppose $|a/b|<c<1$. Then
\BEA
&& \Phi(\beta_1,\beta_2;0)(a,b)=(\II b)^{\beta_1+\beta_2+1} \ .\ \left\{ \frac{\Gamma(-\beta_1-\beta_2)\Gamma(1+\beta_1)}{\Gamma(-\beta_2)} (1-\frac{a}{b})^{\beta_1+\beta_2+1}  \right. \nonumber\\ 
&& \left. +\frac{\beta_1+\beta_2+1}{\beta_1+1} \left(\frac{a}{b}\right)^{1+\beta_1} \ _2 F_1(-\beta_2,1;\beta_1+2;a/b) \right\}
\nonumber\\
\EEA
extends analytically to the domain $\{|a/b|<c<1 \}\cap\{0<\Re a<t\}.$

\item[(iii)]

Suppose $|a/b|>C>1$. Then
\BEA
&& \Phi(\beta_1,\beta_2;0)(a,b)=(\II b)^{\beta_1+\beta_2+1} \ .\ \left\{ \frac{\Gamma(-\beta_1-\beta_2)\Gamma(1+\beta_2)}{\Gamma(-\beta_1)} \left(\frac{b}{a}\right)^{-\beta_1-\beta_2-1} \right. \nonumber\\
&& \left.  (1-\frac{b}{a})^{\beta_1+\beta_2+1}  +\frac{\beta_1+\beta_2+1}{\beta_2+1} \left(\frac{b}{a}\right)^{-\beta_1} \ _2 F_1(-\beta_1,1;\beta_2+2;b/a) \right\}
\nonumber\\
\EEA

extends analytically to the domain $\{|a/b|>C>1 \}\cap\{0<\Re a<t\}.$

\item[(iv)]

Similarly, suppose $|\frac{t-a}{t-b}|<c<1$. Then
\BEA
&& \Phi(\beta_1,\beta_2;t)(a,b)=(-\II(t-b))^{\beta_1+\beta_2+1} \ .\ \left\{ \frac{\Gamma(-\beta_1-\beta_2)\Gamma(1+\beta_1)}{\Gamma(-\beta_2)} \left( \frac{a-b}{t-b} \right)^{\beta_1+\beta_2+1}  \right. \nonumber\\ 
&& \left. +\frac{\beta_1+\beta_2+1}{\beta_1+1} \left(\frac{t-a}{t-b}\right)^{1+\beta_1} \ _2 F_1(-\beta_2,1;\beta_1+2;
\frac{t-a}{t-b}) \right\}
\nonumber\\
\EEA
extends analytically to the domain $\{|\frac{t-a}{t-b}|<c<1 \}\cap\{0<\Re a<t\}.$

\item[(v)]

Suppose $|\frac{t-a}{t-b}|>C>1$. Then
\BEA
&& \Phi(\beta_1,\beta_2;t)(a,b)=(-\II(t-b))^{\beta_1+\beta_2+1} \ .\ \left\{ \frac{\Gamma(-\beta_1-\beta_2)\Gamma(1+\beta_2)}{\Gamma(-\beta_1)} \left(\frac{t-b}{t-a}\right)^{-\beta_1-\beta_2-1}  \right. \nonumber\\ 
&& \left. (\frac{b-a}{t-a})^{\beta_1+\beta_2+1}  +\frac{\beta_1+\beta_2+1}{\beta_2+1} \left(\frac{t-b}{t-a}\right)^{-\beta_1} \ _2 F_1(-\beta_1,1;\beta_2+2;\frac{t-b}{t-a}) \right\}
\nonumber\\
\EEA

extends analytically to the domain $\{|\frac{t-a}{t-b}|>C>1 \}\cap\{0<\Re a<t\}.$
\end{itemize}
}
\label{lemma:2:lemma3}
\end{Lemma}

\begin{Remark}
{\it  Note that $I_-$ behaves as a sum of power functions when $a$ and/or $b$ are in the neighbourhood of
either of the interval ends. All together, one gets the following expressions (assuming to simplify notations
that $a,b\in (0,t)$ are real):
\BEQ \Phi(\beta_1,\beta_2;0)(a,b)\sim C_1 (\max(a,b))^{\beta_1+\beta_2+1}+C_2 a^{\beta_1+1}b^{\beta_2} 1_{a<b}
+C_3 b^{\beta_2+1} a^{\beta_1} 1_{a>b} \EEQ
if $a$ or $b$ is close to $0$, and similarly
\BEA
&&   \Phi(\beta_1,\beta_2;t)(a,b)\sim C_1 (\max(t-a,t-b))^{\beta_1+\beta_2+1} \nonumber\\
&& +C_2
 (t-a)^{\beta_1+1} (t-b)^{\beta_2} 1_{t-a<t-b}
+C_3 (t-b)^{\beta_2+1} (t-a)^{\beta_1} 1_{t-a>t-b} \nonumber\\ \EEA
if $a$ or $b$ is close to $t$. Now sum up the contributions of these two terms to get the exponents of $I_-$.
}
\label{rmk:1}
\end{Remark}

\bigskip

{\bf Proof.}

(i) is proved in \cite{Unt08}, Lemma 4.1 with  slightly different hypotheses. Let us give a self-contained proof with easier
arguments. Decompose the integral $\int_0^t du$ into $\int_0^b du+\int_b^t du$. The first integral writes

\BEA &&  \int_0^b (-\II(u-a))^{\beta_1} (-\II(u-b))^{\beta_2}\ du=b\int_0^1 \left[ \II(vb+(a-b)) \right]^{\beta_1}
(\II vb)^{\beta_2}\ dv. \nonumber\\
\EEA

By hypothesis, $\Im a<\Im b\le 0$ and $0<\Re a,\Re b<t$. Suppose for the moment that  $\Re a< \Re b$. Then  one has
$-\frac{\pi}{2}<\Arg(\II(a-b))<0$ and $0<\Arg(1-\frac{vb}{b-a})<\pi$, hence
\BEQ \left[ \II(vb+(a-b)) \right]^{\beta_1}=(\II(a-b))^{\beta_1} (1-\frac{vb}{b-a})^{\beta_1} \EEQ
and one obtains
\BEQ \int_0^b (-\II(u-a))^{\beta_1} (-\II(u-b))^{\beta_2}\ du=-\II \frac{(\II b)^{\beta_2+1}}{\beta_2+1}
(-\II(b-a))^{\beta_1} \ _2 F_1(-\beta_1,\beta_2+1;\beta_2+2;\frac{b}{b-a}).\EEQ

Similarly,
\BEQ 
\int_b^t (-\II(u-a))^{\beta_1} (-\II(u-b))^{\beta_2}\ du=(t-b)\int_0^1 \left[ -\II(w(t-b)+(b-a)) \right]^{\beta_1}
(-\II(t-b)w)^{\beta_2}\ dw.
\EEQ

Now (still assuming $\Re a<\Re b$)
$-\frac{\pi}{2}<\Arg(-\II(b-a))<0$ and $-\frac{\pi}{2}<\Arg(1-\frac{w(b-t)}{b-a})<\frac{\pi}{2}$, hence
\BEQ \left[ -\II(w(t-b)+(b-a)) \right]^{\beta_1}=(-\II(b-a))^{\beta_1} (1-w\frac{b-t}{b-a})^{\beta_1} \EEQ
and one obtains
\BEQ \int_b^t (-\II(u-a))^{\beta_1} (-\II(u-b))^{\beta_2}\ du=\II \frac{(-\II (t-b))^{\beta_2+1}}{\beta_2+1}
(-\II(b-a))^{\beta_1} \ _2 F_1(-\beta_1,\beta_2+1;\beta_2+2;\frac{b-t}{b-a})\EEQ

Observe one has obtained
\BEQ I_-(\beta_1,\beta_2;t)(a,b)=\frac{\II}{\beta_2+1} (-\II(b-a))^{\beta_1} \left[
F(\beta_1,\beta_2;t)(a,b)-F(\beta_1,\beta_2;0)(a,b)\right] \EEQ
with
\BEQ F(\beta_1,\beta_2;s)(a,b)=(-\II(s-b))^{\beta_2+1} \ _2 F_1(-\beta_1,\beta_2+1;\beta_2+2;\frac{s-b}{a-b}),\quad
s\in(0,t).\EEQ One may lift the restriction $\Re a<\Re b$ since the last expression makes sense for all
$a\in\C$ such that $\Im a<\Im b\le 0$ (namely, $s-b,\frac{1}{a-b}\in\Pi^+$ so that $\frac{s-b}{a-b}\not\in\R_+$).

 Now, by the connection formula (\ref{eq:0:1/z}) (reproducing an argument in \cite{Unt08})
\BEA
&& _2 F_1(-\beta_1,\beta_2+1;\beta_2+2;\frac{s-b}{a-b})=\frac{\beta_2+1}{\beta_1+\beta_2+1} \left(\frac{s-b}{b-a}\right)^{\beta_1}
\nonumber\\ 
&&  _2 F_1(-\beta_1,-\beta_1-\beta_2-1;-\beta_1-\beta_2;\frac{a-b}{s-b}) + \frac{\Gamma(\beta_2+2)\Gamma(-\beta_1-\beta_2-1)}{\Gamma(-\beta_1)} \left(\frac{b-a}{s-b}\right)^{1+\beta_2}.\nonumber\\
\label{connI+-}
\EEA

Now $\left(\frac{b-a}{s-b}\right)^{1+\beta_2}=\frac{(b-a)^{1+\beta_2}}{(s-b)^{1+\beta_2}}$, so 
the second term, multiplied by the prefactor $(-\II(s-b))^{\beta_2+1}$, is independent of $s\in(0,t)$ and hence
makes no contribution to $I_-$.

\bigskip

Now (ii) and (iv), resp. (iii) and (v),  are consequences of the connection formula (\ref{eq:0:1-z}), resp. (\ref{eq:0:1/1-z}).
\hfill \eop

Note (as mentioned in the Introduction) that the fact that $I_-$ is analytic in the parameters $a,b$ away from $0,t$ may
easily be proved by using a deformation of contour in $\Pi^+$.

\bigskip

Let us now turn to the non-analytic case:

\begin{Lemma}

{\it
Let, for $\beta_1,\beta_2\in\R$ with $\beta_2>-1$, 
\BEQ I_+(\beta_1,\beta_2;0,t)(a,b):=\int_0^t (\II(u-a))^{\beta_1} (-\II(u-b))^{\beta_2} \ du,\EEQ
defined a priori, {\em for every fixed complex number $b$ such that $\Im b\le 0$},  as an analytic function of $a$  on $\Pi^+$. We restrict
to $0<\Re a<t$ and $0<\Re b<t$. Let $\Omega^+_t:= \{a\in \Pi^+\ |\ 0<\Re a<t\}$.  Then the following 
formula holds on $\Omega^+_t$:

\BEA 
&& I_+(\beta_1,\beta_2;0,t)(a,b)=\frac{\II}{\beta_1+\beta_2+1} \left[  e^{\II\pi\beta_1}  
 \Phi(\beta_1,\beta_2;t)(a,b)- e^{-\II\pi\beta_1} \Phi(\beta_1,\beta_2;0)(a,b)
\right] \nonumber\\
&& -\frac{\Gamma(\beta_2+1)\Gamma(-\beta_1-\beta_2-1)}{\Gamma(-\beta_1)} \ .\ 2\sin\pi\beta_2 \ .\ (\II(b-a))^{\beta_1+\beta_2+1}.
\nonumber\\ \label{eq:2:I+}
 \EEA

where $\Phi$ is the same analytic function as in Lemma \ref{lemma:2:lemma3}.

Contrary to what happens for the $I_-$ integral, the last term does not admit an analytic extension in $a$ to any neighbourhood
of $b$ (one cannot 'circle' around $b$).
}
\label{lemma:2:lemma3bis}
\end{Lemma}

{\bf Proof.}

A variant of this Lemma is also proved in Lemma 4.1 of \cite{Unt08}, but let us give an independent proof. The beginning
is as in Lemma \ref{lemma:2:lemma3}. Namely,
\BEQ \int_0^b du (\II(u-a))^{\beta_1} (-\II(u-b))^{\beta_2}=b\int_0^1 dv\ \left[ -\II(vb+(a-b))\right]^{\beta_1}
(\II vb)^{\beta_2} \EEQ
and
\BEQ \int_b^t du  (\II(u-a))^{\beta_1} (-\II(u-b))^{\beta_2}+(t-b)\int_0^1  dw\ \left[ \II(w(t-b)+(b-a)) \right]^{\beta_1}
(-\II(t-b)w)^{\beta_2}.\EEQ

Mind $a\in\Pi^+$ this time. Suppose provisorily that  $\Re a< \Re b$: then
\BEQ 0<\Arg(-\II(a-b))<\pi/2, \quad -\pi<\Arg(1-\frac{vb}{b-a})<0 \EEQ
(since $\frac{b}{a-b}=\left(\frac{a}{b}-1\right)^{-1}\in \Pi^-$) and
\BEQ 0<\Arg(-\II(a-b))+\Arg(1-w\frac{b-t}{b-a})<\pi \EEQ
since $0<\Arg(1-w\frac{b-t}{b-a})<\Arg \frac{t-b}{b-a}<\pi$ and $-\II(a-b).\frac{t-b}{b-a}=\II(t-b)\in\Pi^+$.

Hence
\BEA
&&  I_+(\beta_1,\beta_2;t)(a,b)=(\II(b-a))^{\beta_1} \frac{\II}{\beta_2+1} \left\{ (-\II(t-b))^{\beta_2+1}
\ _2 F_1(-\beta_1,\beta_2+1;\beta_2+2;\frac{t-b}{a-b}) \right. \nonumber\\
&& \left. - (\II b)^{\beta_2+1} \ _2 F_1(-\beta_1,\beta_2+1;\beta_2+2;\frac{b}{b-a}) \right\}
\EEA
$(a\in\Pi^+)$ which is the same formula as in the proof of Lemma \ref{lemma:2:lemma3} {\it except} $a\in\Pi^+$ and the prefactor is $(\II(b-a))^{\beta_1}$
instead of $(-\II(b-a))^{\beta_1}$. Apply the connection formula (\ref{eq:0:1/z}). Unfortunately the second terms
in equation (\ref{connI+-}) do not cancel each other this time. Namely, they contribute (up to a constant prefactor) the
following expression:
\BEA
&& J:=(\II(b-a))^{\beta_1} \left[ (-\II(t-b))^{\beta_2+1} \left(\frac{t-b}{b-a}\right)^{-\beta_2-1}
- (\II b)^{\beta_2+1} \left(\frac{b}{a-b}\right)^{-\beta_2-1} \right]  \nonumber\\
&&=(\II(b-a))^{\beta_1} \left[ e^{-\II\frac{\pi}{2}(\beta_2+1)} (b-a)^{\beta_2+1} - e^{\II\frac{\pi}{2}(\beta_2+1)}
(a-b)^{\beta_2+1} \right]. \nonumber\\
\EEA

Since $a\in\Pi^+$, one has $(b-a)^{\beta_2+1}=e^{-\II\frac{\pi}{2}(\beta_2+1)} (\II(b-a))^{\beta_2+1}$  and
$(a-b)^{\beta_2+1}=e^{\II\frac{\pi}{2}(\beta_2+1)} (-\II(a-b))^{\beta_2+1}$, whence
\BEQ J=(\II(b-a))^{\beta_1+\beta_2+1} 2\II \sin\pi\beta_2.\EEQ

Now the first terms in equation (\ref{connI+-}) come up with a supplementary prefactor $e^{\II\pi\beta_1}$ with respect
to the formulas in Lemma \ref{lemma:2:lemma3} because
\BEQ (\II(b-a))^{\beta_1} \ . \left(\frac{t-b}{b-a}\right)^{\beta_1}=(\II(b-a))^{\beta_1}
\left(\frac{\II(t-b)}{\II(b-a)}\right)^{\beta_1}=(\II(t-b))^{\beta_1}=(-\II(t-b))^{\beta_1} e^{\II\pi\beta_1} \EEQ
and similarly
\BEQ (\II(b-a))^{\beta_1}\ . \ \left(\frac{-b}{b-a}\right)^{\beta_1}=(-\II b)^{\beta_1}=(\II b)^{\beta_1} e^{-\II\pi\beta_1}.\EEQ

\hfill \eop

\begin{Remark}
{\it
 {\em Both} integrals $I_{\pm}(\beta_1,\beta_2;0,t;a,b)$ extend analytically to the product of the cut planes
$\left( \C\setminus(\R_-\cup(t+\R_+)) \right)^2$, and  behave as in Lemma \ref{lemma:2:lemma3}  (i.e. 
with the same power functions,
see Remark \ref{rmk:1},
and {\em without} the non-analytic term in $(\II(b-a))^{\beta_1+\beta_2+1}$) if $\Re a\in\R\setminus [0,t]$ or
$\Re b\in\R\setminus [0,t]$. This result is a consequence of Lemmas \ref{lemma:2:lemma3} and \ref{lemma:2:lemma3bis}, and
also Lemmas \ref{lemma:2:lemma1},\ref{lemma:2:lemma2} in the Appendix.  If $\Re a,\Re b<0$ (or similarly if $\Re a,\Re b>t$) then
$I_{\pm}(\beta_1,\beta_2;0,t;a,b)=I_{\pm}(\beta_1,\beta_2;-s,t;a,b)-I_{\pm}(\beta_1,\beta_2;-s,0;a,b)$ for
$-s<\Re a,\Re b$, so the non-analytic term (in the case of $I_+$) disappears. If $\Re a<0$ and $\Re b>t$ for instance (so
$a$ and $b$ are far away) then one may cut the interval of integration into $[0,t/2]\cup [t/2,t]$ and apply Lemmas
\ref{lemma:2:lemma1},\ref{lemma:2:lemma2}. Finally, if $\Re b\in(0,t)$ and, say, $\Re a<0$, then $I_+$ reduces
 to $I_-$  since $(\II(u-a))^{\beta_1}=e^{\II\pi\beta_1}(-\II(u-a))^{\beta_1}$ for every $u\in(0,t)$. Then the formulas appearing
in Lemma \ref{lemma:2:lemma3} may be extended analytically to $\Re a<0$ with a little care.
}
\label{rmk:2}
\end{Remark}

\bigskip

 The extra non-analytic power term proportional to  $a\mapsto (-\II(a-b))^{\beta_1+\beta_2+1}$ in equation (\ref{eq:2:I+}) may in turn be integrated
against $K^{',\pm}$ or $K^{*,\pm}$. Generally speaking, alternate chains of the form $(K^{*,-}(\eta)K^{',-}(\eta))^n 
\left( a\mapsto
(-\II(a-b)+\eta)^{2\alpha} \right)$ or $K^{',-}(\eta)(K^{*,-}(\eta)K^{',-}(\eta))^n \left( a\mapsto
(-\II(a-b)+\eta)^{2\alpha} \right)$ (or conjugate) contain an extra power term with increasing exponent (hence the need for a general
exponent $\beta_2$). Note that integrating against $K^{*,+}(\eta)$ or $K^{',+}(\eta)$
 at some point kills the power term (by Lemma
\ref{lemma:2:lemma3}) and produces a function (depending on $b$) 
 which is analytic on a complex neighbourhood of $(0,t)$. This remark is fundamental to
understand the divergence of the L\'evy area for $\alpha<1/4$, which is due to the non-analytic
power terms only, as we shall show in the next 
paragraph.


\subsection{Convergence of the analytic part of the moments}


Leaving aside the non-analytic extra power terms coming from Lemma \ref{lemma:2:lemma3bis}, one is led 
(see preceding subsection, equation (\ref{eq:2:K*K'}))  to evaluate 
 alternating integral chains
of the form $\cdots K^{*,\pm}_{[0,t]}(\eta)K^{',\pm}_{[0,t]}(\eta)\cdots f$ for some function $f$,
depending on $b$,  which is analytic on a complex neighbourhood 
of $(0,t)$, and multivalued (of power type) at $0$ or $t$. Leaving aside the dependence on the variable $b$
and on the parameter $\eta$,
the results in the Appendix show that such integrals are also analytic on a complex
neighbourhood of $(0,t)$, multivalued of power type (with some exponents) at $0$ and $t$, and allow to 
compute the exponents. Let us make  precise statements.
We first need to define an appropriate class of analytic functions $f(\eta,b,t;u)$
 (called {\em admissible} in the sequel),
of the form $\eta^H b^B u^U F(\eta,b,u,t)$ or $\eta^H b^B (t-u)^U F(\eta,b,u,t)$ (where $u\mapsto
F(\eta,b,t;u)$ is supposed
to be  analytic on a neighbourhood of the closed interval $[0,t]$) 
for which such alternating integral chains make sense. Such functions $F$ may only be defined locally because the
values of the exponents $H,B,U$ depend on the relative positions of $\eta,b,u$. This makes the exact
definition look a little complicated at first:

\begin{Definition}[admissible analytic functions]

{\it
Let $\eta>0$ and $\sigma_f\in\Z,\bar{\sigma}\in\{\pm 1\}$.
Assume $f(\eta,b,t;u):=f_b(\eta,b,t;u)+f_{\eta}(\eta,b,t;u)$ where both $f_b$ and $f_{\eta}$  are 
analytic in $\eta$, $b$ and $u$ for $(\eta,b)$ on a complex neighbourhood $\Omega'$
 of $\{\eta\ge 0,b\ge 0\}$ and
$u$ on the cut domain 
$$\Omega:=\{|\frac{u+\II\sigma_f \eta}{t}|<2\} \setminus\left( \{ \frac{u}{\eta}\in-\II\sigma+\R_-\} \cup \{ \frac{t-u}{\eta}\in-\II\sigma+\R_+\}\right)$$
for some $\sigma\in\{\pm 1\}$.  Let $\Omega'_1=\{\frac{\eta}{b+\II\bar{\sigma}\eta};(\eta,b)\in\Omega'\}$ and
$\Omega'_2=\{\frac{b+\II\bar{\sigma}\eta}{t};(\eta,b)\in\Omega'\}$. 
Assume also that $f_b$  may be written as: 
\begin{itemize}

\item[$(i)_b$] $\sum_{j=1}^J  (b+\II \bar{\sigma}\eta)^{B^{(b)}_j}(u+\II \sigma_f \eta)^{U^{(b)}_j} F_j(\frac{u+\II\sigma_f \eta}{b+\II
\bar{\sigma}\eta},\frac{\eta}{b+\II\bar{\sigma}\eta})$
 on the subdomain $(i)_b$ of $\Omega$ such that   $0<|u+\II\sigma_f \eta|<2|b+\II\bar{\sigma}\eta|/3$;

\item[$(ii)_b$] $\sum_{j=1}^{J'}  (b+\II\bar{\sigma}\eta)^{B^{(b)'}_j} F'_j(\frac{u+\II\sigma_f \eta}{b+\II
\bar{\sigma}\eta},\frac{\eta}{b+\II\bar{\sigma}\eta})$
 on the subdomain $(ii)_b$ of $\Omega$ such that $|b+\II\bar{\sigma}\eta|/3<|u+\II\sigma_f \eta|<3|b+\II\bar{\sigma}\eta|$;

\item[$(iii)_b$] $\sum_{j=1}^{J''}  (b+\II\bar{\sigma}\eta)^{B^{(b)''}_j} 
(u+\II\sigma_f \eta)^{U^{(b)''}_j} F''_j(\frac{b+\II\bar{\sigma}\eta}{u+\II\sigma_f \eta},
\frac{\eta}{b+\II\bar{\sigma}\eta})$
 on the subdomain $(iii)_b$ of $\Omega$ such that  $2|b+\II\bar{\sigma}\eta|<|u+\II\sigma_f \eta|$

\end{itemize}

where $F_j$, $F''_j$ are holomorphic on $B(0,1)\times\Omega'_1$, and $F'_j(w,\zeta)$ are holomorphic on $\{1/3<|w|<3\}\times \Omega'_1$

with $U_j^{(b)}>-1$ for all $j$, and
 $\{(B^{(b)}_j+U^{(b)}_j),j=1\ldots J\}=\{(B^{(b)'}_{j'}), j'=1\ldots J'\}=\{(
B^{(b)''}_{j''}+U^{(b)''}_{j''}),
j''=1\ldots J''\}$,
while $f_{\eta}$ may be written as

\begin{itemize}

\item[$(ii)_{\eta}$] $\sum_{j=1}^{J'} \eta^{H^{(\eta)'}_j} (b+\II \bar{\sigma}\eta)^{B^{(\eta)'}_j}
(u+\II\sigma_f \eta)^{U_j^{(\eta)'}}
  F'_j(\frac{u+\II\sigma_f \eta}{\eta},
\frac{\eta}{b+\II\bar{\sigma}\eta},\frac{b+\II\bar{\sigma}\eta}{t})$ on the  subdomain $(i)_{\eta}$ of $\Omega$
such that  
 $0<|u+\II\sigma_f \eta|<3 \eta$;

\item[$(iii)_{\eta}$]
\BEA &&  \sum_{j=1}^{J''} \eta^{H^{(\eta)''}_j} (b+\II\bar{\sigma}\eta)^{B^{(\eta)''}_j} 
(u+\II\sigma_f \eta)^{U^{(\eta)''}_j}  \nonumber\\
&& \quad  \left( F''_j(\frac{ \eta}{u+\II\sigma_f \eta},
\frac{\eta}{b+\II\bar{\sigma} \eta}, \frac{b+\II\bar{\sigma}\eta}{t}) +
G''_j(\frac{u+\II\sigma_f \eta}{t},\frac{\eta}{b+\II\bar{\sigma}\eta},
\frac{b+\II\bar{\sigma}\eta}{t}) \right) \nonumber\\ \EEA
 on the subdomain $(iii)_{\eta}$ of $\Omega$ such that   $2\eta<|u+\II\sigma_f \eta|<2t/3$;

\item[$\wit{(iii)}_{\eta}$] 
\BEA && \sum_{j=1}^{\tilde{J}''} \eta^{\tilde{H}^{(\eta)''}_j} (b+\II\bar{\sigma}\eta)^{\tilde{B}^{(\eta)''}_j} 
(1-\frac{u+\II\sigma_f \eta}{t})^{\tilde{U}^{(\eta)''}_j} \nonumber\\
&&  \left(\tilde{F}''_j(\frac{\eta/t}{1-(u+\II\sigma_f \eta)/t},
\frac{\eta}{b+\II\bar{\sigma} \eta},\frac{b+\II\bar{\sigma}\eta}{t}) + \tilde{G}''_j(1-\frac{u+\II\sigma_f \eta}{t},
\frac{\eta}{b+\II\bar{\sigma}\eta},
\frac{b+\II\bar{\sigma}\eta}{t}) \right) \nonumber\\ \EEA
 on the subdomain $\wit{(iii)}_{\eta}$ of $\Omega$ such that  
 $2 \eta/t <|1-(u+\II\sigma_f \eta)/t|<2/3$;

\item[$\wit{(ii)}_{\eta}$] 
\BEA && \sum_{j=1}^{\tilde{J}'} \eta^{\tilde{H}^{(\eta)'}_j} (b+\II \bar{\sigma}\eta)^{\tilde{B}^{(\eta)'}_j} (1-(u+\II\sigma_f \eta)/t)^{\tilde{U}_j^{(\eta)'}} \nonumber\\
&& 
 \tilde{F}'_j(\frac{1-(u+\II\sigma_f \eta)/t}{ \eta/t}, \frac{\eta}{b+\II\bar{\sigma}\eta},\frac{b+\II\bar{\sigma}\eta}{t})
\nonumber\\ \EEA
 on the subdomain $\wit{(ii)}_{\eta}$ of $\Omega$
such that   $0<|1-(u+\II\sigma_f \eta)/t|<3\eta/t$ 
\end{itemize}   

where $F'_j(w,\zeta,\xi)$, resp. $\tilde{F}'_j(w,\zeta,\xi)$ are holomorphic on
${\cal B}\times\Omega'_1\times\Omega'_2$, resp. $\tilde{\cal B}\times\Omega'_1\times\Omega'_2$, with
 ${\cal B}:=\{|w|<3,w-\II\sigma_f\not\in -\II\sigma+\R_-\}$, $\tilde{\cal B}:=\{|w|<3,w+\II\sigma_f\not\in \II\sigma+\R_-\}$ and
$F''_j(w,\zeta,\xi)$, $\tilde{F}''_j(w,\zeta,\xi)$, $G''_j(w,\zeta,\xi)$, $\tilde{G}''_j(w,\zeta,\xi)$
 are holomorphic on $B(0,1)\times
\Omega'_1\times\Omega'_2$,

with $U^{(\eta)}_j,\tilde{U}^{(\eta)}_j>-1$ for all $j$ and $\{(B^{(\eta)'}_{j'},
H^{(\eta)'}_{j'}+U_j^{(\eta)'}),j'=1\ldots J'\}=
\{(B^{(\eta)''}_{j''}, H^{(\eta)''}_{j''}+U^{(\eta)''}_{j''}),j''=1\ldots J''\}$ (and similarly for the exponents with a tilde), 
$\{(H^{(\eta)''}_j,B^{(\eta)''}_j),j=1\ldots J''\}=\{(\tilde{H}^{(\eta)''}_{\tilde{j}''},\tilde{B}^{(\eta)''}_{\tilde{j}''}),
 \tilde{j}''=1\ldots \tilde{J}''\}.$

Then one says that $f$ is an {\em admissible analytic function} with $b$-exponents $\{(B^{(b)}_j,U^{(b)}_j)\}$ on
the domain $(i)_b$ and $\{(B^{(b)''}_j,U^{(b)''}_j)\}$ on the domain $(iii)_b$, and
$\eta$-exponents  $\{(H^{(\eta)'}_j,B^{(\eta)'}_j,U^{(\eta)'}_j)\}$ on the domain $(ii)_{\eta}$,
 $\{(H^{(\eta)''}_j,B^{(\eta)''}_j,
U^{(\eta)''}_j)\}$ on the domain $(iii)_{\eta}$ (and similarly for the domains with a tilde).

If $f=f_b$, resp. $f=f_{\eta}$, then one says that $f$ is of {\em $b$-type}, resp. of {\em $\eta$-type}.

In particular, the function $\Phi(\beta_1,\beta_2;0)(u,b)$ appearing in Lemma \ref{lemma:2:lemma3} is an
admissible analytic function of $b$-type with $b$-exponents $\{(B^{(b)}_j,U^{(b)}_j)\}=\{(\beta_1+\beta_2+1,0),(\beta_2,\beta_1+1)\}$,
$\{(B^{(b)''}_j,U^{(b)''}_j)\}=\{(0,\beta_1+\beta_2+1),(\beta_2+1,\beta_1)\}$,  while $\Phi(\beta_1,\beta_2;t)(u,b)$ is of the same type up to the symmetry $u\to
t-u$, $b\to t-b$.

}

\label{def:2:admissible}
\end{Definition}

{\bf Remarks.}

\begin{itemize}

\item The extra conditions $w-\II\sigma_f=\frac{u}{\eta}\not\in -\II\sigma+\R_-$,
resp. $w+\II\sigma_f=\frac{t-u}{\eta}\not\in\II\sigma+\R_-$  (which are true when integrating over the
real axis or, more generally, over
any  deformed contour  $\gamma:0\to t$) avoid considering unnecessary complications (the values
of the exponents are different when $u$ or $t-u$ is  in a neighbourhood of $\pm\II\eta$). They also appear in 
the definition of the cut domain $\Omega$. Note that the functions $F_j$, $F'_j$, $F''_j$ associated to the domains
$(i)_b$, $(ii)_b$, $(iii)_b$ are bounded on their respective domains, and so are the functions $F''_j$, $G''_j$, $\tilde{F}''_j$, $\tilde{G}''_j$ on $(iii)_{\eta}$, $\wit{(iii)}_{\eta}$, while $F'_j$, resp. $\tilde{F}'_j$ are possibly unbounded on
$(ii)_{\eta}$, resp. $\wit{(ii)}_{\eta}$ (but they are bounded on the subdomain $\Omega_{res}:=\Omega\setminus\{B(-\II\eta\sigma,
\eta/3)\cup B(t-\II\eta\sigma,\eta/3))$, see also proof of Theorem \ref{th:2:refinement}, which is
already satisfactory since
one is ultimately interested in the behaviour around $[0,t]\subset\Omega_{res}$).

\medskip

The reason for the appearance of the $\sigma$-sign is that admissible analytic functions are in the
image of the integral transformation $K^{*,\sigma}_{[0,t]}(\eta)$ or $K^{',\sigma}_{[0,t]}(\eta)$ for
some $\sigma\in\{\pm 1\}$. Note that
 $\int_0^t (\II\sigma(v-u)+\eta)^{\beta} g(v)\ dv$, $\beta=2\alpha$ or $2\alpha-2$
(for appropriate functions $g$) is possibly singular on the boundary of $\Omega$, but always
regular on the closure of $\Omega_{res}$.

\item
Note that the family of  $b$-exponents $\{(B^{(b)}_j,U^{(b)}_j), (B^{(b)''}_j,U^{(b)''}_j))\}$
and the family of $\eta$-exponents $\{(H^{(\eta)'}_j,B^{(\eta)'}_j,U^{(\eta)'}_j),(H^{(\eta)''}_j,B^{(\eta)''}_j,U^{(\eta)''}_j)\}$
 (together with the relative exponents with
a tilde)  determine all the exponents
of the function $f$. The conditions on the $U_j$'s  ensure in particular that $f_b$ is integrable when $\eta=0$.

\item The number $\bar{\sigma}$ does not change when one integrates $f$ against $K'(\eta)$ or $K(\eta)$.
  When computing the contribution to the
$2N$-th moment of the $\eta$-L\'evy area of the terms containing the kernel $(\II\sigma_{2N-1}(x_{2N-1}-x_{2N}))^{2\alpha}=
(\II\sigma_{2N-1}(x_{2N-1}-b))^{2\alpha}$,
one  simply  has  $\bar{\sigma}=\sigma_{2N-1}$.

\medskip

On the contrary, the value of $\sigma_f$ is shifted by $\pm 1$ after each integration (see Theorem
\ref{th:2:exponents}).

\item Depending on whether the context is clear of not, we shall sometimes drop the upper indices $(b)$ or $(\eta)$ of the exponents.
\item The constants $1/3$, $2/3, \ldots$ appearing in the definition of the domains are arbitrary and may be replaced by
any other set of positive constants, as long as the subdomains intersect.
\item The splitting of $f$ into $f_b+f_{\eta}$ is essentially a 'pedagogical' artefact. One should actually consider a single
function with different expressions and exponents depending on the relative position of $0,\eta,b,t$ and $u$. That would
make the above Definition even more technical, with the following advantage however: assuming $\sigma_f\not=\sigma$, so
$z:=u+\II\sigma_f \eta$ may be arbitrary close to 0 on $\Omega_{res}$, both functions $f_b$ and
$f_{\eta}$ have a singularity when $z=0$ if some exponent $U_j^{(b)}$ or $U_j^{(\eta)'}$ is negative (which does happen
in our case), while $f_b+f_{\eta}$ is analytic at $z=0$ (by definition). Fortunately this is not a real problem for the
convergence proof.

\end{itemize}

The existence of the $b$-exponents follows naturally from the discussion before Definition \ref{def:2:admissible}. They describe the power
behaviour of the function  $f$ near $0$ and $t$. The complications come from the fact that the parameter $b$ itself may be close to $0$ (on the
contrary, if $b$ is bounded from below, say $b>t/2$, then $b$ may be considered as a constant and the above Definition may be drastically simplified).
The presence of the $\eta$-exponents follows in a less straightforward manner from the fact that one integrates against the {\em  $\eta$-approximation}
of the power kernels $K^{',\pm}$, $K^{*,\pm}$ (see the proof of Theorem \ref{th:2:exponents} for a computational explanation).

\begin{Theorem}[action of the kernels $K^{*,\pm}(\eta)$ of the admissible analytic functions]

{\it
Assume $f$ is an admissible analytic function with exponents as in Definition \ref{def:2:admissible}.
Then $g:=K^{*,\sigma}_{[0,t]}(\eta)f$ $(\sigma=\pm 1)$ is admissible, with $b$-exponents:
\begin{itemize}
\item[$(i)_b$]  
\BEQ
 \{ (B^{(b)}_j+U^{(b)}_j+2\alpha+1,0)_{j=1\ldots J}, (B^{(b)}_j, U^{(b)}_j+2\alpha+1)_{j=1\ldots J}  \}  
\nonumber\\    \label{eq:2:rules-ib}
\EEQ
on  the domain $(i)_b$;
\item[$(iii)_b$]
\BEA &&  \{(B^{(b)}_j+U^{(b)}_j+1,2\alpha)_{j=1\ldots J}\} \nonumber\\
&& \cup \{(B^{(b)''}_j,U^{(b)''}_j+2\alpha+1)_{j=1\ldots J''}  \}    \nonumber\\   \label{eq:2:rules-iiib}
\EEA
on  the domain $(iii)_b$,

\end{itemize}
 while its $\eta$-exponents are:

\begin{itemize}

\item[$(ii)_{\eta}$]
\BEA
&&\{(0,B_j^{(b)},U_j^{(b)}+2\alpha+1)_{j=1\ldots J},(U_j^{(b)}+2\alpha+1,B_j^{(b)},0)_{j=1\ldots J} \} \cup
\{(0,B_j^{(b)''},0)_{j=1\ldots J''} \} \nonumber\\
&& \cup \{ (H_j^{(\eta)'}+U_j^{(\eta)'}+2\alpha+1,B_j^{(\eta)'},0)_{j=1\ldots J'}, (H_j^{(\eta)''},B_j^{(\eta)''},0)_{j=1\ldots J''} \} \nonumber\\
&& \cup \{ (\tilde{H}_j^{(\eta)''},\tilde{B}_j^{(\eta)''},0)_{j=1\ldots J''},
  (\tilde{H}_j^{(\eta)''}+\tilde{U}_j^{(\eta)''}+1,\tilde{B}_j^{(\eta)''},0)_{j=1\ldots J''} \}
\label{eq:2:rules-iieta}
\EEA

\item[$(iii)_{\eta}$]

\BEA
&& \{(U^{(b)}_j+1,B^{(b)}_j,2\alpha)_{j=1\ldots J} \}\cup \{(0,B^{(b)''}_j,0)_{j=1\ldots
J''} \} \nonumber\\
&& \cup \{(H^{(\eta)''}_j+U^{(\eta)''}_j+1,B^{(\eta)''}_j,2\alpha)_{j=1\ldots J} \}  \nonumber\\
&& \cup
\{ (H^{(\eta)''}_j,B^{(\eta)''}_j,U^{(\eta)''}_j+2\alpha+1)_{j=1\ldots J''}, 
  (H^{(\eta)''}_j,B^{(\eta)''}_j,0)_{j=1\ldots J''} \} \nonumber\\
&& \cup  \{(\tilde{H}^{(\eta)''}_j+\tilde{U}^{(\eta)''}_j+1,\tilde{B}^{(\eta)''}_j,0)_{j=1\ldots J''},
(\tilde{H}^{(\eta)''}_j,\tilde{B}_j^{(\eta)''},0)_{j=1\ldots J''} \} \nonumber\\
\label{eq:2:rules-iiieta}
\EEA

\item[$\wit{(iii)}_{\eta}$]

\BEA
&& \{(U^{(b)}_j+1,B^{(b)}_j,0)_{j=1\ldots J}\} \cup \{(0,B^{(b)''}_j,0)_{j=1\ldots
J''}, (0,B_j^{(b)''},2\alpha+1)_{j=1\ldots J''}\} \nonumber\\
&& \cup \{(H^{(\eta)''}_j+U^{(\eta)''}_j+1,B^{(\eta)''}_j,0)_{j=1\ldots J''},
  (H^{(\eta)''}_j,B^{(\eta)''}_j,0)_{j=1\ldots J''}  \} \nonumber\\
&& \cup \{(\tilde{H}^{(\eta)''}_j+\tilde{U}^{(\eta)''}_j+1,\tilde{B}^{(\eta)''}_j,2\alpha)_{j=1\ldots J''} \} \nonumber\\
&&  \cup
\{ (\tilde{H}^{(\eta)''}_j,\tilde{B}^{(\eta)''}_j,\tilde{U}^{(\eta)''}_j+2\alpha+1)_{j=1\ldots J''},
 (\tilde{H}^{(\eta)''}_j,\tilde{B}^{(\eta)''}_j,0)_{j=1\ldots J''} \} \nonumber\\
\label{eq:2:rules-iiitildeeta}
\EEA

\item[$\wit{(ii)}_{\eta}$]

\BEA
&& \{ (U_j^{(b)}+1,B_j^{(b)},0)_{j=1\ldots J} \}
 \cup \{ (0,B_j^{(b)''},0)_{j=1\ldots J''},(0,B_j^{(b)''},2\alpha+1)_{j=1\ldots J''} \} \nonumber\\
&& \cup \{(H^{(\eta)''}_j,B^{(\eta)''}_j,0)_{j=1\ldots J''},(H_j^{(\eta)''}+U_j^{(\eta)''}+1,B^{(\eta)''}_j,0)_{j=1\ldots J''} \} \nonumber\\
&& \cup   \{ (\tilde{H}_j^{(\eta)'}+\tilde{U}_j^{(\eta)'}+2\alpha+1,\tilde{B}_j^{(\eta)'},0)_{j=1\ldots J'}, (\tilde{H}_j^{(\eta)''},\tilde{B}_j^{(\eta)''},0)_{j=1\ldots J''} \} \nonumber\\
\label{eq:2:rules-iitildeeta}
\EEA

\end{itemize}

Furthermore, the integer $\sigma_g$ appearing in Definition \ref{def:2:admissible} is $\sigma_g=\sigma_f+\sigma$.
}

\label{th:2:exponents}
\end{Theorem}

\begin{Theorem}[action of the kernels $K^{',\pm}(\eta)$ of the admissible analytic functions]

{\it
(same hypotheses). Add the following stronger assumption on the exponents: $U^{(b)}_j>-2\alpha$ $(j=1,\ldots,J)$.
  Then $g:=K^{',\pm}_{[0,t]} f$ is admissible, with exponents given
 as in Theorem
\ref{th:2:exponents} except that $\alpha$ should be replaced everywhere by $\alpha-1$.
}

\label{th:2:exponentsbis}
\end{Theorem}

{\bf Proof of Theorems \ref{th:2:exponents} and \ref{th:2:exponentsbis}.}

The proof is exactly the same for both Theorems. Note that the assumption on the exponents $(U_j)$, $(V_j)$ in
the case of Theorem \ref{th:2:exponentsbis} ensures that the exponents $(U^{(b)}_j)$, $(V^{(b)'''}_j)$ of 
$g=K^{',\pm}(\eta)f$
are larger than $-1$. No such assumption is needed in the case of Theorem \ref{th:2:exponents}
. We shall prove for instance Theorem \ref{th:2:exponents} and assume that $\bar{\sigma}=-1$ (otherwise one should take
the complex conjugate of all expressions).

By Lemma \ref{lemma:2:lemma1}, $K^{*,\sigma}_{[0,t]}(\eta)f$ is analytic on $\Omega$. The question is: what are its
exponents ?

Let us first look  at the contribution of $f_{b}$. We shall
use the variable $z:=a+\II(\sigma+\sigma_f)\eta$ throughout the proof. By splitting the interval of integration
$[0,t]$ into three pieces (corresponding to the domains $(i)_b$, $(ii)_b$, $(iii)_b$), there appears spurious singularities
near the inner boundaries of each subdomain, which cancel when summing up all terms. Hence we shall
not consider these boundary regions.

\begin{itemize}
\item[$(i)_b$]

Let 
\BEA
&&h(a) = \int_0^{\frac{b-\II\eta}{2}} (\II\sigma(u-a)+\eta)^{2\alpha}  (b-\II\eta)^B (u+\II\sigma_f \eta)^U
 F(\frac{u+\II\sigma_f \eta}{b-\II\eta},\frac{\eta}{b-\II\eta})\ du\nonumber\\
&&=  (b-\II\eta)^{B+U+2\alpha+1} \int_{\frac{\II\sigma_f \eta}{b-\II\eta}}^{\half \frac{b+\II\eta(2\sigma_f-1)}{b-\II\eta} } 
\left( \II\sigma(w-\frac{z}{b-\II\eta})  \right)^{2\alpha} w^U F(w,\frac{\eta}{b-\II\eta})\ dw
\nonumber\\
\EEA
(with $w=\frac{u+\II\sigma_f \eta}{b-\II\eta}$)
for some  $U>-1$. It is the sum of two terms (see  Lemma \ref{lemma:2:lemma5}). We may assume $|b/\eta|>c>0$ (otherwise apply simply Lemmas \ref{lemma:2:lemma1}, \ref{lemma:2:lemma2}),
so $\left| \half \frac{b+\II\eta(2\sigma_f-1)}{b-\II\eta} \right|> \left| \frac{\II\sigma_f \eta}{b-\II\eta} \right|$.
Apply now Lemma \ref{lemma:2:lemma5}, with the 'crude' version of case (iv) -- see last Remark after Definition
\ref{def:2:admissible} on the non-optimality of the $(f_b,f_{\eta})$-splitting --.
Then $h(a)=:h_b(z)+h_{\eta}(z)$ with :

\medskip

$\bullet$ ($b$-exponents)

\medskip  
One finds 
\BEA
 h_b(z) &=&   
  (b-\II\eta)^{B+U+2\alpha+1} G(\frac{z}{b-\II\eta},\frac{\eta}{b-\II\eta}) \nonumber\\
&+&
(b-\II\eta)^B z^{U+2\alpha+1} H(\frac{z}{b-\II\eta}, \frac{\eta}{b-\II\eta})   \nonumber\\
\EEA
if $\left| \frac{z}{\half[b+\II\eta(2\sigma_f-1)]} \right|<c<1$, and
\BEQ h_b(z)= (b-\II\eta)^{B+U+1} z^{2\alpha} H(\frac{b-\II\eta}{z},\frac{\eta}{b-\II\eta}) \EEQ
 if  $\left| \frac{z}{\half[b+\II\eta(2\sigma_f-1)]} \right|>C>1.$
 In other
words, one has the following sets of $b$-exponents : $\{(B+U+2\alpha+1,0),(B,U+2\alpha+1)\}$ on
 the domain $(i)_b$, and $\{(B+U+1,2\alpha)\}$ on the
domain $(iii)_b$.

\medskip

$\bullet$ ($\eta$-exponents)

\medskip

Suppose $\sigma_f\not=0$ (otherwise $h_{\eta}\equiv 0$).
One finds (on the cut domain $\Omega$)
\BEA
h_{\eta}(z)
&=&  (b-\II\eta)^{B+U+2\alpha+1}   \left[
  \left(\frac{z}{b-\II\eta}\right)^{U+2\alpha+1} F(\frac{ z}{ \eta},\frac{\eta}{b-\II\eta}) \right. \nonumber\\
&+ & \left.
\left(\frac{ \eta}{b-\II\eta}\right)^{U+2\alpha+1}
  G(\frac{ z}{ \eta},\frac{\eta}{b-\II\eta}) \right] \nonumber\\
\EEA
if $0<|\frac{z}{\II \eta}|<c<1$, and
 \BEQ h_b(z)= (b-\II\eta)^{B+U+2\alpha+1} \left( \frac{\II\sigma_f \eta}{b-\II\eta}\right)^{U+1}
 \left(\frac{z}{b-\II\eta}\right)^{2\alpha} F(\frac{ \eta}{z},\frac{\eta}{b-\II\eta})
\EEQ
 if $|\frac{z}{\II \eta}|>C>1$.
The function $F$ may go to infinity when $z\to \II\sigma_f \eta$, i.e. $a\to -\II\sigma\eta$, on the boundary
of the cut domain $\Omega$.  In other
words, one has the following sets of $\eta$-exponents :  $\{(0,B,U+2\alpha+1),(U+2\alpha+1,B,0)\}$
on the domain $(ii)_{\eta}$; $\{(U+1,B,2\alpha)\}$ on the
domain $(iii)_{\eta}$ and $\{(U+1,B,0)\}$ on the relative tilded domains $\wit{(iii)}_{\eta}$, $\wit{(ii)}_{\eta}$.

\item[$(ii)_b$]

Let 
\BEA
h(a)&=& \int_{\frac{b-\II\eta}{2}}^{b-\II\eta/2}
 ( \II\sigma(u-a)+\eta)^{2\alpha}   (b-\II\eta)^{B'} F(\frac{u+\II\sigma_f \eta}{b-\II\eta},\frac{\eta}{b-\II\eta})\ du \nonumber\\
&=&  (b-\II\eta)^{B'+2\alpha+1} \int_{\half \frac{b+\II\eta(2\sigma_f-1)}{b-\II\eta} }^{ \frac{b+\II\eta(\sigma_f-1/2)}{b-\II\eta} }
  (\II\sigma(w-\frac{z}{b-\II\eta}))^{2\alpha} F(w,\frac{\eta}{b-\II\eta})\ dw \nonumber\\
\EEA
with the same change of variables. 
Applying Lemma \ref{lemma:2:lemma5} (or Lemma \ref{lemma:2:lemma1})  to the above integral, one gets:
\BEQ h(a)= (b-\II\eta)^{B'+2\alpha+1} G(\frac{z}{b-\II\eta},\frac{\eta}{b-\II\eta}) \EEQ
if $\left|\frac{z}{\half[b+\II\eta(2\sigma_f-1)]}\right|<c<1$
 and 
\BEQ h(a)= (b-\II\eta)^{B'+1} z^{2\alpha} G(\frac{b-\II\eta}{z},\frac{\eta}{b-\II\eta}) \EEQ
 if $\left|\frac{z}{b+\II\eta(\sigma_f-1/2)]}\right|>C>1$. Hence $h$ is admissible of $b$-type, with
the following sets of $b$-exponents: $\{(B'+2\alpha+1,0)\}$ on the domain $(i)_b$,
$\{(B'+1,2\alpha)\}$ on the domain $(iii)_b$.

\item[$(iii)_b$]

Let 
\BEA
&& h(a)= \int_{b-\II\eta/2}^{t} ( \II\sigma(u-a)+\eta)^{2\alpha} 
 (b-\II\eta)^{B''} (u+\II\sigma_f \eta)^{U''} F(\frac{b-\II\eta}{u+\II\sigma_f \eta},
\frac{\eta}{b-\II\eta})\ du \nonumber\\
&&=  (b-\II\eta)^{B''+U''+1} z^{2\alpha}  \nonumber\\
&& \quad \int_{ \frac{b-\II\eta}{t+\II\sigma_f \eta}}^{ \frac{b-\II\eta}{b+\II\eta(\sigma_f-1/2)}}
 (\II\sigma(\frac{(b-\II\eta)}{z}-w))^{2\alpha}
 w^{-2\alpha-2-U''} F(w,\frac{\eta}{b-\II\eta})\ dw \nonumber\\
\EEA
(set $w=\frac{b-\II\eta}{u+\II\sigma_f \eta}$). By Lemma \ref{lemma:2:lemma5}, it is the sum of two terms,
$h=h_b+h_{\eta}$. One gets:

$\bullet$ ($b$-exponents)

\BEQ h_b(a)=(b-\II\eta)^{B''+U''+2\alpha+1}
 F(\frac{z}{b-\II\eta},\frac{\eta}{b-\II\eta}) \EEQ
if $\left| \frac{z}{b+\II\eta(\sigma_f-\half)} \right|<c<1$;
\BEA
 h_b(a)&=&(b-\II\eta)^{B''} z^{U''+2\alpha+1} F(\frac{b-\II\eta}{z},\frac{\eta}{b-\II\eta})
 \nonumber\\
&+&(b-\II\eta)^{B''+U''+1} z^{2\alpha} G(\frac{b-\II\eta}{z},
\frac{\eta}{b-\II\eta})
\EEA
if   $\left| \frac{z}{b+\II\eta(\sigma_f-\half)} \right|>C>1$.

In other words, $h_b$ is admissible of $b$-type, with  the following sets of $b$-exponents: $\{(B''+U''+2\alpha+1,0)\}$
 on the domain $(i)_b$, and 
$\{ (B'',U''+2\alpha+1), (B''+U''+1,2\alpha) \}$ on the domain $(iii)_b$.

$\bullet$ ($\eta$-exponents)

One finds:

\BEQ h_{\eta}(a)=(b-\II\eta)^{B''} H(\frac{b-\II\eta}{t+\II\sigma_f \eta},\frac{z}{t+\II\sigma_f \eta},\frac{\eta}{b-\II\eta}) \EEQ
if $\left| \frac{z}{t+\II\sigma_f \eta} \right|<c<1$;

\BEA
&&  h_{\eta}(a)=(b-\II\eta)^{B''} \left( G(\frac{b-\II\eta}{t+\II\sigma_f \eta},1-\frac{t+\II\sigma_f \eta}{z})
\right. \nonumber\\
&& \left. 
+(1-\frac{t+\II\sigma_f \eta}{z})^{2\alpha+1} H(\frac{b-\II\eta}{t+\II\sigma_f \eta},1-\frac{t+\II\sigma_f \eta}{z}) \right) \EEA
if $|1-z/t|<c<1.$

Hence $h_{\eta}$ is admissible of $\eta$-type, with the following sets of $\eta$-exponents: $\{(0,B'',0)\}$
on $(ii)_{\eta}$ and $(iii)_{\eta}$; $\{ (0,B'',0),(0,B'',2\alpha+1)\}$ on $\widetilde{(iii)}_{\eta}$,
$\widetilde{(ii)}_{\eta}$. 

\end{itemize}

\bigskip

There remains to analyze the contribution to $K^{*,\sigma}_{[0,t]}(\eta)f$ of $f_{\eta}$. These terms are simpler since 
(as we shall presently see) they only
contribute admissible functions of $\eta$-type.

\begin{itemize}

\item[$(ii)_{\eta}$]

Consider
\BEA h(a)&=& \int_0^{2\eta} (\II\sigma(u-a)+\eta)^{2\alpha} \eta^{H'} (b-\II\eta)^{B'}  (u+\II\sigma_f \eta)^{U'}
\nonumber\\
&& \quad 
F(\frac{u+\II\sigma_f \eta}{\sigma_f \eta},\frac{\eta}{b-\II\eta},\frac{b-\II\eta}{t}) \ du \nonumber\\
&=& (b-\II\eta)^{B'} \eta^{H'+U'+2\alpha+1} \int_{\II\sigma_f}^{2(1+\II\sigma_f/2)} (\II\sigma(w-\frac{z}{\eta}))^{2\alpha}  w^{U'} 
F(w,\frac{\eta}{b-\II\eta},\frac{b-\II\eta}{t}) \ dw
\nonumber\\
\EEA
where one has set $w=\frac{u+\II\sigma_f \eta}{\eta}$.
 By  Lemma \ref{lemma:2:lemma2}
\BEQ h(a)=(b-\II\eta)^{B'}  \eta^{H'+U'+2\alpha+1}  G(\frac{z}{\eta},\frac{\eta}{b-\II\eta},\frac{b-\II\eta}{t})
 \EEQ
if $|\frac{z}{2(1+\II\sigma_f/2)\eta}|<c<1$,  where $G$ becomes possibly infinite when $z\to \II\sigma_f \eta$;
\BEQ h(a)=(b-\II\eta)^{B'} \eta^{H'+U'+1} z^{2\alpha} H(\frac{\eta}{z},\frac{\eta}{b-\II\eta},\frac{b-\II\eta}{t}) \EEQ
if $|\frac{z}{2(1+\II\sigma_f/2)\eta}|>C>1$,

 whence $h$ is of $\eta$-type with $\eta$-exponents $\{(H'+U'+2\alpha+1,B',0\}$
on $(ii)_{\eta}$, 
 $\{(H'+U'+1,B',2\alpha)\}$ on $(iii)_{\eta}$, and
$\{(H'+U'+1,B',0)\}$ on the relative tilded  domains $\wit{(iii)}_{\eta}$, $\wit{(ii)}_{\eta}$.

\medskip

\item[$(iii)_{\eta}$]

Consider
\BEA && h(a) = 
\int_{2\eta}^{t/2} (\II\sigma(u-a)+\eta)^{2\alpha} \eta^{H''} (b-\II\eta)^{B''}
(u+\II\sigma_f \eta)^{U''} \nonumber\\
&& \left(  F(\frac{\eta}{u+\II\sigma_f \eta},\frac{\eta}{b-\II\eta},
\frac{b-\II\eta}{t}) + G(\frac{u+\II\sigma_f \eta}{t+\II\sigma_f \eta},\frac{\eta}{b-\II\eta},\frac{b-\II\eta}{t}) \right)
\ du \nonumber\\
&&= \eta^{H''+U''+1} (b-\II\eta)^{B''} z^{2\alpha} \int_{\frac{2\eta}{t+2\II\sigma_f \eta}}^{\frac{1}{1+\II\sigma_f/2}} (\II\sigma(\frac{\eta}{z}-w))^{2\alpha}
w^{-2\alpha-2-U''} F(w,\frac{\eta}{b-\II\eta},\frac{b-\II\eta}{t})\ dw   \nonumber\\
&& +\eta^{H''} (b-\II\eta)^{B''} \int_{\frac{(2+\II\sigma_f)\eta}{t+\II\sigma_f \eta}}^{\frac{\frac{t}{2}+\II\sigma_f \eta}{t+\II
\sigma_f \eta}} (\II\sigma(w-\frac{z}{t+\II\sigma_f \eta}))^{2\alpha}  w^{U''} G(w,\frac{\eta}{b-\II\eta},\frac{b-\II\eta}{t})\ dw
\nonumber\\
&&=: h_1(a)+h_2(a)
\EEA
with $w=\frac{ \eta}{u+\II\sigma_f \eta}$, resp. $\frac{u+\II\sigma_f \eta}{t+\II\sigma_f \eta}$. 

Using
Lemma \ref{lemma:2:lemma5} -- with the 'refined' version of case (iv) --, one obtains:
\BEA
 h_1(a) &=& \eta^{H''} (b-\II\eta)^{B''} H(\frac{z}{t},\frac{\eta}{b-\II\eta},\frac{b-\II\eta}{t}) \nonumber\\
 &+& \eta^{H''+U''+2\alpha+1} (b-\II\eta)^{B''} 
 \tilde{G}(\frac{z}{\eta},\frac{\eta}{b-\II\eta},\frac{b-\II\eta}{t}) \EEA
if  $|\frac{z}{(1+\II\sigma_f/2)\eta}|<c<1$; 
\BEA
 h_1(a)&=& \eta^{H''} (b-\II\eta)^{B''} z^{U''+2\alpha+1} F(\frac{\eta}{z},\frac{\eta}{b-\II\eta},\frac{b-\II\eta}{t}) \nonumber\\
&+&
\eta^{H''+U''+1} (b-\II\eta)^{B''}
z^{2\alpha} G(\frac{\eta}{z},\frac{\eta}{b-\II\eta},\frac{b-\II\eta}{t})  \nonumber\\
&+&  \eta^{H''} (b-\II\eta)^{B''}
H( \frac{z}{t}, \frac{\eta}{b-\II\eta},\frac{b-\II\eta}{t})  \EEA
if $\left| \frac{z}{\frac{t}{2}+\II\sigma_f \eta} \right|<c<1$,
 $\left| \frac{z}{(1+\II\sigma_f/2)\eta}\right|>C>1$, and
\BEA
 h(a)&=& \eta^{H''+U''+1} (b-\II\eta)^{B''}  G(\frac{\eta}{z},\frac{\eta}{b-\II\eta},\frac{b-\II\eta}{t})
\nonumber\\
&+&  \eta^{H''} (b-\II\eta)^{B''} F(1-\frac{z}{t},\frac{\eta}{b-\II\eta},\frac{b-\II\eta}{t}) \EEA
if $\left| \frac{z}{\frac{t}{2}+\II\sigma_f \eta} \right|>C>1$, $|z/t|<2$. (Note that $\frac{\eta}{z}=\frac{\eta}{t}
\frac{1}{1-(1-z/t)}$ is an analytic function of $\frac{\eta}{t}=\frac{\eta}{b-\II\eta} \frac{b-\II\eta}{t}$ and 
$1-\frac{z}{t}$, while the powers in $z$ in factor in the last equation may be skipped since $\left(\frac{z}{t}\right)^{\beta}=(1-(1-\frac{z}{t}))^{\beta}$ is an analytic function of $\frac{z}{t}$).

\bigskip

Similarly,

\BEA
&& h_2(a)=\eta^{H''} (b-\II\eta)^{B''} G(\frac{z}{t},\frac{\eta}{b-\II\eta},\frac{b-\II\eta}{t}) \nonumber\\
&& + \eta^{H''+U''+2\alpha+1} (b-\II\eta)^{B''} H(\frac{z}{\eta},\frac{\eta}{b-\II\eta},\frac{b-\II\eta}{t}) \nonumber\\
\EEA
if $\left| \frac{z}{(2+\II\sigma_f)\eta}\right|<c<1$;
\BEA
&& h_2(a)=\eta^{H''}(b-\II\eta)^{B''}z^{U''+2\alpha+1} F(\frac{z}{t},\frac{\eta}{b-\II\eta},\frac{b-\II\eta}{t})\nonumber\\
&& +\eta^{H''} (b-\II\eta)^{B''} G(\frac{z}{t},\frac{\eta}{b-\II\eta},\frac{b-\II\eta}{t}) \nonumber\\
&&+\eta^{H''+U''+1} (b-\II\eta)^{B''} z^{2\alpha} H(\frac{\eta}{z},\frac{\eta}{b-\II\eta},\frac{b-\II\eta}{t}) \nonumber\\
\EEA
if  $\left| \frac{z}{\frac{t}{2}+\II\sigma_f \eta} \right|<c<1$,
 $\left| \frac{z}{(1+\II\sigma_f/2)\eta}\right|>C>1$ ;

\BEA
&& h_2(a)=\eta^{H''} (b-\II\eta)^{B''} \tilde{G}(1-\frac{z}{t},\frac{\eta}{b-\II\eta},\frac{b-\II\eta}{t})\nonumber\\
&&+\eta^{H''+U''+1} (b-\II\eta)^{B''} H(\frac{\eta}{z},\frac{\eta}{b-\II\eta},\frac{b-\II\eta}{t}) \nonumber\\
\EEA

if  $\left| \frac{z}{\frac{t}{2}+\II\sigma_f \eta} \right|>C>1$, $|z/t|<2$.

Hence one gets the following $\eta$-exponents: $\{(H'',B'',0),(H''+U''+2\alpha+1,B'',0)\}$
on $(ii)_{\eta}$; 
  $\{(H'',B'',U''+2\alpha+1),(H''+U''+1,B'',2\alpha),(H'',B'',0)\}$
on $(iii)_{\eta}$;  $\{(H''+U''+1,B'',0),(H'',B'',0)\}$ on $\wit{(iii)}_{\eta}$ and $\wit{(ii)}_{\eta}$.

\item[$\wit{(ii)}_{\eta}$,$\wit{(iii)}_{\eta}$]

The expansions may be derived very simply from the previous ones by using the symmetry $u\longleftrightarrow t-u$. 

\end{itemize}

\hfill \eop

\begin{Corollary}[integrability and $b$-exponents]

{\it
\begin{enumerate}
\item (integrability)
Suppose $f$ is an admissible analytic function of $b$-type  with exponents $(U^{(b)}_j)_{j=1\ldots J}$
 given by $U^{(b)}_1=0$,\ $U^{(b)}_2=2\alpha-1$. Then $(K_{[0,t]}^{',\pm}(\eta)K_{[0,t]}^{*,\pm}(\eta))^m f$, resp. $K_{[0,t]}^{*,\pm}(\eta) (K_{[0,t]}^{',\pm}(\eta)K_{[0,t]}^{*,\pm}(\eta))^m f$ are
well-defined admissible analytic functions for any $m=0,1,\ldots$, with exponents $\{(U_j)\}\subset
\{4\alpha n,4\alpha n+2\alpha-1\ |\ n=0,\ldots,m \} $,  resp.
 $\{(U^{(b)}_j)\}\subset \{0\}\cup
\{4\alpha n+2\alpha+1,4\alpha n+4\alpha\ |\ n=0,\ldots,m \} $.  In particular, the exponents $(U^{(b)}_j)_{j=1\ldots J}$  of\\
 $K_{[0,t]}^{*,\pm}(\eta) (K_{[0,t]}^{',\pm}(\eta)K_{[0,t]}^{*,\pm}(\eta))^m f$, resp.
$(K_{[0,t]}^{',\pm}(\eta)K_{[0,t]}^{*,\pm}(\eta))^m f$
 are all non-negative, resp. $\ge 2\alpha-1$.

Similarly, if $f$ is an admissible analytic function of $b$-type with exponents $(U^{(b)}_j)_{j=1\ldots J}$
 given by $U^{(b)}_1=0$,\ $U^{(b)}_2=2\alpha+1$ then $(K_{[0,t]}^{*,\pm}(\eta) K_{[0,t]}^{',\pm}(\eta))^m f$,
resp. $K_{[0,t]}^{',\pm}(\eta)  (K_{[0,t]}^{*,\pm}(\eta) K_{[0,t]}^{',\pm}(\eta))^m f$, are well-defined admissible analytic functions for any $m\ge 0$,
and the exponents  $(U^{(b)}_j)_{j=1\ldots J}$ of $ (K_{[0,t]}^{*,\pm}(\eta)K_{[0,t]}^{',\pm}(\eta))^m f$, resp.
$K_{[0,t]}^{',\pm}(\eta)  (K_{[0,t]}^{*,\pm}(\eta)K_{[0,t]}^{',\pm}(\eta))^m f$  are all non-negative, resp. $\ge 2\alpha-1$.

\item (on some other $b$-exponents)
Let $n\ge 1$. Suppose $f$ is an admissible analytic function of $b$-type with exponents $(B^{(b)}_j,U^{(b)}_j)_{j=1\ldots J}$ such that
$B^{(b)}_j+U^{(b)}_j=4\alpha n-1$, resp. $4\alpha n+2\alpha$  for all $j$ and
 $\{(B^{(b)''}_j)\}=\{0,4\alpha(n-1)+2\alpha+1\}$, resp.
$\{0,4\alpha n\}$. Then 
the $B^{(b)''}$-exponents of  $K_{[0,t]}^{*,\pm}(\eta) (K_{[0,t]}^{',\pm}(\eta)K_{[0,t]}^{*,\pm}(\eta))^m f$, resp.
$(K_{[0,t]}^{*,\pm}(\eta) K_{[0,t]}^{',\pm}(\eta))^m f$  $(m\ge 0)$ are all non-negative, and
so are all the elements of the set  $\{B^{(b)}_j+U^{(b)}_j,j=1\ldots J\}=\{B^{(b)''}_j+U^{(b)''}_j,j=1\ldots J''\}$.

 Also,
the $B^{(b)''}$-exponents of $(K_{[0,t]}^{',\pm}(\eta)K_{[0,t]}^{*,\pm}(\eta))^m f$, resp. $K_{[0,t]}^{',\pm}(\eta) (K_{[0,t]}^{*,\pm}(\eta)K_{[0,t]}^{',\pm}(\eta))^m f$ are all non-negative, and
the elements of the set $\{B^{(b)}_j+U^{(b)}_j,j=1\ldots J\}=\{B^{(b)''}_j+U^{(b)''}_j,j=1\ldots J''\}$ are all $\ge 2\alpha-1$.

\item More precisely, if one of these sums of exponents ($X_j$, say) satisfies instead a strict inequlity, namely,
$X_j>0$, resp. $X_j>2\alpha-1$, then $X_j\ge 2\alpha$, resp. $X_j\ge 0$.

\end{enumerate}
}
\label{cor:2:integrability}
\end{Corollary}

{\bf Proof.}

\begin{enumerate}

\item Elementary induction on $n$ using Theorems \ref{th:2:exponents} and \ref{th:2:exponentsbis}.
\item The inequalities hold true for $f$. Then one may prove  that if the $b$-exponents of any admissible function $f$
(of $b$-type) satisfy the relations $B^{(b)''}_j,B^{(b)''}_j+U^{(b)''}_j \ge 0$, then the $b$-exponents
of $K^{',\pm}_{[0,t]}(\eta)f$ satisfy the corresponding relations $B^{(b)''}_j\ge 0$, $B^{(b)''}_j+U^{(b)''}_j
\ge 
2\alpha-1$, and vice versa if one considers $K^{*,\pm}_{[0,t]}(\eta)f$ instead.
\item is proved along the same lines as 2.
\end{enumerate}

 \hfill \eop

\begin{Corollary}[$\eta$-exponents]

{\it 
Let $f$ be  an admissible analytic function of $b$-type with exponents $(B^{(b)}_j,U^{(b)}_j)_{j=1\ldots J}$ such that
$B^{(b)}_j+U^{(b)}_j=4\alpha n-1$, resp. $4\alpha n+2\alpha$  for all $j$ and
 $\{(B^{(b)''}_j)\}=\{0,4\alpha(n-1)+2\alpha+1\}$, resp.
$\{0,4\alpha n\}$ (see Corollary \ref{cor:2:integrability}, point 2). Then:

\begin{enumerate}
\item
Let $g_0:=K^{*,\pm}_{[0,t]}(\eta) (K^{',\pm}(\eta)_{[0,t]} K^{*,\pm}_{[0,t]}(\eta))^m f$, resp. $(K^{*,\pm}_{[0,t]}(\eta) K^{',\pm}_{[0,t]}(\eta))^m f$ for
some $m\ge 0$, $g_1:=\left( K^{',\pm}_{[0,t]}(\eta) \left( (g_0)_b \right) \right)_{\eta}$ be the $\eta$-part of the integral of
the $b$-part of $g_0$ against $K^{',\pm}(\eta)$, and $h:= (K^{*,\pm}_{[0,t]}(\eta) K^{',\pm}_{[0,t]}(\eta) )^n K^{*,\pm}_{[0,t]}(\eta) g_1$, resp. $(K^{',\pm}_{[0,t]}(\eta) K^{*,\pm}_{[0,t]}(\eta))^n g_1$ for some $n\ge 0$. 
Then $h$ is of $\eta$-type, with $\eta$-exponents such that
\BEQ  H^{(\eta)''}_j,\tilde{H}^{(\eta)''}_j,  H^{(\eta)''}_j+B^{(\eta)''}_j,\tilde{H}^{(\eta)''}_j+\tilde{B}^{(\eta)''}_j
\ge 0  \label{eq:2:eta-exponents-condition1} \EEQ and
\BEQ U_j^{(\eta)'}, \tilde{U}_j^{(\eta)'}, H^{(\eta)''}_j+B^{(\eta)''}_j+U^{(\eta)''}_j,\tilde{H}^{(\eta)''}_j+\tilde{B}^{(\eta)''}_j+\tilde{U}^{(\eta)''}_j,
H^{(\eta)''}_j+U^{(\eta)''}_j,\tilde{H}^{(\eta)''}_j+\tilde{U}^{(\eta)''}_j\ge 0 \label{eq:2:eta-exponents-condition2}, \EEQ
resp. (depending on $h$)
\BEQ U_j^{(\eta)'}, \tilde{U}_j^{(\eta)'}, H^{(\eta)''}_j+B^{(\eta)''}_j+U^{(\eta)''}_j,\tilde{H}^{(\eta)''}_j+\tilde{B}^{(\eta)''}_j+\tilde{U}^{(\eta)''}_j,
H^{(\eta)''}_j+U^{(\eta)''}_j,\tilde{H}^{(\eta)''}_j+\tilde{U}^{(\eta)''}_j\ge 2\alpha-1 \label{eq:2:eta-exponents-condition2bis}. \EEQ

More precisely, if one of these sums of exponents ($X_j$, say) satisfies instead a strict inequlity, namely,
$X_j>0$, resp. $X_j>2\alpha-1$, then $X_j\ge 2\alpha$, resp. $X_j\ge 0$.

\item
Similarly, let $g_0:= (K^{',\pm}(\eta)_{[0,t]} K^{*,\pm}_{[0,t]}(\eta))^m f$, resp. $K^{',\pm}_{[0,t]}(\eta) (K^{*,\pm}_{[0,t]}(\eta) K^{',\pm}_{[0,t]}(\eta))^m f$ for
some $m\ge 0$, $g_1:=\left( K^{*,\pm}_{[0,t]}(\eta) \left( (g_0)_b \right) \right)_{\eta}$ be the $\eta$-part of the integral of
the $b$-part of $g_0$ against $K^{',\pm}(\eta)$, and $h:= (K^{*,\pm}_{[0,t]}(\eta) K^{',\pm}_{[0,t]}(\eta) )^n  g_1$,
resp. $K^{',\pm}_{[0,t]}(\eta)(K^{*,\pm}_{[0,t]}(\eta) K^{',\pm}_{[0,t]}(\eta))^n g_1$  for some $n\ge 0$. 
Then $h$ is of $\eta$-type, with $\eta$-exponents $(H^{(\eta)''}_j,B^{(\eta)''}_j,U^{(\eta)''}_j)$ and
 $(\tilde{H}^{(\eta)''}_j,\tilde{B}^{(\eta)''}_j,\tilde{U}^{(\eta)''}_j)$ satisfying the same relations.
\end{enumerate}
}
\label{cor:2:eta-exponents}
\end{Corollary}

{\bf Proof.}

\begin{enumerate}
\item

$\bullet$ Let $(B_{g_0}^{(b)},U_{g_0}^{(b)})$  be the $(B^{(b)},U^{(b)})$-exponent
 of any term in $(g_0)_b$. By Corollary \ref{cor:2:integrability},
 $U_{g_0}^{(b)}$, $U_{g_0}^{(b)}+B_{g_0}^{(b)},B_{g_0}^{(b)''}\ge 0$.
Now rules (\ref{eq:2:rules-iieta}), (\ref{eq:2:rules-iitildeeta}),
 (\ref{eq:2:rules-iiieta}) and (\ref{eq:2:rules-iiitildeeta}) in Theorem \ref{th:2:exponents} imply that
 the $\eta$-exponents of $g_1$ satisfy relations (\ref{eq:2:eta-exponents-condition1}) and (\ref{eq:2:eta-exponents-condition2bis}). Now one may check very easily that $K^{*,\pm}_{[0,t]}(\eta)f_{\eta}$ satisfies relations 
 (\ref{eq:2:eta-exponents-condition1}) and (\ref{eq:2:eta-exponents-condition2}) if $f_{\eta}$ satisfies relations
 (\ref{eq:2:eta-exponents-condition1}), (\ref{eq:2:eta-exponents-condition2bis}) and vice-versa if one considers
$K^{',\pm}_{[0,t]}(\eta)f_{\eta}$ instead.

\item

The proof is similar, with initial relations $U_{g_0}^{(b)},B_{g_0}^{(b)}+U_{g_0}^{(b)}\ge 2\alpha-1$ and $B_{g_0}^{(b)''}\ge 0$ this time.

\end{enumerate}
\hfill \eop

{\bf Remark.} Note that the exponents $U_j^{(\eta)'}$ and $\tilde{U}_j^{(\eta)'}$ are simply $0$, except those which come
from an admissible function of $b$-type (and are actually only due to the $(f_b,f_{\eta})$-splitting, which is somewhat
unfortunate in this respect).


\subsection{Asymptotic behaviour of the moments of the L\'evy area}


Now the computation of the exponents is over, we may study the singularities of $\esper[({\cal A}_{s,t}(\eta))^{2N}]$
when $\eta\to 0$. It turns out eventually that only one term of the $2N$-th connected moment is singular
(see Theorem \ref{th:2:divergence-moments} below, cited in the Introduction).  This term comes from the only closed bipartite diagram
with alternating simple and double lines, see end of subsection 1.1, and obtained by iterating $I_+$-type
integrals, see comments at the end of subsection 1.2.

Elementary arguments relying on Lemma \ref{lemma:1:exp} allow then to deduce the asymptotic behaviour
of $\esper[ ({\cal A}_{s,t}(\eta))^{2N}]$ from that of the connected moments.

The proof of Theorem \ref{th:2:divergence-moments} requires first a separate analysis of all terms coming
from the splitting of the kernel $K$ (see end of subsection 1.1).

\begin{Lemma}

{\it
 Let $I_n(\eta,b,t;u)=(K^*_{[s,t]}(\eta)K'_{[s,t]}(\eta))^n \left( u\mapsto (\II\bar{\sigma}
(u-b)+\eta)^{2\alpha} \right)(u)$, $\bar{\sigma}\in\{\pm 1\}$.
 Then $I_n$ 
is the sum of an admissible analytic function of type $b$  with $b$-exponents
such that $U^{(b)}_j,B^{(b)''}_j,B^{(b)}_j+U^{(b)}_j\ge 0$ for all possible indices $j$, of an admissible analytic function
of type $\eta$ with $\eta$-exponents such that
\BEQ  H_j^{(\eta)''},H_j^{(\eta)''}+U_j^{(\eta)''}, H_j^{(\eta)''}+B_j^{(\eta)''},H_j^{(\eta)''}+B_j^{(\eta)''}
+U_j^{(\eta)''}, U_j^{(\eta)'}\ge 0   \label{eq:2:HBU>0} \EEQ
(plus the same inequalities with a tilde),  \\
 and of the non-analytic function
\BEQ R_n(b,u)=C_n\Re (\II(b-u)+(2n+1)\eta)^{2\alpha+4\alpha n}\EEQ
with \BEQ C_n=\frac{1}{2\pi} \left( \frac{\pi/2}{\cos\pi\alpha \Gamma(-2\alpha)} \right)^{2n} \sin \pi\alpha \ 
\Gamma(2\alpha+1) \Gamma(-2\alpha-4\alpha n).\EEQ
}

\label{lemma:2:term-Fig2}
\end{Lemma}

{\bf Remark.} As appears clearly  in the proof below, the non-analytic function $R_N$
doesn't show up when one considers the moments of the  L\'evy area of the analytic process $\Gamma$ (see
Proposition \ref{prop:1}), namely
$\int_0^t d\Gamma^{(1)}_x(\eta) \int_0^x d\Gamma^{(2)}_y(\eta).$

\medskip

{\bf Proof.}

By definition (choosing for instance $\bar{\sigma}=-1$), 
\BEA
&&  I_n=\sum_{\sigma_{2(N-n)},\ldots,\sigma_{2N}\in \{\pm 1\}} K^{*,\sigma_{2(N-n)}}_{[0,t]}(\eta) 
K^{',\sigma_{2(N-n)+1}}_{[0,t]}(\eta) \cdots
K^{*,\sigma_{2N-2}}_{[0,t]}(\eta) K^{',\sigma_{2N-1}}_{[0,t]} (\eta) \nonumber\\
&& 
\left( u\mapsto(-\II(u-(b-\II\eta)))^{2\alpha} \right) (u). \EEA

 Iterating the non-analytic term $(\pm\II(b-a))^{\beta_1+\beta_2+1}$ appearing
in Lemma \ref{lemma:2:lemma3bis}, one obtains (up to a certain coefficient) $(\pm\II(u-b)+(2n+1)\eta)^{2\alpha+4\alpha m}$
after $2m$ integrations, resp.  $(\pm\II(u-b)+(2m+2)\eta)^{4\alpha-1+4\alpha m}$ after $2m+1$ integrations.  This is possible
only if the  $2m$, resp.  $2m+1$ iterated integrals are of $I_+$-type, i.e. if $-1=\sigma_{2N}=\ldots=\sigma_{2(N-m)}$,
resp. $-1=\sigma_{2N}=\ldots=\sigma_{2(N-m)-1}$. 
Integrate once again and look instead at the analytic term this time: by  Definition \ref{def:2:admissible}, it is
 an admissible analytic function with exponents $\{(B_j,U_j)\}=\{(4\alpha-1+4\alpha m,0),
(2\alpha+4\alpha m,2\alpha-1)\}$,
 $\{(B''_j,U''_j)\}=\{(0,4\alpha-1+4\alpha m),(2\alpha+1+4\alpha m,2\alpha-2)\}$,
 resp. $\{(B_j,U_j)\}=\{(6\alpha+4\alpha m,0),
(4\alpha-1+4\alpha m,2\alpha+1)\}$,
 $\{(B''_j,U''_j)\}=\{(0,6\alpha+4\alpha m)),(4\alpha+4\alpha m,2\alpha)\}$, possibly
up to the symmetry $u\to t-u$ (so $b$ must be replaced with $t-b$ in that case).
Now Theorem \ref{th:2:exponents}  shows that integrating such an admissible function
alternatively against the kernels $K^{*,\pm}(\eta)$ and $K^{',\pm}(\eta)$ yields admissible functions. The lower
bounds on the exponents come from Corollaries \ref{cor:2:integrability} and \ref{cor:2:eta-exponents} in the 'good' case
where exponents are non-negative.

Suppose now {\em all} signs $\sigma_{2N},\ldots,\sigma_{2(N-n)}$ are equal. Then one  obtains in the end  a non-analytic term
by iterating $2n$ times Lemma \ref{lemma:2:lemma3bis}. The coefficient before that term may be checked
by an easy induction using the complement formula for the Gamma function, namely,
$\Gamma(x)\Gamma(1-x)=\frac{\pi}{\sin\pi x}$.

\hfill \eop

\vskip 1 cm

We may now proceed to estimate the connected diagrams.
Let us first analyze the contribution of the {\em open diagrams} (see end of subsection 1.1). The following Theorem
shows that they are all regular in the limit $\eta\to 0$.

\begin{Theorem}

{\it The contribution of the open diagrams writes $C_N t^{4n\alpha}+O(\eta^{2\alpha})$ for some constant $C_N$.}

\label{th:2:open-diagrams}

\end{Theorem}

{\bf Proof.}

 Recall  that the open diagrams are products of terms
of three types, $(\emptyset \emptyset)$, $(\emptyset \newmoon)$ and $(\newmoon \newmoon)$. Let us analyze these
three cases separately.

\bigskip

\underline{Case $(\newmoon \newmoon)$}. 

Let 
\BEA 
&&  I^{\newmoon \newmoon}_t := \sum_{\sigma_0,\ldots,\sigma_{2n-1}} 
\int_0^t \ldots \int_0^t (\II\sigma_0 x_1+\eta)^{2\alpha} (\II\sigma_1(x_1-x_2)+\eta)^{2\alpha-2} \ldots \nonumber\\
&&  (\II\sigma_{2n-1}(x_{2n-1}-x_{2n})+\eta)^{2\alpha-2}  (\II\sigma_{2n} x_{2n}+\eta)^{2\alpha} dx_1\ldots
dx_{2n} \EEA
with $\sigma_0,\ldots,\sigma_{2n}\in\{\pm 1\}$. Then $I_t^{\newmoon\newmoon}=I_n(\eta,0,t;0)$
in the notation of Lemma \ref{lemma:2:term-Fig2}. The non-analytic part $R_n(0,0)$ of $I_n$ is negligible, of
order $O(\eta^{2\alpha+4\alpha n})$, so let us consider the analytic part. Assume first $\sigma_f\not=0$, then
$\eta^H (b+\II\bar{\sigma}\eta)^B (u+\II\sigma_f \eta)^U=C\eta^{H+B+U}$ with
\BEA && (H,B,U)\in \{ (0,B_j^{(b)},U_j^{(b)})_j, (0,B_j^{(b)'},0)_j, (0,B_j^{(b)''},U_j^{(b)''})_j,  \nonumber\\
&&  (H_j^{(\eta)'},B_j^{(\eta)'},U_j^{(\eta)'})_j, (H_j^{(\eta)''},B_j^{(\eta)''},U_j^{(\eta)''})_j\}.
\EEA

In any case, by Lemma \ref{lemma:2:term-Fig2}, $I_n(\eta,0,t;0)$ converges when $\eta\to 0$. More
precisely,  Corollaries \ref{cor:2:integrability} and \ref{cor:2:eta-exponents} imply:
$I_n(\eta,0,t;0)-I_n(0,0,t;0)=O(\eta^{2\alpha}).$

If now $\sigma_f=0$, then $\eta^H (b+\II\bar{\sigma}\eta)^B (u+\II\sigma_f \eta)^U=0$ except if $U=0$, with $U=U_j^{(b)}$
or $U_j^{(\eta)'}$ (by Lemma \ref{lemma:2:term-Fig2}, these exponents are non-negative). The rest of the proof is identical.

\bigskip

\underline{Case $(\emptyset \newmoon)$}.

Let 
\BEA
&&  I^{\emptyset \newmoon}_t:=\sum_{\sigma_2,\ldots,\sigma_{2n-1}} \int_0^t \ldots \int_0^t  (\II\sigma_1(x_1-x_2)+\eta)^{2\alpha-2} \ldots \nonumber\\
&&  (\II\sigma_{2n-1}(x_{2n-1}-x_{2n})+\eta)^{2\alpha-2} (\II\sigma_{2n} x_{2n}+\eta)^{2\alpha} dx_1\ldots
dx_{2n} \EEA
with $\sigma_1,\ldots,\sigma_{2n}\in\{\pm 1\}$. Then
\BEQ I_t^{\emptyset \newmoon}=\int_0^t (K^{',\sigma_1}_{[0,t]}(\eta) I_{n-1}(\eta,0,t;.))(u)\ du.\EEQ
The non-analytic part of $(K^{',\sigma_1}_{[0,t]}(\eta) I_{n-1}(\eta,0,t;.)(u)$ writes (up to a coefficient)
$(\pm\II u+2n\eta)^{4\alpha-1+4\alpha(n-1)}$ which is uniformly integrable in $\eta$; the integral
writes $C+O(\eta^{4\alpha n})+O(\eta)$ $(n\ge 1)$. As for the analytic part, one must integrate
$(b+\II\bar{\sigma}\eta)^B (u+\II\sigma_f \eta)^U=C\eta^B (u+\II\sigma_f \eta)^U$ (with the $b$-exponents), resp.
$\eta^H (b+\II\bar{\sigma}\eta)^B (u+\II\sigma_f \eta)^U=C \eta^{H+B} (u+\II\sigma_f \eta)^U$ (with the $\eta$-exponents).
The integral over $(i)_B: 0<|u+\II\sigma_f \eta|<\eta/3$, resp. $(iii)_b:2\eta<|u+\II\sigma_f \eta|$ yields an 
expression bounded by $O(\eta^{B_j^{(b)}+U_j^{(b)}+1})$, resp. $O(\eta^{B_j^{(b)''}})+O(\eta^{B_j^{(b)''}+U_j^{(b)'''}+1})$.
Similar statements hold for the $\eta$-exponents (one must essentially replace $B$ with $H+B$). Now use the relations
given in Corollaries \ref{cor:2:integrability} and \ref{cor:2:eta-exponents} in the 'bad' case where exponents may
be negative.

\bigskip

\underline{Case $(\emptyset \emptyset)$}. 

Let
\BEA
&&   I^{\emptyset \emptyset}_t:= \nonumber\\
&& \sum_{\sigma_2,\sigma_{2n-2}} \int_0^t \ldots \int_0^t  (\II\sigma_1(x_1-x_2)+\eta)^{2\alpha-2} \ldots (\II\sigma_{2n-1}(x_{2n-1}-x_{2n})+\eta)^{2\alpha-2}  dx_1\ldots
dx_{2n} \nonumber\\
&&= C \int_0^t dx_{2n-1} \left(  \int_0^t J_{n-2}(\eta,x_{2n-1},t;x_1)\ dx_1\right) \nonumber\\
&& \left[ (\II\sigma_{2n-1} (x_{2n-1}-t)+\eta)^{2\alpha-1}-(\II\sigma_{2n-1}x_{2n-1}+\eta)^{2\alpha-1} \right] 
\EEA
($\sigma_1,\ldots,\sigma_{2n-1}\in\{\pm 1\}$), where ($n\ge 2$)
\BEQ J_{n-2}(\eta,x_{2n-1},t;x_1)=\left( K^{',\sigma_1}_{[0,t]}(\eta) I_{n-2}(\eta,x_{2n-1},t;.)\right)(x_1). \EEQ

Set $b=x_{2n-1}$ for convenience. The non-analytic part of $J_{n-2}$ writes (up to a coefficient) $(\pm\II(x_1-b)+2(n-1)\eta)^{4\alpha-1+4\alpha(n-2)}$, hence its contribution to $I_t^{\emptyset \emptyset}$ is of the form
\BEA &&  C_1 \int_0^t db\ (\pm\II b+2(n-1)\eta)^{4\alpha(n-1)} (\pm\II b+\eta)^{2\alpha-1} \nonumber\\
 && \quad +C_2 \int_0^t db\ (\pm\II(t-b)+2
(n-1)\eta)^{4\alpha-1} (\pm\II b+\eta)^{2\alpha-1} 
\EEA
plus two similar terms (that reduce to the previous ones by the symmetry $b\leftrightarrow t-b$). By splitting the integral
into $\int_0^{\eta}\ db+\int_{\eta}^{t/2}\ db+\int_{t/2}^{t-\eta}\ db+\int_{t-\eta}^t\ db$, it is easy to prove
that the limit when $\eta\to 0$ writes $C+O(\eta^{2\alpha})$.

Turning to the analytic part, one must integrate in $u:=x_1$ separately over each of the 7 domains $(i)_b,\ldots,\tilde{(ii)}_{\eta}$. The precise dependence in $\eta$ may be computed by rewriting the proof of Theorem \ref{th:2:exponents} with $\alpha=0$, which makes things essentially trivial. Let us write in details for instance the case $u\in(i)_b$. One has

\BEA && \int_0^{\frac{b-\II\eta}{2}} J_{n-1}(\eta,b,t;u)\ du=(b-\II\eta)^{B_j^{(b)}+U_j^{(b)}+1} \int_{\frac{\II\sigma_f \eta}{b-\II\eta}}^{\half\frac{b+\II\eta(2\sigma_f-1)}{b-\II\eta}} w^{U_j^{(b)}} F(w,\frac{\eta}{b-\II\eta})\ dw \nonumber\\
&&= (b-\II\eta)^{B_j^{(b)}+U_j^{(b)}+1} \left( G(\frac{\eta}{b-\II\eta}+\left(\frac{\eta}{b-\II\eta}\right)^{U_j^{(b)}+1}
H(\frac{\eta}{b-\II\eta}) \right). 
\EEA

Now integrate in $b$ the product of this function by $(\pm\II(t-b)+\eta)^{2\alpha-1}-(\mp \II b+\eta)^{2\alpha-1}$ by
the same method as for the non-analytic term.

The other domains are left to the reader.
\hfill \eop

\vskip 1 cm

We now turn to the contribution of the unique {\em closed diagram}  $\int_0^t du K'(\eta,b,u) 
I_{N-1}(\eta,b,t;u).$ The function
$I_{N-1}$ has been analyzed in Lemma  \ref{lemma:2:term-Fig2}. There only 
remains to integrate it  against the kernel $K^{',\pm}(\eta)$.

\begin{Lemma}

{\it 
The double integral
\BEQ F(\eta;t):=\int_0^t db \int_0^t du\ (\pm\II(u-b)+\eta)^{2\alpha-2} I_{N-1}(\eta,b,t;u)  \EEQ
tends to a finite limit $F(t)$ when $\eta\to 0$. More precisely, $F(\eta;t)=F_0(t)+O(\eta^{2\alpha})$ where $F_0$
is independent of $\eta$.
}

\label{lemma:2:admissible-double-integral}
\end{Lemma}

{\bf Proof.}

Integrating with respect to $u$ gives (applying Theorem \ref{th:2:exponentsbis} with $u=b$)
\BEQ F(\eta;t)=\int_0^t db\ g_b(\eta,b,t)+\int_0^t db\ g_{\eta}(\eta,b,t)\EEQ
where
\BEQ g_b(\eta,b,t)=\sum_{j=1}^{J'} (b-\II\eta)^{B_j^{(b)'}} G_j(\frac{\eta}{b-\II\eta}),\EEQ
and
\BEQ g_{\eta}(\eta,b,t)=\sum_{j=1}^{J'} \eta^{H_j^{(\eta)'}} (b-\II\eta)^{B_j^{(\eta)'} +U_j^{(\eta)'}} H_j(\frac{\eta}{b-\II\eta}, \frac{b-\II\eta}{t}) \EEQ

if $|b/\eta|<C$,

\BEQ g_{\eta}(\eta,b,t)=\sum_{j=1}^{J''} \eta^{H_j^{(\eta)''}}
 (b-\II\eta)^{B_j^{(\eta)''}+U_j^{(\eta)''}} F_j(\frac{\eta}{b-\II\eta},\frac{b-\II\eta}{t}) \EEQ
if $|b/\eta|>C>0$, $|1-b/t|>C'>0$, and similarly for the domains with a tilde. Adequate lower bounds for the exponents are
given in Corollaries \ref{cor:2:integrability} and \ref{cor:2:eta-exponents}. 

There remains to integrate over $b$. Consider for instance $g_b$. One has
\BEA && \int_0^t db\ (b-\II\eta)^B G(\frac{\eta}{b-\II\eta})=\int_0^{\eta} db\ (b-\II\eta)^B G(\frac{\eta}{b-\II\eta})
+\int_{\eta}^t db\ (b-\II\eta)^B G(\frac{\eta}{b-\II\eta}) \nonumber\\
&& =C \eta^{B+1} (1+O(\eta))+C'(1+O(\eta)) \EEA
(use for instance a series expansion for $G$ in the second integral). In the present case, $B=B_j^{(b)'}=B_j^{(b)}+U_j^{(b)}\ge 2\alpha-1$. The other cases are similar. \hfill \eop

\bigskip

We may now prove the main result of this section.

\newpage

\begin{Theorem}

{\it
The $2N$-th connected moment of the $\eta$-approximation of the L\'evy area $\phi^{(c)}_{2N}(\eta;t)$ is given
by the sum of two terms: the first one is regular in the limit $\eta\to 0$ and equal to $C_{reg,N} t^{4N\alpha}+O(\eta^{2\alpha})$ for
some constant $C_{reg,N}$; the second one is equal to $C_{irr,N} t \eta^{4N\alpha-1}$ with
\BEQ C_{irr,N}=\left( \frac{\pi/2}{\cos \pi\alpha \Gamma(-2\alpha)} \right)^{2(N-1)} \sin \pi\alpha \frac{\Gamma(2\alpha+1)}{\Gamma(2-2\alpha)} \Gamma(1-4\alpha N) (2N)^{4\alpha N-1}.\EEQ
}
\label{th:2:divergence-moments}
\end{Theorem}

{\bf Proof.}

The first term is obtained by summing the contribution of all admissible functions, see Theorem \ref{th:2:open-diagrams} and
Lemma \ref{lemma:2:admissible-double-integral}. The second one is obtained from the single irregular term (see
Lemma \ref{lemma:2:term-Fig2})
\BEQ \frac{\alpha(1-2\alpha)}{2\cos\pi\alpha} \int_0^t du \int_0^t db (\pm\II(u-b)+\eta)^{2\alpha-2} R_{N-1}(b,u). \EEQ

Up to a coefficient (essentially $C_{N_1}$, see Lemma \ref{lemma:2:term-Fig2}), and forgetting the regular $I_-$ integral,
this is
 \BEQ 2\Re \int_0^t db\ I_+(2\alpha-2,2\alpha+4\alpha(N-1);0,t)(b+\II\eta,b-\II(2N-1)\eta),\EEQ
whose irregular part is (see Lemma \ref{lemma:2:lemma3bis})
\BEQ -2 \frac{\Gamma(2\alpha+1+4\alpha(N-1))\Gamma(1-4\alpha N)}{\Gamma(2-2\alpha)} \ . \ 2\sin\pi[2\alpha+4\alpha(N-1)] \ .
(2N\eta)^{4\alpha N-1},\EEQ
hence the result.

\hfill \eop

{\bf Remark.} Stirling's formula implies: $|C_{irr,N}|\le C'^N$ for some constant $C'$.

\begin{Corollary}

{\it
The $2N$-th moment of the $\eta$-approximation of the L\'evy area $\phi_{2N}(\eta;t)$ writes
\BEQ \esper[({\cal A}_{0,t}(\eta))^{2N}]=(2N-1)!! C_{irr,1}^N t^N \eta^{(4\alpha-1)N} (1+O(\eta^{1-4\alpha})).\EEQ
}
\label{cor:2:2Nth-moment}
\end{Corollary}

{\bf Proof.}

Decompose $\esper[({\cal A}_{0,t}(\eta))^{2N}]$ into a sum of diagrams as in subsection 1.1. The
$(2N-1)!!$ totally disconnected diagrams (with $N$ trivial components) contribute the main term in
the expansion, namely
\BEQ (2N-1)!! \left( C_{irr,1} t \eta^{4\alpha-1}+O(1) \right)^N=(2N-1)!! C_{irr,1}^N t^N \eta^{(4\alpha-1)N}
(1+O(\eta^{1-4\alpha})).\EEQ
Then 'almost' totally disconnected diagrams (i.e. with $N-2$ trivial components and only one bipartite
polynomial line with $4$ legs) contribute 
\BEQ c_N \left( C_{irr,1} t\eta^{4\alpha-1}+O(1)\right)^{N-2} \left( C_{irr,2} t\eta^{8\alpha-1} +O(1) \right)
\sim_{\eta\to 0} c'_N t^{N-1} \eta^{(4\alpha-1)N+1} \EEQ
for some constants $c_N,c'_N$. It is clear that the other diagrams are even less irregular. \hfill \eop

{\bf Remark.} 
The coefficient $(2N-1)!!$ in front of $(t\eta^{4\alpha-1})^N$ is equal to $\esper[X^{2N}]$ if $X$ is a standard
Gaussian variable. This simple remark leads after some more computations to Theorem A, our main result (see
Introduction)  which states
that the rescaled L\'evy area converges to Brownian motion.


\vskip 3cm

We shall also need the following refinement of Theorem \ref{th:2:divergence-moments} in section 2 to establish the
existence of a uniform exponential moment for the rescaled L\'evy area.

\begin{Theorem}

{\it
Let $t\in [0,T]$, $T$ fixed. Then there exists a constant $C$ (depending only on $\alpha$) such that the regular
part $C_{reg,N} t^{4N\alpha}+O(\eta^{2\alpha})$ of the $2N$-th connected moment $\phi_{2N}^{(c)}(\eta;t)$ satisfies:
$|\phi_{2N}^{(c)}(\eta;t)|\le C^N$ for all $N\ge 1$.
}

\label{th:2:refinement}
\end{Theorem}

The proof relies mainly on the following Lemma.

\begin{Lemma}

{\it 
Assume $t\in [0,T]$ for some fixed $T>0$, and $\eta>0$ small enough. Let also
\BEA
 && \Omega':=\{u\in\C |\ |u/t|<2\}\setminus  \nonumber\\
&& \quad \left( \{\frac{u}{\eta}\in-\II\sigma+\R_-\} \cup \{\frac{u-t}{\eta}\in
-\II\sigma+\R_+\}\cup B(-\II\sigma\eta,\frac{\eta}{3})\cup B(t-\II\sigma\eta,\frac{\eta}{3}) \right). \nonumber\\
\EEA

 Let $f\in C([0,t])$ be a continuous function such that $f(z)=z^{\beta_1}\phi_1(z)$ if $z\in[0,t/2]$, $z^{\beta_2}\phi_2(z)$
if $z\in[t/2,t]$, where $\beta_1,\beta_2\ge 0$, and $\phi_1,\phi_2$ are analytic and bounded on $\Omega'$.
Finally, let $\sigma_1,\sigma_2\in\{\pm 1\}$. Then
$ \left(K^{*,\sigma_1}_{[0,t]}(\eta) K^{',\sigma_2}_{[0,t]}(\eta) \right)f$ is analytic in $\Omega'$ and furthermore
\BEQ \sup_{z\in\Omega'} | \left(K^{*,\sigma_1}_{[0,t]}(\eta) K^{',\sigma_2}_{[0,t]}(\eta) \right)f|\le C 
\sup_{z\in\Omega'}
|f(z)| \EEQ
where $C$ depends only on  $\alpha$.

}

\label{lemma:2:refinement}
\end{Lemma}

Note that $\Omega'$ is essentially the domain $\Omega_{res}$ defined in the first Remark after Definition 
\ref{def:2:admissible}.

{\bf Proof.}

Consider first the easier case $\sigma_1=-\sigma_2$: let for instance
\BEA
g(z)&:=&\left(K^{*,-}_{[0,t]}(\eta) K^{',+}_{[0,t]}(\eta) f\right)(z)=\int_0^t (-\II(z-u)+\eta)^{2\alpha} (K^{',+}_{[0,t]}(\eta)f)(u)\ du \nonumber\\
&=& C' \int_0^t  (-\II(z-u)+\eta)^{2\alpha} du \int_0^t (\II(u-v)+\eta)^{2\alpha-2}\ f(v)\ dv \nonumber\\
&=& C'\int_0^t f(v)\ dv \ .\ \int_0^t  (\II(u-z)+\eta)^{2\alpha} (\II(u-v))^{2\alpha-2}\ d\bar{w} \nonumber\\
&=& \int_0^t f(v)\ dv \ .\ \overline{I_-(2\alpha-2,2\alpha;0,t)(v+\II\eta,z+\II\eta)} \quad (z\in\Pi^+)
\EEA

The $I_-$-function extends analytically to $\Omega'$; it is bounded when its arguments $v+\II\eta$, $z+\II\eta$ are bounded
away from $0$ and $t$, and has an integrable singularity (with negative exponents $2\alpha-1$, $4\alpha-1$, see Lemma
\ref{lemma:2:lemma3} and the Remark following it) near $0$ and $t$. Hence $|g(z)|\le C\sup_{z\in[0,t]} |f(z)|$.

\bigskip

Consider now, for $z\in \Pi^+$,
\BEA
g(z)&:=&(K^{*,-}_{[0,t]}(\eta) K^{',-}_{[0,t]}(\eta) f)(z)=\int_0^t (-\II(z-u)+\eta)^{2\alpha} (K^{',-}_{[0,t]}(\eta)f)(u)\ du \nonumber\\
&=& C' \int_{\Gamma_1} (-\II(z-w)+\eta)^{2\alpha} dw \int_0^t (-\II(w-v)+\eta)^{2\alpha-2}\ f(v)\ dv \nonumber\\
&=& C'\int_0^t f(v)\ dv \ .\ \int_{\Gamma_1} (\II(w-z)+\eta)^{2\alpha} (-\II(w-v)+\eta)^{2\alpha-2}\ dw \nonumber\\
&=& C' \int_0^t f(v)\ dv \ .\ \overline{I_+(2\alpha-2,2\alpha;0,t)(v-\II\eta,z+\II\eta)}
\EEA
where $\Gamma_1:0\to t$, $\Gamma((0,t))\subset\Pi^+$ and $\Gamma$ passes below $z$.

The same conclusions hold as in the previous case except (see Lemma \ref{lemma:2:lemma3bis}) for the supplementary
non-analytic term of the form $C (\II(v-z))^{4\alpha-1}$. This term may be integrated against $f$ as in Lemma
\ref{lemma:2:lemma1}, yielding for $z\in\Pi^-$ an integral
\BEQ \int_{\bar{\Gamma}_2} f(\bar{w}) (\II(\bar{w}-z))^{4\alpha-1}\ d\bar{w},\EEQ
where $\bar{\Gamma}_2$ passes below $z$ while staying in $\Omega$ when it leaves the real axis. Now the singularity
of the kernel $(\II(\bar{w}-z))^{4\alpha-1}$ is integrable, hence the result. \hfill \eop

\medskip

We may now prove briefly Theorem \ref{th:2:refinement}. Consider, as in the proof of Lemma \ref{lemma:2:term-Fig2},
a $(2m)$- or $(2m+1)$-iterated non-analytic term, and then the analytic term obtained by integrating once more this
non-analytic term. Integrate against $K^*$ once more if the last integration was against $K'$, so the result, $f$,
satisfies the hypotheses of Lemma \ref{lemma:2:refinement}. Namely, it has positive exponents $U_j^{(b)}, U_j^{(\eta)'},
\tilde{U}_j^{(\eta)'}$ by Lemma \ref{lemma:2:term-Fig2}, hence it is bounded on $\Omega'$ (see first remark
after Definition \ref{def:2:admissible}), with an overall constant $|C_{irr,m}|\le C'^m$ by the remark following
Theorem
\ref{th:2:divergence-moments}. One may now iterate Lemma \ref{lemma:2:refinement} till the last two iterated integrals.

 \hfill \eop


\section{Convergence in law of the rescaled L\'evy area}


Using the analysis of singularities developed in Section 1, we mainly aim to prove in this section the following Theorem.

\begin{Definition}

{\it
Let $\tilde{\cal A}_{s,t}(\eta):=\eta^{\half(1-4\alpha)} {\cal A}_{s,t}(\eta)$ be the rescaled L\'evy area.
}

\label{def:3:rescaled}
\end{Definition}

{\bf Theorem A}

{\it The three-dimensional  process $(B^{(1)}(\eta),B^{(2)}(\eta),\tilde{\cal A}(\eta))$  converges in law in the Skohorod
topology to
 $(B^{(1)},B^{(2)},\sqrt{C_{irr,1}} \del W)$
where $\del W_{s,t}:=W_t-W_s$ are the increments of  a standard one-dimensional Brownian motion {\em
independent} from $B^{(1)}$ and $B^{(2)}$.
}

\bigskip

Theorem A is a consequence of the following Theorem which generalizes the asymptotics obtained in section 1
for the moments of ${\cal A}_{s,t}$ to the case of the moments of finite-dimensional distributions.

\begin{Theorem}

{\it
Let $(W_t)_{t\in \R}$ be a (two-sided) standard one-dimensional Brownian motion. Then, for every $s_1<t_1,\ldots,
 s_{n}<t_{n}$, and \\ $u_{(1),1},\ldots, u_{(1),k_1},u_{(2),1},\ldots,u_{(2),k_2}\in\R$,
\BEA
&& \esper \left[\left( B^{(1)}_{u_{(1),1}}(\eta)\ldots B^{(1)}_{u_{(1),k_1}}(\eta) \right)
\left( B^{(2)}_{u_{(2),1}}(\eta)\ldots B^{(2)}_{u_{(1),k_2}}(\eta) \right) 
\tilde{\cal A}_{s_1,t_1}(\eta) \ldots \tilde{\cal A}_{s_{n},t_{n}}(\eta) \right]\to_{\eta\to 0} \nonumber\\
&& 
C_{irr,1}^{n/2} \esper[(W_{t_1}-W_{s_1})\ldots (W_{t_{n}}-W_{s_{n}})].\ \esper[B^{(1)}_{u_{(1),1}}\ldots
B^{(1)}_{u_{(1),k_1}}]\ \esper[B^{(2)}_{u_{(2),1}}\ldots
B^{(2)}_{u_{(2),k_2}}]. \nonumber\\
\EEA

\label{th:3:convergence-moments}
}

\end{Theorem}

The fact that the moments of $\tilde{\cal A}(\eta)$ and $B(\eta)$ 'factorize' in the limit $\eta\to 0$ is of course an
indication of the asymptotic independence of $\del W$ and $B$.

For the proof of Theorem  \ref{th:3:convergence-moments}, we shall need the following two  Lemmas.

\begin{Lemma}

{\it
Let $s_1<t_1,s_2<t_2$, then
\BEQ \esper[ {\cal A}_{s_1,t_1}(\eta) {\cal A}_{s_2,t_2}(\eta)]=C_{irr,1}\  \lambda([s_1,t_1]\cap
[s_2,t_2])\  \eta^{4\alpha-1}+O(1) \EEQ
($\lambda$=Lebesgue measure) when $\eta\to 0$.
}
\label{lemma:3:convergence-moments1}
\end{Lemma}

{\bf Proof.}

Consider first the case when the intervals $[s_1,t_1]$ and $[s_2,t_2]$ are disjoint,
say, $s_1<t_1\le s_2<t_2$. Then
\BEA
&&\left| \esper[ {\cal A}_{s_1,t_1}(\eta)  {\cal A}_{s_2,t_2}(\eta)]  \right| \nonumber\\
&& \le C \int\int_{s_1\le x_1,y_1\le t_1} dx_1 \ dy_1 \int_{s_2\le x_2,y_2\le t_2} dx_2\ dy_2
(x_2-x_1)^{2\alpha-2} (y_2-y_1)^{2\alpha-2} \nonumber\\
&& = C  \left[ \int_{s_1\le x_1\le t_1, s_2\le x_2\le t_2} (x_2-x_1)^{2\alpha-2} \ dx_1\ dx_2 \right]^2
\nonumber\\
&&
= C' \left[ (t_2-t_1)^{2\alpha}+(s_2-s_1)^{2\alpha}-(t_2-s_1)^{2\alpha}-(s_2-t_1)^{2\alpha} \right]
<\infty. \nonumber\\
\EEA

For the general case, set $[s_1,t_1]\cap[s_2,t_2]:=[s,t]$ ($s<t$) so that $[s_i,t_i]=[s_i,s]\cup [s,t]\cup
[t,t_i]$ $(i=1,2)$ and $\lambda([s_1,t_1]\cap[s_2,t_2])=t-s$. Write (using the multiplicative property (\ref{eq:0:mult}))
\BEA
  {\cal A}_{s_i,t_i}(\eta)&=& {\cal A}_{s_i,s}(\eta)+{\cal A}_{s,t}(\eta)+{\cal A}_{t,t_i}(\eta) \nonumber\\
&+& (B_s^{(1)}(\eta)-B_{s_i}^{(1)}(\eta))(B^{(2)}_t(\eta)-B^{(2)}_{s}(\eta))+(B^{(1)}_{t}(\eta)-
B_{s_i}^{(1)}(\eta))(B_{t_i}^{(2)}(\eta)-B_t^{(2)}(\eta)).  \nonumber\\
\EEA

Forgetting about the products of increments appearing on the last line, the only singular term comes from
$\esper[({\cal A}_{s,t}(\eta))^2]=C_{irr,1} |t-s| \eta^{4\alpha-1}$. Now the covariance between the area terms and
the products of increments is regular in the limit $\eta\to 0$, as follows from the general arguments in
Lemma \ref{lemma:3:convergence-moments2} below (which does not use the result of this Lemma). \hfill \eop

\begin{Lemma}

{\it
Let $s_1<t_1,\ldots,s_{2N}<t_{2N}$, and \\ $u_{(1),1},\ldots, u_{(1),k_1},u_{(2),1},\ldots,u_{(2),k_2}\in\R$
as in Theorem \ref{th:3:convergence-moments}. Consider the closed connected diagrams coming from the evaluation
of 
\BEQ
\esper \left[\left( B^{(1)}_{u_{(1),1}}(\eta)\ldots B^{(1)}_{u_{(1),k_1}}(\eta) \right)
\left( B^{(2)}_{u_{(2),1}}(\eta)\ldots B^{(2)}_{u_{(1),k_2}}(\eta) \right) 
{\cal A}_{s_1,t_1}(\eta) \ldots {\cal A}_{s_{2N},t_{2N}}(\eta) \right].
\EEQ
 Then their sum is regular in the limit $\eta\to 0$ unless $k_1=k_2=0$. In the latter case, their sum  writes
$(2N-1)!\ C_{irr,N} t\eta^{4N\alpha-1}$ plus a regular term in the limit $\eta\to 0$, where $t=\lambda([s_1,t_1]\cap
\ldots \cap [s_{2N},t_{2N}]).$

The same expression with an {\em odd} number of ${\cal A}$'s is always regular in the limit $\eta\to 0$.
}

\label{lemma:3:convergence-moments2}
\end{Lemma}

{\bf Proof.} 

Set $s_{2N+1},\ldots,s_{2N+k_1+k_2}=0$ and
 $$(t_{2N+1},\ldots,t_{2N+k_1})=(u_{(1),1},\ldots, u_{(1),k_1}), \quad
(t_{2N+k_1+1},\ldots,t_{2N+k_1+k_2})=(u_{(2),1},\ldots, u_{(2),k_2}).$$
Decompose $[s_1,t_1]\cup\ldots\cup[s_{2N+k_1+k_2},t_{2N+k_1+k_2}]$ into a finite union of  intervals 
with disjoint interiors $I_1,\ldots, I_m$, with $m$ minimal,
so that for each $j$ and $k$, one has either $I_k\subset [s_j,t_j]$ or $I_k\cap [s_j,t_j]=\emptyset$.
Let $B_I:=B_t-B_s$ and ${\cal A}_I:={\cal A}_{s,t}$ if $I=[s,t]$.
Then the multiplicativity property for the area, see equation (\ref{eq:0:mult}), implies:

\BEA
&&  \esper \left[\left( B^{(1)}_{u_{(1),1}}(\eta)\ldots B^{(1)}_{u_{(1),k_1}}(\eta) \right)
\left( B^{(2)}_{u_{(2),1}}(\eta)\ldots B^{(2)}_{u_{(2),k_2}}(\eta) \right) 
{\cal A}_{s_1,t_1}(\eta) \ldots {\cal A}_{s_{2N},t_{2N}}(\eta) \right]
\nonumber\\
&& = 
\esper\left[\prod_{j=1}^{2N} \left( {\cal A}_{I_{l_{j,1}}}(\eta)+ {\cal A}_{I_{l_{j,2}}}(\eta)+\ldots+
P_j(s_1,\ldots,s_{2N+k_1+k_2},t_1,\ldots,t_{2N+k_1+k_2}) \right) \right. \nonumber\\
&& \left. \prod_{k=1}^{k_1} (B^{(1)}_{I_{l^{(1)}_{k,1}}}(\eta)+B^{(1)}_{I_{l^{(1)}_{k,2}}}(\eta)+\ldots)
\prod_{k=1}^{k_2} (B^{(2)}_{I_{l^{(2)}_{k,1}}}(\eta)+B^{(2)}_{I_{l^{(2)}_{k,2}}}(\eta)+\ldots)
  \right] \EEA
for some polynomials  $P_j$ in the variables \\
 $B_{s_1}(\eta),\ldots, B_{s_{2N+k_1+k_2}}(\eta),
B_{t_1}(\eta),\ldots,B_{t_{2N+k_1+k_2}}(\eta)$.   Expanding the product, one gets terms of the type
$\esper[ {\cal A}_{I_{i_1}}(\eta)\ldots {\cal A}_{I_{i_{k}}}(\eta) Q]$ for some polynomial $Q$ as above.

Assume to begin with that $Q=1$. Then the singularities of $\esper[{\cal A}_{I_{i_1}}(\eta)\ldots
{\cal A}_{I_{i_k}}(\eta)]$ may be investigated as in section 1, with the only difference that the
non-analytic terms come from the iterated integrals
\BEQ K^{*,-}_{I_{i'_1}}(\eta) K^{',-}_{I_{i'_2}}(\eta)\ldots K^{',-}_{I'_{i_{2k'}}}(\eta)(u\mapsto
(-\II(u-(b-\II\eta)))^{2\alpha})(u) \EEQ
or
\BEQ K^{',-}_{I_{i'_1}}(\eta) K^{*,-}_{I_{i'_2}}(\eta)\ldots K^{',-}_{I'_{i_{2k'+1}}}(\eta)(u\mapsto
(-\II(u-(b-\II\eta)))^{2\alpha})(u) \EEQ
(and their conjugates) for all possible choices of subsets of intervals $\{I_{i'_1},\ldots,I_{i'_{k}}\}\subset
\{I_1,\ldots,I_m\}$. Now  Remark \ref{rmk:2} after the proof of Lemma \ref{lemma:2:lemma3bis} proves that the
corresponding  non-analytic
term is $0$ except if $I_{i'_1}=\ldots=I_{i'_{k}}$. This gives the singular term $(2N-1)!\ C_{irr,N} t\eta^{4N\alpha-1}$
by Theorem \ref{th:2:divergence-moments}. 

 If now $Q$ is of degree $\ge 1$, then $Q$ may be written as $$\left(\int_{I_{i'_1}} dB_{s_1}^{(1)}(\eta)\int_{I_{j'_1}}
dB_{t_1}^{(2)}(\eta)\right) \ldots \left(\int_{I_{i'_{k}}} dB_{s_k}^{(1)}(\eta)\int_{I_{j'_{k}}}
dB_{t_k}^{(2)}(\eta)\right), \quad k\ge 1$$
 for some intervals $I_{i'_l},I_{j'_l}$ (indeed, the expectation is simply zero if there isn't the same
number of $B^{(1)}$'s and $B^{(2)}$'s).
Consider any connected diagram and evaluated it by the same method as usual (see subsection 1.1).
 The integration $\int_{I_{j'_1}} dt_1$ yields a result which
does not depend on $s_1$, hence the next integration $\int_{I_{i'_1}} ds_1$ contains no non-analytic
term. The results of subsection 1.4 show then that $\esper[{\cal A}_{I_{i_1}}(\eta)\ldots{\cal A}_{I_{i_k}}(\eta)Q]$
is regular in the limit $\eta\to 0$.

The proof is the same for an odd number of ${\cal A}$'s.

\hfill \eop

\bigskip

{\bf Proof of Theorem \ref{th:3:convergence-moments}}. 

Assume $n=2N$ (the proof is the same for $n$ odd).
Decompose $[s_1,t_1]\cup\ldots\cup[s_{2N+k_1+k_2},t_{2N+k_1+k_2}]$ into a finite union of disjoint intervals $I_1,\ldots, I_m$
as in the proof of Lemma \ref{lemma:3:convergence-moments2}. 
Then the multiplicativity property for the area  implies once again:
\BEA
&&  \esper \left[\left( B^{(1)}_{u_{(1),1}}(\eta)\ldots B^{(1)}_{u_{(1),k_1}}(\eta) \right)
\left( B^{(2)}_{u_{(2),1}}(\eta)\ldots B^{(2)}_{u_{(2),k_2}}(\eta) \right) 
\tilde{\cal A}_{s_1,t_1}(\eta) \ldots \tilde{\cal A}_{s_{2N},t_{2N}}(\eta) \right]
\nonumber\\
&& = 
\esper\left[\prod_{j=1}^{2N} \left( \tilde{\cal A}_{I_{l_{j,1}}}(\eta)+ \tilde{\cal A}_{I_{l_{j,2}}}(\eta)+\ldots+
\eta^{\half(1-4\alpha)}
P_j \right) \right. \nonumber\\
&& \left. \prod_{k=1}^{k_1} (B^{(1)}_{I_{l^{(1)}_{k,1}}}(\eta)+B^{(1)}_{I_{l^{(1)}_{k,2}}}(\eta)+\ldots)
\prod_{k=1}^{k_2} (B^{(2)}_{I_{l^{(2)}_{k,1}}}(\eta)+B^{(2)}_{I_{l^{(2)}_{k,2}}}(\eta)+\ldots)
  \right] \EEA
  Expanding the product, one gets terms of the type
\BEQ \esper[ {\cal A}_{I_{i_1}}(\eta)\ldots {\cal A}_{I_{i_{k}}}(\eta) R] 
\EEQ
 for some polynomial $R$ as above. Use a diagrammatic expansion now. Lemma \ref{lemma:3:convergence-moments2} shows that
connected diagrams involving both L\'evy areas and a non-trivial product of increments are regular, hence go to $0$ in
the limit $\eta\to 0$ because of the rescaling. This implies the 'factorization' of the moments of 
$\tilde{\cal A}(\eta)$ and $B(\eta)$
in the limit $\eta\to 0$. Now the same Lemma implies that (if $R=1$) 
 only the totally
disconnected diagrams survive in the limit $\eta\to 0$, still because of the rescaling
(see proof of Corollary \ref{cor:2:2Nth-moment} for the power-counting).
 Then Lemma \ref{lemma:3:convergence-moments1} allows one to evaluate such diagrams, and the result follows
 now from a simple combinatorial argument by 're-gluing' together the intervals $I_1,\ldots,I_m$. \hfill \eop

\bigskip

We may now prove Theorem A. 

\bigskip

{\bf Proof of Theorem A.}

The moments of the finite-dimensional distributions of $(B^{(1)}(\eta),B^{(2)}(\eta),\tilde{\cal A}(\eta))$
converge to those of $(B^{(1)},B^{(2)},\del W)$ as Theorem \ref{th:3:convergence-moments} shows. Furthermore,
\BEQ \esper[ ||B_t(\eta)-B_s(\eta)||^2]\le C|t-s|^{2\alpha} \EEQ
for all $\eta$ (which results from a simple computation using the explicit formula for
the covariance, or from \cite{Unt08}, Lemma 1.5) 
and (assuming for instance $s<t\le u<v$, but similar estimates hold in all cases), using
once again the multiplicative property (\ref{eq:0:mult})
\BEA 
&& \esper[ (\tilde{\cal A}_{s,t}(\eta)-\tilde{\cal A}_{u,v}(\eta))^2] \nonumber\\
&&=\esper[(\tilde{\cal A}_{s,u}(\eta)-\tilde{\cal A}_{t,v}(\eta)+\eta^{\half(1-4\alpha)}(\del B_{s,u}^{(1)}
\del B_{u,t}^{(2)}
-\del B_{u,t}^{(1)}\del B_{t,v}^{(2)}))^2] \nonumber\\
&& \le C\left( |u-s|+|v-t|+\eta^{1-4\alpha} |u-s|^{2\alpha} |v-t|^{2\alpha}\right) \nonumber\\
\EEA
(where $\del B_{x,y}:=B_y-B_x$), 
which implies that the sequence of processes is tight (by standard arguments for processes in a Gaussian chaos
of finite order, see \cite{Billingsley}). Since Gaussian laws are uniquely characterized by their moments, this implies Theorem A. \hfill \eop

\medskip

{\bf Remark.}  By using the central limit theorem due to D.Nualart and G. Peccati \cite{NuaPec05}, it would have been
enough to compute the second and fourth moments of the joint process to prove the convergence
in law. But we feel that these particular cases are not much
easier than the general case, and that the general and powerful asymptotic analysis given in Section 1 may be used
for a wide range of applications.

\bigskip

Let us end this paragraph by giving  an estimate of the characteristic function of the rescaled
 L\'evy area $\tilde{\cal A}_{s,t}(\eta)$ for $\alpha\in(\frac{1}{8},\frac{1}{4})$, implying 
a uniform exponential bound.

\begin{Lemma}

{\it 
Let $\alpha\in(\frac{1}{8},\frac{1}{4})$, $0<s,t<T$ ($T$ fixed) and
 $\tilde{\phi}_{s,t}(\eta;\lambda)=\esper[e^{\II\lambda \tilde{\cal A}_{s,t}(\eta)}]$ be the characteristic function of the L\'evy
area. Then  there exist  two constants $\lambda_0,C_0>0$ such that, for all $\lambda\in\C$ satisfying
 $|\lambda|\le \lambda_0 \eta^{-\half(1-4\alpha)}$, and for all $\eta$ small enough,
\BEQ |\tilde{\phi}_{s,t}(\eta;\lambda)-e^{-\half C_{irr,1}|t-s|\lambda^2}|\le C_0
\lambda^2 \eta^{1-4\alpha}  e^{-\half C_{irr,1}|t-s|\lambda^2}.\EEQ
}

\label{lemma:3:characteristic-function}
\end{Lemma}

{\bf Proof.}

Theorems \ref{th:2:divergence-moments} and  \ref{th:2:refinement} yield
\BEA
&&  \tilde{\phi}^{(c)}_{s,t}(\eta;\lambda)=\phi^{(c)}_{s,t}(\eta;\eta^{\half(1-4\alpha)}\lambda) \nonumber\\
&& =-\half \lambda^2 \left( C_{irr,1} |t-s| +C_{reg,1}(\eta;t-s)\eta^{1-4\alpha} \right)
+\ldots+ \nonumber\\
&& + \frac{\lambda^{2N}}{2N}   (-1)^N  \left(  C_{irr,N} \eta^{N-1}|t-s| +
 C_{reg,N}(\eta;t-s) \eta^{(1-4\alpha)N} \right) + \ldots \nonumber\\
\EEA
where $|C_{reg,N}(\eta;t-s)|=|C_{reg,N}|t-s|^{4N\alpha}+O(\eta^{2\alpha})|\le C^N$ is the regular part of the $2N$-th connected moment.
Recall also $C_{irr,N}\le C'^N$ for some constant $C'$. Note that the condition
$\alpha>\frac{1}{8}$ implies: $\eta^{N-1}< \eta^{(1-4\alpha)N}$ for all  $N\ge 2$. 

 Let $\lambda_0=\frac{1}{2\sqrt{\max(C,C')}}$,
 then the series converges for $|\lambda|\le \lambda_0\eta^{-\half(1-4\alpha)}$ and 
$|\tilde{\phi}^{(c)}_{s,t}(\eta;\lambda)+ \half C_{irr,1} |t-s|\lambda^2| $ is bounded by a constant
times $\lambda^2 \eta^{1-4\alpha}$ when $\eta\to 0$.
Now $\phi_{s,t}(\eta;\lambda)=\exp \phi^{(c)}_{s,t}(\eta;\lambda)$ (see Lemma \ref{lemma:1:exp}) 
which yields the result. \hfill \eop

\begin{Corollary}[uniform exponential moment]

{\it Fix  $T>0$ and $\alpha\in(\frac{1}{8},\frac{1}{4})$. 
Then there exist constants $\lambda_0,C_0$ such that, for every $0<s,t<T$
 and $0\le \lambda\le \lambda_0 \eta^{-\half(1-4\alpha)}$,
\BEQ \esper \left[\exp \lambda \tilde{\cal A}_{s,t}(\eta) \right] \le C_0 e^{\half C_{irr,1} |t-s|\lambda^2}. \EEQ
}
\label{cor:2:uem}
\end{Corollary}

{\bf Proof.} Straightforward.  \hfill \eop

For instance, this implies in particular by Markov's inequality:
\BEQ \proba[{\cal A}_{s,t}(\eta)\ge A]\le C_0 e^{-\half \frac{A^2}{C_{irr,1}|t-s|}} \EEQ
for every $A\le \lambda_0 C_{irr,1}|t-s| \eta^{-\half(1-4\alpha)}.$


\section{Appendix}


Assume $\eta=0$ to begin with.
Recall from Definition \ref{def:2:K'K+-ab}
that $K^{',\pm}_{[a,b]}$ and $K^{*,\pm}_{[a,b]}$ are integral operators from $L^1((a,b))$ to $Hol(\Pi^{\mp})$ (the space of holomorphic
functions on one of the half-planes)  defined by

\BEQ
 (K^{',\pm}_{[a,b]}f)(z) =\frac{\alpha(1-2\alpha)}{2\cos\pi\alpha} \int_a^b f(u) (\pm \II(z-u))^{2\alpha-2}\ du, \quad z\in \Pi^{\mp}
\EEQ
and
\BEQ
 (K^{*,\pm}_{[a,b]}f)(z)=-\frac{1}{4\cos\pi\alpha}  \int_a^b f(u) (\pm \II(z-u))^{2\alpha}\ du, 
\quad z\in \Pi^{\mp}.\EEQ



 The function $u\mapsto (\pm\II(z-u))^{2\alpha-2}$,
$u\in[a,b]$ depends analytically on $z$ if $z$ belongs to the domain $\C\setminus\{s\pm\II y\ |\ s\in[a,b],y\ge 0\}$. Outside this
domain, the functions $(\pm\II(z-u))^{2\alpha-2}$ or $(\pm\II(z-u))^{2\alpha}$ are multivalued and admit singularities.

One of the first (and easiest) results established in this Appendix (see Lemma \ref{lemma:2:lemma1}) is that $K^{',\pm}_{[a,b]}f$,
$K^{*,\pm}_{[a,b]}f$ may be extended analytically to a neighbourhood of any point $u\in(a,b)$ in a neighbourhood of which $f$ is
analytic. Hence (supposing $f$ is analytic on a neighbourhood of $(a,b)$)  problems of multivaluedness and singularities
are concentrated at the ends $a,b$ of the interval of integration.

Lemmas \ref{lemma:2:lemma1}, \ref{lemma:2:lemma2} give the local behaviour of $K^{',\pm}_{[a,b]}f$ and $K^{*,\pm}_{[a,b]}f$
around $a$ and $b$ under some hypotheses on $f$. 
Then Lemmas \ref{lemma:2:lemma4}, \ref{lemma:2:lemma5}, \ref{lemma:2:lemma5bis} generalize the previous results to
the case where $a$ or $b$ is a varying parameter going to $0$, with $f$ possibly possessing
a non-integrable singularity  at $0$.

The generalization to $\eta>0$ is straightforward in principle, but leads (see subsection 1.3)
to some complications in practice (to be specific, they necessitate the introduction of the $\eta$-exponents, see Definition \ref{def:2:admissible}).

The proofs of the Lemmas in this Appendix depend crucially on the properties of Gauss' hypergeometric function $_2 F_1$ recalled in the Introduction.

Let us start with a technical Lemma.

\begin{Lemma}

{\it
Let, for $2\alpha\in\R\setminus\Z$, $\beta>-1$ and $n=0,1,\ldots$,
\BEQ F_n(\alpha,\beta;t;z):=\int_0^t \frac{(t-u)^n}{n!} u^{\beta} (-\II(z-u))^{2\alpha-2}\ du,\quad
z\in \Pi^+ \EEQ
be the $n$-times iterated integral of the function $u\to u^{\beta} (-\II(z-u))^{2\alpha-2}$.

Then
$F_n(\alpha,\beta;t;z)$ has an analytic extension to $\C\setminus \left( \R_-\cup \{t-\II y\ |\ y\ge 0\} \right)$
given by 

\BEA
&& F_n(\alpha,\beta;t;z)=\nonumber\\
&& \frac{\II}{z} t^{\beta+n+1} \ .\ \left\{
\frac{\Gamma(1+\beta)\Gamma(n-1+2\alpha)}{n! \Gamma(n+2\alpha+\beta)} \II e^{-\II\pi\alpha}
z^{2\alpha-1} \ _2 F_1(2-2\alpha,1+\beta;-n+2-2\alpha;1-t/z) \right. \nonumber\\
&& \left. + \frac{\Gamma(-n+1-2\alpha)}{\Gamma(2-2\alpha)} (-\II(z-t))^{2\alpha-1} 
(1-t/z)^n \ _2 F_1(n+2\alpha+\beta,n+1;n+2\alpha;1-t/z) \right\}. \nonumber\\
\EEA
with restriction to $B(0,t)\setminus(-t,0]=\{z\in\C\ |\ |z|<t\}\setminus(-t,0]$ given by
\BEA
&& F_n(\alpha,\beta;t;z)=\nonumber\\
&& -t^{\beta+n+1} \left[ e^{\II \pi\alpha} \frac{\Gamma(2\alpha+\beta-1)}{\Gamma(2\alpha+\beta+n)} t^{2\alpha-2} \ _2 F_1(2-2\alpha,1-2\alpha-\beta-n;2-2\alpha-\beta;z/t) \right. \nonumber\\
&&\left. +\frac{\Gamma(1+\beta)\Gamma(1-2\alpha-\beta)}{\Gamma(2-2\alpha)n!} e^{-\II\pi(\alpha+\beta+1)}
t^{-1-\beta} z^{2\alpha+\beta-1} \ _2 F_1(1+\beta,-n;2\alpha+\beta;z/t) \right]. \nonumber\\
\EEA
}
\label{lemma:2:lemma0}
\end{Lemma}

{\bf Proof.}

Suppose $v\in(0,1)$. If $0<\Arg(z)<\pi/2$ then $-\pi/2<\Arg(-\II z)<0$ and $0<\Arg(1-\frac{t}{z}v)<\pi$.

If $\pi/2<\Arg(z)<\pi$ then $0<\Arg(-\II z)<\pi/2$ and $0<\Arg(1-\frac{t}{z}v)<\pi/2$.

In both cases (hence for any $z\in\Pi^+$) one has $|\Arg(-\II z)+\Arg(1-\frac{t}{z}v)|<\pi$. Hence

\BEA
F_n(\alpha,\beta;t;z)&=& \frac{t^{\beta+n+1}}{n!} \int_0^1 v^{\beta} (1-v)^n
(-\II z)^{2\alpha-2} (1-\frac{t}{z}v)^{2\alpha-2}\ dv \nonumber\\
&=& t^{\beta+n+1} (-\II z)^{2\alpha-2} \frac{\Gamma(1+\beta)}{\Gamma(n+2+\beta)}
\ _2 F_1(2-2\alpha,1+\beta;n+2+\beta;\frac{t}{z}). \nonumber\\ \label{eq:Fn}
\EEA

By the connection formula (\ref{eq:0:1-z}),
\BEA
&& _2 F_1(2-2\alpha,1+\beta;n+2+\beta;\frac{t}{z})\nonumber\\
&& =\frac{\Gamma(n+2+\beta)\Gamma(n-1+2\alpha)}{\Gamma(n+2\alpha+\beta) n!}
\ _2 F_1(2-2\alpha,1+\beta;-n+2-2\alpha;1-\frac{t}{z}) \nonumber\\
&& +(1-\frac{t}{z})^n (1-\frac{t}{z})^{2\alpha-1} \frac{\Gamma(n+2+\beta)\Gamma(-n+1-2\alpha)}{\Gamma(1+\beta)\Gamma(2-2\alpha)}
\ _2 F_1(n+2\alpha+\beta,n+1;n+2\alpha;1-\frac{t}{z}). \nonumber\\
\EEA

If $0<\Arg(z)<\frac{\pi}{2}$ then $-\frac{\pi}{2}<\Arg(-\II(z-t))<\frac{\pi}{2}$ and $0<\Arg(\frac{\II}{z})<\frac{\pi}{2}.$

If $\frac{\pi}{2}<\Arg z<\pi$ then $0<\Arg(-\II(z-t))<\frac{\pi}{2}$ and $-\frac{\pi}{2}<\Arg(\frac{\II}{z})<0.$

In both cases $|\Arg ( -\II (z-t))+\Arg (\frac{\II}{z})|<\pi$. Hence, if $z\in\Pi^+$,
$$ (1-\frac{t}{z})^{2\alpha-1}=\left(\frac{\II}{z}\right)^{2\alpha-1} (-\II(z-t))^{2\alpha-1}.$$
Whence the first result.

Alternatively, if $z\in B(0,t)\setminus(-t,0]$, then the functions $-e^{-\II\pi\alpha} z^{2\alpha-2}$, resp.
$e^{-\II\pi\gamma} (z/t)^{\gamma}$, extend the function $(-\II z)^{2\alpha-2}$, resp. $(-z/t)^{\gamma}$ defined on $\Pi^+$,
whence the second result by applying the connection formula (\ref{eq:0:1/z}) to (\ref{eq:Fn}) . \hfill \eop

\bigskip

We now look at the case where $[a,b]=[0,t]$ is a {\em fixed} interval and present two Lemmas. It is important to understand that we omit
 the dependence
in $t$ of the results since  $t$ is assumed to be a non-zero constant (one might just as well have assumed that $t=1$, but this
would make the use of the lemmas somewhat awkward). This remark is valid for the whole Appendix (and for the whole article).

\begin{Lemma}

{\it
Let $f\in L^1([0,t],\C)$ and $\phi:z\mapsto \int_0^t (-\II(z-u))^{2\alpha-2} f(u)\ du$ with $2\alpha\in\R\setminus\Z$.
(The result applies in particular to $K^{',\pm}_{[0,t]}(\eta)f$ and $K^{*,\pm}_{[0,t]}(\eta)f$).

\begin{enumerate}
\item Assume $f$ is analytic in a (complex) neighbourhood $\Omega$ of $s\in (0,t)$. Then $\phi$ has
an analytic extension to a complex neighbourhood of $s$.

\item Assume $f$ is analytic in a complex neighbourhood $\Omega$ of $0$. Then $\phi$ may be written
on a small enough neighbourhood of $0$ as
\BEQ \phi(z)=(-\II z)^{2\alpha-1} F(z)+G(z) \EEQ
where both $F$ and $G$ are analytic. The function $F$ has the following expression near $0$:
\BEQ F(z)= \frac{\II}{2\alpha-1} \sum_{n\ge 0} a_n \frac{n!}{(2\alpha)_n} z^n. \EEQ
where $(2\alpha)_n=\frac{\Gamma(2\alpha+n)}{\Gamma(2\alpha)}$ is the Pochhammer symbol.
\item (no analyticity assumption is required here) Assume $|z/t|>C>1$.
 Then $\phi(z)=z^{2\alpha-2} F(t/z)$ for some
analytic function $F$.
\end{enumerate}
}
\label{lemma:2:lemma1}
\end{Lemma}

{\bf Proof.}

\begin{enumerate}
\item
Assume $\Omega\supset B(s,2r)$ for some $r>0$. Then the contour of integration may be deformed into
$\Gamma=[0,s-r]\cup \{s+re^{-\II\phi}, \phi:\pi\to 0\} \cup [s+r,t]$. If $z\in B(s,r/2)$ then
$-\II(z-\bar{w})\not\in ]-\infty,0]$ for any $\bar{w}\in \Gamma$, so $\phi(z)=\int_{\Gamma}
f(\bar{w}) (-\II(z-\bar{w}))^{2\alpha-2} \ d\bar{w}$ is well-defined and analytic on $B(s,r/2)$.

\item
Suppose $f(\bar{w})$ is given by the convergent series $\sum_{n\ge 0} a_n \bar{w}^n$ on $B(0,4r)$. 
A first integration by parts 
\BEA
&& \int_0^r f(u) (-\II (z-u))^{2\alpha-2} \ du \nonumber\\
&&=-\frac{\II}{2\alpha-1} \left[ f(r) (-\II (z-r))^{2\alpha-1} -f(0) (-\II z)^{2\alpha-1}\right] \nonumber\\
&& +\frac{\II}{2\alpha-1} \int_0^r f'(u) (-\II(z-u))^{2\alpha-1}\ du \nonumber\\
\EEA
yields a function $\psi_r(z)=\int_0^r f'(u) (-\II(z-u))^{2\alpha-1}\  du$ which has a continuous
extension to the real axis.  Then successive integrations by parts yield

\BEA
&&\int_0^r u^n (-\II(z-u))^{2\alpha-1}\ du  \nonumber\\
&&=-\frac{\II}{2\alpha} r^n (-\II(z-r))^{2\alpha} +\frac{\II}{2\alpha} \int_0^r nu^{n-1} (-\II(z-u))^{2\alpha}
\ du \nonumber\\
&&=\ldots=-\sum_{m=1}^n \frac{\II^m}{(2\alpha)_m} n(n-1)\ldots (n-m+2) r^{n+1-m} (-\II(z-r))^{2\alpha-1+m} \nonumber\\
&& \quad
+\frac{\II^n}{(2\alpha)_n} n! \int_0^r (-\II(z-u))^{2\alpha-1+n}\ du \nonumber\\
&&=-\sum_{m=0}^n \frac{\II^{m+1}}{(2\alpha)_{m+1}} n(n-1)\ldots (n-m+1) r^{n-m} (-\II(z-r))^{2\alpha+m} \nonumber\\
&& \quad
+\frac{\II^{n+1}}{(2\alpha)_{n+1}} n! (-\II z)^{2\alpha+n} \nonumber\\
\EEA

so
\BEA
 \psi_r(z) &=& (-\II z)^{2\alpha} \left( \sum_{n\ge 0} (n+1)a_{n+1} \ .\ \II^{n+1} \frac{n!}{(2\alpha)_{n+1}}
(-\II z)^n \right)  \nonumber\\
& - &    (-\II(z-r))^{2\alpha}   \sum_{n\ge 0} (n+1)a_{n+1} \ . \nonumber\\
&&  \sum_{m=0}^n \frac{\II^{m+1}}{(2\alpha)_{m+1}}
n(n-1)\ldots (n-m+1)r^{n-m} (-\II(z-r))^m \nonumber\\
\EEA

The first series is easily seen to be convergent for $|z|<4r$ since
 $\frac{\Gamma(n+2\alpha)}{\Gamma(2n)}\sim_{n\to\infty} n^{2\alpha}$. The multivalued function $z\mapsto
(-\II(z-r))^{2\alpha}$ is well-defined on $B(0,r)$. The double series converges in the
supremum norm $||\ ||_{\infty}$ on $B(0,r)$ since $|z-r|<2r$ (hence
$ ||r^{n-m} (-\II(z-r))^{m} ||_{\infty}\le (2r)^n)$  and
$$\sum_{m=0}^n \frac{n(n-1)\ldots (n-m+1)}{(2\alpha)_{m+1}}\le \Gamma(2\alpha)\sum_{m=0}^n \left(\begin{array}{c}
 n\\ m\end{array}\right) =\Gamma(2\alpha) \ 2^n.$$

Finally, $\int_r^t f(u) (-\II(z-u))^{2\alpha-2} \ du$ is analytic on $B(0,r).$

\item
Expand $(-\II(z-u))^{2\alpha-2}$ $(u\le t,|z/t|>C>1)$ into $-e^{-\II\pi\alpha} z^{2\alpha-2} \sum_{k\ge 0}
\frac{(2-2\alpha)_k}{k!} \left(\frac{u}{z}\right)^k$.

\end{enumerate}

 \hfill \eop

\bigskip

This Lemma has the following generalization:

\begin{Lemma}

{\it
Let $f\in L^1([0,t],\C)$ be analytic in a  neighbourhood $\Omega$ of $0$, $2\alpha\in(0,4)\setminus\{1,2,3\}$,
 $\beta>-1$ and
$$\phi : z\mapsto \int_0^t u^{\beta} f(u) (-\II(z-u))^{2\alpha-2} \ du\quad (z\in\Pi^+).$$
(The result applies in particular to $K^{',\pm}_{[0,t]}(\eta) (u\mapsto u^{\beta}f(u))$ and
$K^{*,\pm}_{[0,t]}(\eta) (u\mapsto u^{\beta}f(u))$).

Then $\phi$ may be written on a small neighbourhood of $0$ as $z^{2\alpha+\beta-1} F(z)+G(z)$, where
both $F$ and $G$ are analytic.
}
\label{lemma:2:lemma2}
\end{Lemma}

{\bf Proof.}

Suppose   $f(\bar{w})$ is given by the convergent series $\sum_{n\ge 0} a_n \bar{w}^n$ on $B(s,4r)$.
Here again, $\int_r^t u^{\beta}f(u)(-\II(z-u))^{2\alpha-2}\ du$ is readily shown to be analytic
on $B(0,r)$. A first integration by parts
\BEQ \int_0^r f(u) u^{\beta} (-\II(z-u))^{2\alpha-2}\ du=f(r) F_0(\alpha,\beta;r;z)-\int_0^r
f'(u) F_0(\alpha,\beta;u;z)\ du \EEQ
(see Lemma \ref{lemma:2:lemma0}
 for the definition of the functions $F_n$) yields a function $\psi_r(z)=\int_0^r f'(u) F_0(\alpha,\beta;u;z)\ du$ which has a continuous
extension to $\R_+^*$.

Then successive integrations by parts yield

\BEA
&&\int_0^r u^n F_0(\alpha,\beta;u;z)\ du= r^n F_1(\alpha,\beta;r;z)-\int_0^r nu^{n-1}
F_1(\alpha,\beta;u;z)\ du \nonumber\\
&&=\ldots=\sum_{m=0}^{n-1} (-1)^m r^{n-m} n(n-1)\ldots (n-m+1) F_{1+m}(\alpha,\beta;r;z) \nonumber\\
&& +(-1)^n n!
\int_0^r F_n(\alpha,\beta;u;z)\ du \nonumber\\
&&=\sum_{m=0}^n (-1)^m r^{n-m} n(n-1)\ldots(n-m+1) F_{m+1}(\alpha,\beta;r;z)
\EEA
so $$\psi_r(z)=\sum_{n\ge 0} (n+1)a_{n+1} \sum_{m=0}^n (-1)^m r^{n-m} n(n-1)\ldots(n-m+1) F_{m+1}(\alpha,\beta;
r;z) $$
is the sum of two terms (see last formula in Lemma \ref{lemma:2:lemma0}) for $|z|<r$, $z\not\in\R_-$. 

\begin{itemize}
\item
The first one is
\BEA
\psi_r^{(1)}(z)&:=&-e^{\II\pi\alpha} \Gamma(2\alpha+\beta-1) r^{2\alpha+\beta}
 \sum_{n\ge 0} (n+1) a_{n+1} r^n   
 \sum_{m=0}^n \frac{(-1)^m}{\Gamma(2\alpha+\beta+m+1)} \nonumber\\
&& n(n-1)\ldots
(n-m+1)\    _2 F_1(2-2\alpha,-2\alpha-\beta-m;2-2\alpha-\beta;z/r) \nonumber\\
\EEA

\bigskip

Let us assume $2\alpha\in(0,2)$ (otherwise apply once or twice the formula $_2 F_1(a,b,c;z)=\frac{ab}{c}
\int_0^z \ _2 F_1(a+1,b+1,c+1;u)\ du$ \cite{Abra84} to reduce to the following computations).

Suppose $|w|<1$ and $\beta<0$, then
\BEA
&&  _2 F_1(2-2\alpha,-2\alpha-\beta-m;2-2\alpha-\beta;w) \nonumber\\
&& =\frac{\Gamma(2-2\alpha-\beta)}{\Gamma(2-2\alpha)\Gamma(-\beta)}
\int_0^1 t^{1-2\alpha} (1-t)^{-1-\beta}
(1-tw)^{2\alpha+\beta+m}\ dt. \nonumber\\ 
\EEA
If $2\alpha+\beta+m\ge 0$ (in particular as soon as $m\ge 1$) then $|(1-tw)^{2\alpha+\beta+m}|\le 2^{2\alpha+\beta+m}$, so
$$|_2 F_1(2-2\alpha,-2\alpha-\beta-m;2-2\alpha-\beta;w)|\le C. 2^m$$
for all $|w|<1$. On the other hand, if $2\alpha+\beta<0$, then $|(1-tw)^{2\alpha+\beta}|\le (1-t)^{2\alpha+\beta}$,
so $|_2 F_1(2-2\alpha,-2\alpha-\beta;2-2\alpha-\beta;w)|\le C$.

Now, if $\beta\ge 0$, the preceding arguments must be slightly adapted. Let $p:=\lceil \beta \rceil+1$, then
(by \cite{Abra84}, formula (15.2.4))
\BEA 
&&  _2 F_1(2-2\alpha,-2\alpha-\beta-m;2-2\alpha-\beta;w) 
 =\frac{w^{\alpha+\beta-1}}{(2-2\alpha-\beta)_p}
\left( \frac{d}{dw}\right)^p \nonumber\\
&&  \left[ w^{2-2\alpha-\beta+\lceil \beta\rceil} \ . \
 _2 F_1(2-2\alpha,-2\alpha-\beta-m;2-2\alpha-\beta+\lceil \beta \rceil +1;w) \right] \nonumber\\
\EEA
with
$$_2 F_1(2-2\alpha,-2\alpha-\beta-m;2-2\alpha-\beta+\lceil \beta\rceil+1;w)=
\int_0^1 t^{1-2\alpha}(1-t)^{\lceil \beta\rceil -\beta} (1-tw)^{2\alpha+\beta+m}\ dt.$$
The same kind of estimates also apply to the latter hypergeometric function, together with its derivatives up to order $p$.

All together one has proved in all cases:
\BEQ |_2 F_1(2-2\alpha,-2\alpha-\beta-m;2-2\alpha-\beta;w)|\le C. 2^m\EEQ
for all $|w|<1$. 

\bigskip

Now
\BEQ \sum_{m=0}^n \frac{n(n-1)\ldots(n-m+1)}{\Gamma(2\alpha+\beta+m+1)} 2^m\le \sum_{m=0}^n
\left(\begin{array}{c} n\\ m\end{array}\right) m^{C_1} 2^m \le C_2 n^{C_1}  3^n \EEQ
for some constants $C_1,C_2,C_3$. Hence  the series giving $\psi_r^{(1)}$ converges in the supremum
norm $||\ ||_{\infty}$ on $|z|<r$ to an analytic function $F(z)$;

\item
The second one is
\BEA
&& \psi_r^{(2)}(z):=-\frac{\Gamma(1+\beta)\Gamma(1-2\alpha-\beta)}{\Gamma(2-2\alpha)}
e^{-\II\pi(\alpha+\beta+1)}  \left[ \sum_{n\ge 0} (n+1)a_{n+1} r^{n+1} \right. \nonumber\\
&& \left.  \sum_{m=0}^n \frac{(-1)^m}{(m+1)!}  n(n-1)\ldots (n-m+1) \ _2 F_1(1+\beta,-m-1;2\alpha+\beta;
z/r) \right] z^{2\alpha+\beta-1} \nonumber\\
\EEA
which (using the same type of estimates) converges in the supremum norm to a multivalued
function $z^{2\alpha+\beta-1} G(z)$ ($G$ analytic) on $|z|<r$.

\end{itemize}
\hfill \eop


\bigskip

We now generalize the previous results to the case when $[a,b]=[0,\eps]$ or $[\eps,t]$ where
$\eps$ and $t$ are assumed to be {\em complex}. (By definition, $[z,w]:=\{(1-s)z+sw\ |\ s\in[0,1]\}$
if $z,w\in\C$). We assume $t$ is bounded and bounded away from $0$, i.e.  $0<c<|t|<C$, and
 $|\eps|<|t|$ may be arbitrarily close to $0$ (actually, the following Lemmas are meaningful only when $|\eps|$ is small,
otherwise they are redundant with Lemmas \ref{lemma:2:lemma1} and \ref{lemma:2:lemma2}).

\medskip

Recall $\Omega\subset\C$ is star-shaped with respect to $0$ if $(z\in\Omega\Rightarrow \lambda z\in\Omega \forall
\lambda\in[0,1])$.

\begin{Lemma}

{\it

Let $f$ be analytic on a fixed complex star-shaped  neighbourhood $\Omega_f$ of $0$, 
 $\beta>-1$, $2\alpha\in(0,4)\setminus\{1,2,3\}$, and, for $\eps\in\Omega_f\setminus\R_-$
such that $0<|\eps|<C$,
\BEQ g(\eps;z):=\int_0^{\eps} (-\II(z-u))^{2\alpha-2} u^{\beta} f(u)\ du. \EEQ
  The function $z\mapsto g(\eps;z)$,
 initially
defined as an analytic function on $\Im z>max(0,\Im \eps)$, may be extended into an analytic
function on the cut domain $\Omega:=\Omega_f\setminus (\eps(1-\II\R_+)\cup \R_-)$.
 The behaviour of the function $g$ on $\Omega$ is given
as follows.

\begin{itemize}
\item[(i)] Suppose $z\in\Omega$, $|z/\eps|>C>1$: then 
\BEQ g(\eps,z)=\eps^{\beta+1} z^{2\alpha-2} F(\eps,\eps/z) \EEQ
where $F(\eps,\zeta)$ is holomorphic on the domain $\{(\eps,\zeta)\in\C^2\ |\ \eps\in\Omega_f,
|1/\zeta|>C\}$;
\item[(ii)] Suppose $z\in\Omega$, $|1-z/\eps|<c<1$: then
\BEQ g(\eps,z)=\eps^{2\alpha-1+\beta} \left[  F(\eps,1-\frac{z}{\eps})
 + (-\II(\frac{z}{\eps}-1))^{2\alpha-1}
 G(\eps,1-\frac{z}{\eps}) \right]\EEQ
where $F(\eps,\zeta)$ is holomorphic on the domain $\{(\eps,\zeta)\in\C^2\ |\ \eps,\eps(1-\zeta)\in\Omega_f,
|\zeta|<c\}$;
\item[(iii)] Suppose $z\in\Omega$, $|z/\eps|<c<1$, then
\BEQ g(\eps,z)=z^{2\alpha-1+\beta}
F(\eps,\frac{z}{\eps})+ \eps^{2\alpha-1+\beta} G(\eps,\frac{z}{\eps})   \label{eq:2:0z/eps} \EEQ
where $F(\eps,\zeta),G(\eps,\zeta)$ are holomorphic on the domain $\{(\eps,\zeta)\in\C^2\ |\ \eps,\eps\zeta\in\Omega_f,
|\zeta|<c\}$;
\item[(iv)]
Suppose $z\in\Omega$ is in the {\em cut $\eps$-ring} $\Omega_{\eps,c,c'}:=\{0<c'<\left|\frac{z}{\eps}\right|<C',
\ |1-\frac{z}{\eps}|>c>0 \} $
($0<c'<1<C'$), then
\BEQ g(\eps,z)=\eps^{2\alpha-1+\beta} F(\eps,z/\eps)  \EEQ
where $F(\eps,\zeta)$ is holomorphic on the domain $\{(\eps,\zeta)\in\C^2\ |\ \eps,\eps\zeta\in\Omega_f,
c'<|\zeta|<C',|1-\zeta|>c\}$;

\end{itemize}

}
\label{lemma:2:lemma4}
\end{Lemma}

{\bf Remark.}

We shall need in the sequel to define the above function $g$ on the complex plane cut along
non-intersecting half-lines in a general position. The exponents remain of course the same, but the
determination of the power functions should be chosen in a different way.

{\bf Proof.}

Rewrite $g(\eps;z)$ in the following form:
\BEQ g(\eps;z)=\eps^{2\alpha-1+\beta} \int_0^1 (-\II(\frac{z}{\eps}-v))^{2\alpha-2} v^{\beta} f(\eps v)\ dv,
\label{eq:lemma4} \EEQ
a priori valid for $-\II z\in[A;+\infty)$, $A$ large enough,
which defines it as an analytic function on $\{z\in\C\ | \ \frac{z}{\eps}\not\in[0,1]-\II\R_+\}$.

The previous results show that $g$ may be extended analytically to a cut domain $\Omega\subset\Omega_f$
 excluding
two non-intersecting half-lines ending at $0$ and $\eps$, for instance,
 $\Omega=\Omega_f\setminus (\eps(1-\II\R_+)\cup \R_-)$.

Now apply Lemmas \ref{lemma:2:lemma1}, \ref{lemma:2:lemma2} (more precisely, 
the obvious extension of these Lemmas to the case of a function $f(\eps;v):=f(\eps v)$ depending analytically
on a parameter $\eps$). \hfill \eop

\begin{Lemma}[integration against the infinitesimal kernel on a general interval\\]

{\it
Let, for $\alpha\in(0,\frac{1}{2})$, $-1<\beta<0$ and $p=0,1,\ldots$, and for some fixed function $f$ analytic
on a neighbourhood $\Omega_f$ of the triangle $\cal T$ with vertices $\{0,\eps,t\}$,
\BEQ h(\eps,t;z)=\int_{\eps}^t (-\II(z-u))^{2\alpha-2} u^{\beta-p} f(u)\ du\EEQ
where $\eps,t\in\C$,  $0<|\eps|<|t|$, $c<|t|<C$. Choose three non-intersecting half-lines
$D_1,D_2,D_3$ in the exterior of $\cal T$, with endings at the vertices $0,\eps,t$. The function
$z\to h(\eps,t;z)$, initially defined as an analytic function for $\Im z>\max(\Im \eps,\Im t)$,
may be extended into an analytic function on $\Omega=\C\setminus(D_1\cup D_2\cup D_3)$. The behaviour
of $h$ on $\Omega$ is given as follows, where the functions $F,G,\tilde{G},H,\tilde{H}$ are assumed to be
holomorphic:

\begin{itemize}
\item[(i)]
Suppose $z\in\Omega$, $|z/t|>C>1$. Then
\BEQ h(\eps,t;z)=z^{2\alpha-2} \left(t^{\beta-p+1} \tilde{G}(t,t/z)+
 \eps^{\beta-p+1} H(\eps,\eps/z)  \right).
\label{eq:2:exponent1} \EEQ
\item[(ii)]
Suppose $z\in\Omega$, $|z/\eps|>C>1$, $|z/t|<c<1$: then
\BEQ h(\eps,t;z)= \left( z^{2\alpha-1+(\beta-p)} 
F(z)+G(z) \right) + \eps^{(\beta-p)+1} z^{2\alpha-2} H(\eps,\eps/z). \label{eq:2:exponent2} \EEQ
\item[(iii)]
Suppose $z\in\Omega$,  $|1-z/\eps|<c<1$, $|z/t|<c'<1$: then 
\BEQ h(\eps,t;z)=G(z)+ \eps^{2\alpha-1+(\beta-p)} \left[ \tilde{H}(\eps,1-\frac{z}{\eps})+
 (-\II(\frac{z}{\eps}-1))^{2\alpha-1}  H(\eps,1-\frac{z}{\eps}) \right]. \label{eq:2:exponent3}
\EEQ
\item[(iv)]
Suppose $z\in\Omega$, $|z/\eps|<c<1$: then
\BEQ h(\eps,t;z)=
G(z)+ \eps^{2\alpha-1+(\beta-p)} H(\eps,z/\eps). \label{eq:2:exponent4} \EEQ
\item[(v)]
Assume $z\in\Omega$ is in the cut {\em $\eps$-ring}, i.e. $0<c'<|z/\eps|<C'<|t/\eps|$
$(0<c'<1<C')$, $|1-z/\eps|>c>0$. Then
\BEQ h(\eps,t;z)=\left( z^{2\alpha-1+(\beta-p)} F(z)+G(z)\right)+\eps^{2\alpha-1+(\beta-p)} H(\eps,z/\eps). 
\label{eq:2:exponent5} \EEQ
\end{itemize}

}

\label{lemma:2:lemma5}
\end{Lemma}

{\bf Remarks.}
\begin{enumerate}
\item Domains of holomorphy for the functions $F,G,\ldots$ follow from Lemma \ref{lemma:2:lemma4}.
\item  We do not give the behaviour of $h$ in the   $t$-ring $|\eps/t|<c<|z/t|<C$ because
we shall not need it.
\item The results for cases $(i),\ (ii),\ (iii)$ and $(v)$ follow essentially from splitting the integral into
$\int_0^t \ du -\int_0^{\eps}\ du$ (with some extra care when $p\ge 1$ since the singularity is not integrable at $0$).
The method also applies in case $(iv)$ but yields a spurious extra term of the form $z^{2\alpha-1+\beta-p} \tilde{H}(\eps,
z/\eps)$ which is singular when $z=0$. Yet we use this so-called 'crude' verions of case (iv) at some places in the course
of the proof of Theorem 1.1 because this splitting yields the required $(f_b,f_{\eta})$-splitting. Eq. (\ref{eq:2:exponent4})
-- called the 'refined' version of case (iv) -- follows from a different splitting which avoids integrating around $0$. 
\end{enumerate}

{\bf Proof.}

Let us first prove briefly the 'refined' version of case (iv) as stated in the Lemma, see eq. (\ref{eq:2:exponent4}). 
Write
\BEQ h(\eps,t;z)=C\int_{\eps}^t u^{2\alpha-2+\beta-p} (1-z/u)^{2\alpha-2}\ f(u)\ du\EEQ
and expand $(1-z/u)^{2\alpha-2}=\sum_{k\ge 0} \frac{(2-2\alpha)_k}{k!} \left(\frac{z}{u}\right)^k$, $f(u)=\sum_{n\ge 0}
a_n u^n$ ($|u|\le c$). Splitting $\int_{\eps}^t \ du$ into $\int_{\eps}^{c/2}\ du+\int_{c/2}^t\ du$ and exchanging the
order of summation for the first integral leads to the expression $G(z)+\eps^{2\alpha-1+\beta-p}H(\eps,z/\eps)$, while
the second integral is trivially analytic in $z$. (Easy details are left to the reader). 

Let us now prove the 'crude' version of the Lemma (which does not differ from the 'refined' version, except for
case (iv)). 

Suppose first $p=0$. Then
\BEQ h(\eps,t;z)=\int_0^t (-\II(z-u))^{2\alpha-2} u^{\beta} f(u)\ du-\int_0^{\eps} (-\II(z-u))^{2\alpha-2}
u^{\beta} f(u)\ du.\EEQ
The first integral is estimated in Lemma \ref{lemma:2:lemma2}, and the second one in Lemma \ref{lemma:2:lemma4}. One gets:
\BEQ h(\eps,t;z)=\left( z^{2\alpha-1+\beta} F(z)+G(z)\right) +\eps^{\beta+1}z^{2\alpha-2} H(\eps,\eps/z)\EEQ
if $|z/\eps|>C>1$, $|z/t|<c<1$;
\BEQ h(\eps,t;z)=\left( z^{2\alpha-1+\beta} F(z)+G(z)\right)
  +\eps^{2\alpha-1+\beta} \left( \tilde{H}(\eps,1-z/\eps) + (-\II(z/\eps-1))^{2\alpha-1}
H(\eps,1-z/\eps) \right)\EEQ
if $|1-z/\eps|<c<1$, $|z/t|<c'<1$ (note that the function $F$ may be discarded since
$z^{2\alpha-1+\beta}=\eps^{2\alpha-1+\beta} (1-(1-\frac{z}{\eps}))^{2\alpha-1+\beta}=\eps^{2\alpha-1+\beta}
H(1-\frac{z}{\eps})$); 
\BEQ h(\eps,t;z)=\left(z^{2\alpha-1+\beta}F(z)+G(z)\right) + \left( z^{2\alpha-1+\beta} \tilde{H}(\eps,z/\eps)+
 \eps^{2\alpha-1+\beta} H(\eps,z/\eps) \right)\EEQ
if $ |z/\eps|<c<1$ (once again, one may discard the function $F$); and
\BEQ h(\eps,t;z)=\left(z^{2\alpha-1+\beta} F(z)+G(z)\right)+\eps^{2\alpha-1+\beta} H(\eps,z/\eps) \EEQ
in the cut $\eps$-ring.  On the other hand, if $|z/t|>C>1$, then 
\BEQ h(\eps,t;z)=z^{2\alpha-2} \left( t^{\beta+1} F(t,t/z)+\eps^{\beta+1} G(\eps,\eps/z) \right).\EEQ

\bigskip

Suppose now $p>0$. Let $\bar{B}(0,2\rho|t|)\subset \Omega_f$, $\rho\in(0,1]$ maximal, and
set $t'=\rho t$, $\eps'=\max\{ s\in[0,1]\ |\ s\eps\in \bar{B}(0,2\rho|t|)\}\ .\ \eps$, so that
$\eps',t'$ are proportional to $\eps,t$ and contained in $\bar{B}(0,2\rho|t|)$. Then the integral
$\int_{\eps}^t$ splits into $I_1+I_2$, where
$I_1=\int_{\eps}^{\eps'}+\int_{t'}^t$ and $I_2=\int_{\eps'}^{t'}$. The domain of integration of
 $I_1$ is bounded away from the origin, hence $I_1$ does not 'feel' the singularity at $0$ and
its behaviour may be deduced from Lemmas \ref{lemma:2:lemma1} and \ref{lemma:2:lemma2}. The
multivalued terms in $I_2$ at $\eps',t'$ are compensated by those of $I_1$. So one may just
as well assume that $\eps=\eps'$ and $t=t'$, which we do in the sequel.

Write $f(u)=\sum_{k=0}^{p-1} a_k u^k + u^p \tilde{f}(u)$ with $\tilde{f}$ holomorphic
in a neighbourhood of $0$, so that
\BEQ h(\eps,t;z)=\int_{\eps}^t (-\II(z-u))^{2\alpha-2} u^{\beta} \tilde{f}(u)\ du+
\sum_{k=0}^{p-1} a_k \int_{\eps}^t (-\II(z-u))^{2\alpha-2} u^{\beta-p+k}\ du. \label{eq:2:psum} 
\EEQ
The first integral in the right hand side is estimated as in the case $p=0$. As for the sum,
\BEA
 && \int_{\eps}^t (-\II(z-u))^{2\alpha-2} u^{\beta-p+k}\ du  \nonumber\\
&& = \int_{\eps}^{\infty} (-\II(z-u))^{2\alpha-2}
u^{\beta-p+k}\ du - \int_t^{\infty} (-\II(z-u))^{2\alpha-2} u^{\beta-p+k}\ du \nonumber\\
&& =:   g_{p-k}(\eps)-g_{p-k}(t).
\EEA
(Note that the integrals converge since $2\alpha-2+\beta<0$).

Now
\BEA
 &&  g_{p-k}(s)=  -e^{\II\pi\alpha} s^{2\alpha-1+\beta-p+k} \int_0^1 (1-\frac{z}{s}v)^{2\alpha-2}
v^{-2\alpha-\beta+p-k} \ dv \nonumber\\
&&
 =\frac{ e^{\II\pi\alpha}}{2\alpha-1+\beta-p+k}  s^{2\alpha-1+\beta-p+k} \nonumber\\
&&  _2 F_1(2-2\alpha,-(2\alpha-1+\beta-p+k);1-(2\alpha-1+\beta-p+k);z/s) \nonumber\\
\EEA
If $|1-z/s|<c<1$  then the connection formula (\ref{eq:0:1-z}) yields
\BEA
&&  g_{p-k}(s)=\frac{ e^{\II\pi\alpha}}{2\alpha-1+\beta-p+k} \frac{\Gamma(2-2\alpha-\beta+p-k)\Gamma(2\alpha-1)}{\Gamma(-\beta+p-k)}
  z^{2\alpha-1+\beta-p+k} \nonumber\\
&& -\frac{1}{1-2\alpha} e^{\II\pi\alpha} s^{2\alpha-1+\beta-p+k}
(1-z/s)^{2\alpha-1}\ _2 F_1(-\beta+p-k,1;2\alpha;1-z/s). \nonumber\\ \label{eq:2:g1}
\EEA

If $|z/s|>C>1$, then the connection formula (\ref{eq:0:1/z}) yields 
\BEA
&& g_{p-k}(s)=\frac{ e^{\II\pi\alpha}}{2\alpha-1+\beta-p+k} s^{2\alpha-1+\beta-p+k}  \ .\ \nonumber\\
&& \left\{
\frac{2\alpha-1+\beta-p+k}{1+\beta-p+k} ( -z/s)^{2\alpha-2} \   _2 F_1(2-2\alpha,1+\beta-p+k;
2+\beta-p+k;s/z) \right. \nonumber\\
&& \left.  +\frac{\Gamma(2-2\alpha-\beta+p-k)\Gamma(1+\beta-p+k)}{\Gamma(2-2\alpha)} (-z/s)^{2\alpha-1+\beta-p+k} 
\right\}. \label{eq:2:g2} \EEA

Using the formulas  (\ref{eq:2:g1}),(\ref{eq:2:g2}) gives the result (the most singular terms
are obtained for $k=0$).

\hfill \eop

Exactly the same results hold when one integrates against the kernel  $K^{*,\pm}(\eta)$, but the proof is different.

\begin{Lemma}[integration against the integrated kernel on a general interval]

{\it

Let, for $\alpha\in(0,\frac{1}{2})$, $-1<\beta<0$ and $p=0,1,\ldots$, and for some function $f$ analytic
on an $\eps$-independent neighbourhood of $0$,
\BEQ h(\eps,t;z)=\int_{\eps}^t (-\II(z-u))^{2\alpha} u^{\beta-p} f(u)\ du.\EEQ
Then the  results of Lemma  \ref{lemma:2:lemma5} hold if one replaces $\alpha$ with $\alpha+1$.
}

\label{lemma:2:lemma5bis}
\end{Lemma}

{\bf Proof.}

The proof is the same as for Lemma \ref{lemma:2:lemma5} 
except for the computation of $\int_{\eps}^t (-\II(z-u))^{2\alpha} u^{\beta-p+k}\ du$.
The number  $2\alpha+\beta$ is not necessarily negative, so one cannot integrate to infinity.  Use this time
\BEQ \int_{\eps}^t (-\II(z-u))^{2\alpha} u^{\beta-p+k}\ du=\left( \int_z^t -\int_z^{\eps} \right)
(-\II(z-u))^{2\alpha} u^{\beta-p+k}\ du=: g_{p-k}(t)-g_{p-k}(\eps).\EEQ
Setting $v:=\frac{u-z}{s-z}$, one obtains:
\BEA
&& g_{p-k}(s)=e^{\II\pi\alpha} (s-z)^{2\alpha+1} z^{\beta-p+k} \int_0^1 v^{2\alpha} (1-(1-s/z)v)^{\beta-p+k}\ dv \nonumber\\
&&=e^{\II\pi\alpha} \frac{(s-z)^{2\alpha+1}}{2\alpha+1} z^{\beta-p+k} \ _2 F_1(-\beta+p-k,2\alpha+1;2\alpha+2;1-s/z).
\nonumber\\
\EEA

If $|s/z|<c<1$ then the connection formula (\ref{eq:0:1-z}) entails
\BEA
&&  _2 F_1(-\beta+p-k,2\alpha+1,2\alpha+2;1-s/z)=\frac{\Gamma(2\alpha+2)\Gamma(\beta+1-p+k)}{\Gamma(2\alpha+2+
\beta-p+k)} \left(1-s/z\right)^{-2\alpha-1} \nonumber\\
&&  +\left(\frac{s}{z}\right)^{\beta+1-p+k}
\frac{2\alpha+1}{-\beta-1+p-k} \ _2 F_1(2\alpha+2+\beta-p+k,1;2+\beta-p+k;s/z) \nonumber\\ \EEA

On the other hand, if $|s/z|>C>1$ then the connection formula (\ref{eq:0:1/1-z}) entails
\BEA
&&   _2 F_1(-\beta+p-k,2\alpha+1,2\alpha+2;1-s/z)= \nonumber\\
&&
\frac{\Gamma(2\alpha+2)\Gamma(-2\alpha-1-\beta+p-k)}{\Gamma(-\beta+p-k)} (s/z-1)^{-2\alpha-1}  \nonumber\\
&& + \frac{2\alpha+1}{2\alpha+1+\beta-p+k} (z/s)^{-\beta+p-k}
\ _2 F_1(-\beta+p-k,1;-2\alpha-\beta+p-k;z/s)  \nonumber\\
 \EEA

One may check that these expansions lead to the same leading exponents as in Lemma \ref{lemma:2:lemma5}.
\hfill \eop



\newpage

{\small
\begin{thebibliography}{999}
\bibitem{Abra84} M. Abramowitz, A. Stegun, M. Danos, J. Rafelski. {\it
Handbook of mathematical functions}, Harri Deutsch, Frankfurt (1984).
\bibitem{BreMaj83} P. Breuer, P. Major. {\it Central limit theorems for nonlinear functionals of
Gaussian fields}, J. Multivariate Anal. 13 (3), 425--441 (1983).
\bibitem{Billingsley} P. Billingsley. {\it Convergence of probability measures}, Wiley (1968). 
\bibitem{CheNua05} P. Cheridito, D. Nualart, {\it Stochastic integral of divergence type with respect to 
fractional Brownian motion with Hurst parameter $H\in(0,\half)$},   Ann. Inst. H. Poincar\'e B41(6), 1049 (2005).
\bibitem{CQ02} L. Coutin, Z. Qian, {\it Stochastic analysis, rough path analysis and
fractional Brownian motions}, Probab. Theory Relat. Fields 122, 108-140 (2002).
\bibitem{GraNouRusVal04} M. Gradinaru, I. Nourdin, F. Russo, P. Vallois, {\it $m$-order integrals and generalized
It\^o's formula: the case of a fractional Brownian motion with any Hurst index}, Ann. Inst. H. Poincar\'e B41(4),
781 (2004).
\bibitem{HuNu05} Y. Hu, D. Nualart. {\it Renormalized self-intersection local time for fractional Brownian motion},
Ann. Prob. 33 (3), 948--983 (2005).
\bibitem{LeBellac} M. Le Bellac. {\it From critical phenomena to gauge fields}, InterEditions, Paris (1988).
\bibitem{Lej03} A. Lejay. {\it An introduction to rough paths}, S\'eminaire de Probabilit\'es XXXVII, 1--59, Lecture Notes in Math., 1832 (2003).
\bibitem{LLQ02} M. Ledoux, T. Lyons, Z. Qian, {\it L\'evy area of Wiener processes in
 Banach spaces}, Annals of Probability, vol. 30(2), 546-578 (2002).
\bibitem{Lyo98} T. Lyons, {\it Differential equations driven by rough signals}, Rev. Mat. Ibroamericana 14 (2), 215-310 (1998).
\bibitem{LyoQia02} T. Lyons, Z. Qian. {\it System Control and Rough Paths}, Oxford Mathematical Monographs, Oxford
University Press (2002).
\bibitem{Nou08} I. Nourdin, {\it A change of variable formula for the 2D fractional Brownian motion of
index $1/4$} (preprint).
\bibitem{NouPec08} I. Nourdin, G. Peccati. {\it Stein's method on Wiener chaos} (preprint).
\bibitem{NouPec08bis} I. Nourdin, G. Peccati. {\it Stein's method and exact Berry-Ess\'een asymptotics
for functionals of Gaussian fields} (preprint). 
\bibitem{Nua95} D. Nualart, {\it The Malliavin calculus and related topics}, Probability and its applications, Springer Verlag, New-York (1995).
\bibitem{NuaOrt08} D. Nualart, S. Ortiz-Latorre, {\it Central limit theorems for multiple stochastic
integrals and Malliavin calculs}, Stoch. Proc.  Appl. 118, 614--628 (2008).
\bibitem{NuaPec05} D. Nualart, G. Peccati. {\it Central limit theorems for sequences of multiple
stochastic integrals}, Ann. Prob. 33(1), 177-193 (2005).
\bibitem{PipTaq00} V. Pipiras, M. Taqqu, {\it Integration questions related to fractional Brownian motion}, Probab. Theory Related Fields 118 (2), 251 (2002).
\bibitem{RusVal93} F. Russo, P. Vallois. {\it Forward, backward and symmetric stochastic integration}, Prob. Th. Relat. Fields 97, 403-421 (1993).
\bibitem{RusVal00} F. Russo, P. Vallois. {\it Stochastic calculus with respect to continuous finite quadratic variation processes}, Stochastics and stochastics reports 70, 1-40 (2000). 
\bibitem{Unt08} J. Unterberger. {\it Stochastic calculus for fractional Brownian motion with Hurst exponent $H>1/4$:
a rough path method by analytic extension.} To appear in Ann. Prob.
\end{thebibliography}}
\end{document}